\documentclass[a4,10pt]{article}
\usepackage{bbm}
\usepackage{amsmath}
\usepackage{amssymb}
\usepackage{amsfonts}
\usepackage{cite}
\usepackage[dvips]{graphicx}
\usepackage{color}

\numberwithin{equation}{section}

\newtheorem{theorem}{Theorem}[section]

\newtheorem{corollary}[theorem]{Corollary}
\newtheorem{defn}[theorem]{Definition}

\newtheorem{lemma}[theorem]{Lemma}
\newtheorem{prop}[theorem]{Proposition}

\newtheorem{remark}[theorem]{Remark}

\def \df{\square}
\def \mc{\mathcal}
\def \inv{^{-1}}

\def \v{\vskip 0.1in}
\def \n{\noindent}

\def \real{\mathbb{R}}
\def \cplane{\mathbb{C}}
\def \mff{\mathsf}
\def \indexn{^{(k)}}
\def \integer{\mathbb{Z}}
\def \detf{\det(f_{s\bar{t}})}
\def \al{\frac{1}{3}}
\def \six{\frac{1}{6}}
\def \half{\frac{1}{2}}

\def \bb {\mathbb}
\def \b{\bar}
\def \p{\partial}

\def \halfplane{\mathsf{h}^\ast}
\def \indexm{_k}
\def \fkz{\mathfrak z}

 \def\CHART{\mathsf{U}}
\def\t{\mathfrak{t}}
\begin{document}

\title{Affine techniques on extremal metrics on toric surfaces}

\author{Bohui Chen, An-Min Li\footnote{Corresponding author.}, Li Sheng}

\maketitle
{\abstract This paper consists of a package of real and complex affine techniques for
the Abreu equation on toric surfaces/manifolds. In particular, as an application of this package we give an interior estimate for   the  Ricci
tensor of  toric surfaces.
\endabstract}

\tableofcontents

\v Studying the extremal metrics on compact toric manifolds has
been an interesting project in  complex geometry over the past
decade. The project was mainly
formulated by Donaldson in \cite{D1}. Since then, some important
progress has been made. In particular, Donaldson   \cite{D2}
completes the study of the metrics of  
 constant scalar curvature
on toric surfaces. However, finding the extremal metrics on toric
surfaces is still open. In this paper and its  sequal \cite{CLS},
we will study the metrics on toric surfaces with prescribed ``scalar
curvature" function $K$ on the Delzant polytopes. In \cite{CLS}, we
prove that, for any smooth function $K$ on the Delzant polytope with
a mild assumption (cf.\;Definition \ref{definition_1}) and satisfying
 a certain stability condition, there exists a K\"{a}hler metric on the
toric surface whose ``scalar curvature" is $K$. In particular, this
includes most  extremal metrics.

The purpose of this paper is to develop a package of techniques  for the Abreu equation. Since these techniques are motivated by the study of affine geometry by the second author and his collaborators, we call them the affine techniques. One of the highlights in this paper is that we also extend the affine techniques to the complex case (cf.\;\S\ref{sect_5} and \S\ref{sect_6}).  This enables us to obtain an estimate for  the Ricci tensor (cf.\;Theorem \ref{theorem_7.0.1}), which is one of the most important steps for \cite{CLS}.
So far, people  lack tools to achieve estimates for Ricci tensors in general. We hope
this paper may shed  light on this important issue in K\"abler geometry.

For  toric manifolds, the equation  of  extremal metrics can be reduced to a real 4th-order partial differential equation on the Delzant polytope in  a real space. This equation is known as the Abreu equation. The second author with his collaborators  has developed a framework to study  certain  types of 4-th order PDEs (cf.\cite{L-J,L-J-1,L-J-2,L-J-3,L-J-4,L-X1,L-X2,L-X-S-J,CLS1,CLS2,L-S})
which  include  the Abreu equation. This is explained in  \S\ref{sect_2}-\S\ref{sect_4}. The whole package of these techniques includes the differential inequality  for the norm of the Tchebychev vector field  $\Phi$  (\S\ref{sect_2.3}), the
convergence theorems (cf.\;Theorem \ref{theorem_2.4.2},  Theorem
\ref{theorem_3.2.7}  and \S \ref{sect_3.2}), the Bernstein properties
(\S\ref{sect_3.3}) and the affine blow-up analysis (\S\ref{sect_4}). This package is very powerful  not only for the interior estimates  inside
 polytopes but also for the behaviour
of the estimates approaching boundaries. The latter
one is very crucial for the boundary estimates.  For example, the estimates in \S\ref{sect_4} are  of this type.

The challenging part of solving the extremal metric problem is to study the
boundary estimates of the Abreu equation. Since the boundary of polytopes  can be thought as the
interior of the complex manifold, it  is then natural  to ask whether the real affine techniques can be  extended
 to the complex case or not. In this paper,
we make an attempt on this direction. These techniques  we develop are
enough for us to study the Abreu equations and will be used in
\cite{CLS}   to completely solve the extremal metric problem for toric surfaces. We call the techniques we develop through
 \S\ref{sect_5}-\S\ref{sect_7} the {\em complex affine techniques}.  These techniques also   include
 the differential inequalities for affine-type
invariants in \S\ref{sect_5}.
As applications of these differential inequalities, we
 establish the interior estimates in \S\ref{sect_6} and prove  a convergence theorem (Theorem
\ref{theorem_7.4.1}) in  \S\ref{sect_7}.

Finally, with the aid of blow-up analysis, we are able
to establish the estimate of the norm $\mc K$ of the Ricci tensor (\eqref{eqn_1.1c}) in a neighborhood of points on edges of
near the boundary (cf.\;Theorem \ref{theorem_7.0.1}). This is the main theorem in this paper.

 Acknowledge. We would like to thank Xiuxiong Chen, Yongbin
Ruan and Gang Tian for their  interest in the project and constant
supports. We would also like to thank Bo Guan  and Qing Han  for
 valuable discussions.  The first author is partially supported
by NSFC 10631050 and NKBRPC (2006CB805905), and the second author is
partially supported by NSFC 10631050, NKBRPC (2006CB805905) and
RFDP (20060610004).


\section{K\"{a}hler geometry on toric surfaces}\label{sect_1}

In this section, we review  the K\"{a}hler geometry of toric surface and introduce the notations to be used in this paper. We assume that the readers are familiar with some basic knowledge of toric varieties.

A toric manifold is a 2n-dimensional K\"{a}hler manifold
$(M,\omega)$ that admits an $n$-torus (denoted $\bb T^n$) Hamitonian action.
Let $\tau: M\to \mathfrak{t}^\ast$ be the moment map, where
$\t^\ast$, identified with $  \real^n$,   is the  dual of the  $\mathfrak{t},$ which is the Lie algebra of $\bb T^n$.  The image $\bar\Delta=\tau(M)$ is known to be
a convex polytope  \cite{De}.
In the literature, people use $\Delta$ for the image of the moment map.
However,  it is  more convenient in this paper to always assume that {\em $\Delta$ is an open polytope.}   Note that $\Delta$ determines a fan $\Sigma$ in $\t$. The converse is not true: $\Sigma$ determines $\Delta$ only up to  a certain similarity. $M$ can be reconstructed from either $\Delta$ or $\Sigma$ (cf.\cite{Fu} and \cite{Gu}). Moreover, the class of $\omega$ can be read from $\Delta$. Hence, the uncertainty of $\Delta$ reflects the non-uniqueness of K\"{a}hler classes. Two different constructions are related via Legendre transformations.
The K\"{a}hler geometry appears naturally when considering the
transformation between two different constructions. This was
explored by Guillemin   in \cite{Gu}. We will summarize these facts in
this section. For simplicity, we only consider the toric surfaces,
i.e, $n=2$.

\subsection{Toric surfaces and coordinate charts}\label{sect_1.1}

Let $\Sigma$ and $\Delta$ be a pair  consisting of a fan and a polytope  for a toric surface $M$. For simplicity, we focus on compact toric
surfaces. Then $\Delta$ is a Delzant polytope in $\t^\ast$ and
its closure is compact.

We use the notations in \S 2.5   of \cite{Fu}  to describe the fan. Let
$\Sigma$ be a fan given by a sequence of lattice points
$$\{
\nu_0,\nu_1,\ldots, \nu_{d-1},\nu_d=\nu_0 \}$$ in  the counterclockwise order, in
$N=\integer^2\subset \t$, such that successive pairs generate
$N$.
Suppose that the vertices and  edges of $\Delta$ are denoted by
$$
\{\vartheta_0,\ldots,\vartheta_d=\vartheta_0\}, \;\;\;
\{\ell_0,\ell_1,\ldots, \ell_{d-1},\ell_d=\ell_0\}.
$$
Here $\vartheta_i=\ell_i\cap \ell_{i+1}$.

By saying that $\Sigma$ is dual to $\Delta$ we mean that $\nu_i$ is
the inward pointing normal vector to $\ell_i$ of $\Delta$. Hence,
$\Sigma$ is determined by $\Delta$. Suppose that the equation for
$\ell_i$ is
\begin{equation}\label{eqn_1.1}
l_i(\xi):=\langle\xi,\nu_i\rangle- \lambda_i=0.
\end{equation}
 Then  we have
$$
\Delta=\{\xi|l_i(\xi)> 0,\;\;\; 0\leq
i\leq d-1\}
$$

There are three types of cones in $\Sigma$: a 0-dimensional cone
$\{0\}$  denoted by $C_\Delta$;  1-dimensional cones
generated by $\nu_i$ and denoted by  $C_{\ell_i}$;
2-dimensional cones generated by $\{\nu_{i},\nu_{i+1}\}$ and
denoted by $C_{\vartheta_i}$.
It is known that for each cone of $\Sigma$, one can associate to it a
complex coordinate chart  of $M$ (cf.\;\S1.3 and \S 1.4 in
\cite{Fu}).
 Let $\CHART_\Delta,
\CHART_{\ell_i}$ and $\CHART_{\vartheta_i}$ be the coordinate
charts. Then
$$
\CHART_{\Delta}\cong (\cplane^\ast)^2;\;\;\;
\CHART_{\ell_i}\cong \cplane\times\cplane^\ast ;\;\;\;
\CHART_{\vartheta_i}\cong \cplane^2.
$$
In particular, in each $\CHART_{\ell_i}$ there is a divisor
$\{0\}\times \cplane^\ast$. Its closure is a divisor in $M$, we
denote it by $Z_{\ell_i}$.

\begin{remark}\label{remark_1.1.1}
$\cplane^\ast$ is called a complex torus and denoted by $\bb T^c$. Let $z$ be its natural
coordinate.

In this paper, we introduce another complex coordinate by
considering the following identification
\begin{equation}\label{eqn_1.1a}
 \bb T^c\to \real\times 2\sqrt{-1}\bb T;\;\;\;
 w=\log z^2.
\end{equation} We call $w=x+ \sqrt{-1}y$  the
 log-affine complex coordinate (or  log-affine coordinate) of $\cplane^\ast$.

When $n=2$, we have
$$
(\cplane^\ast)^2\cong \t\times 2\sqrt{-1}\bb T^2.
$$
Then $(z_1,z_2)$ on the left hand side is  the usual complex
coordinate; while $(w_1,w_2)$ on the right hand side is  the
 log-affine coordinate.
Write $w_i=x_i+ \sqrt{-1}y_i$, $y_{i}\in [0,4\pi]$.
Then $(x_1,x_2)$ is the coordinate of $\t$.
\end{remark}

We make the following convention.
\begin{remark}\label{remark_1.1.2}
On different  types of coordinate chart  we use different
coordinate systems as follows:
\begin{itemize}
\item on $\CHART_\vartheta\cong \cplane^2$, we use the coordinate
$(z_1,z_2)$; \item on $\CHART_\ell\cong
\cplane\times\cplane^\ast$, we use the coordinate $(z_1,w_2)$;
\item on $\CHART_\Delta\cong (\cplane^\ast)^2$, we use the
coordinate $(w_1,w_2)$, \end{itemize}
 where $z_i=e^{\frac{w_i}{2}},i=1,2$.
\end{remark}

\begin{remark}\label{remark_1.1.3}
Since we study the $\bb T^2$-invariant geometry on $M$, it is
useful to  specify a representative point of each $\bb T^2$-orbit.
Hence for $(\cplane^\ast)^2$, we let the points on $\t\times
2\sqrt{-1}\{1\}$  be the representative points.
\end{remark}

\subsection{K\"{a}hler geometry on toric surfaces}\label{sect_1.2}
  Guillemin     in \cite{Gu} constructed  a natural $T^2$-invariant K\"{a}hler  form
$\omega_o$ on  $M$ from the polytope $\Delta$.   We take this form  as a reference point in  the class
 $[\omega_o]$ and call  the associated K\"{a}hler metric  the Guillemin metric.

For each $T^2$-invariant K\"{a}hler form $\omega\in [\omega_o]$, on each coordinate chart  associated to a
cone of the fan $\Sigma$, there
is a    K\"{a}hler  potential function (unique up to
 linear functions).  Write the collection of the potential functions as
$$
\mathsf f=\{f_\bullet\}:=\{f_\Delta, f_{\ell_0},\ldots, f_{\ell_{d-1}}, f_{\vartheta_0},\ldots,
f_{\vartheta_{d-1}} \}
$$ as  the   K\"{a}hler potential function with respect to the coordinate charts
$$\{
\CHART_\Delta,
\CHART_{\ell_i}, \CHART_{\vartheta_i} | i=1, \cdots, d-1\}.$$
We write $\omega=\omega_{\mathsf f}$    to indicate the associated potential functions.

 Let
$ \mathsf{g}=\{g_\bullet \}$
be the collection of potential functions for $\omega_o$.  We can  realize $\mathsf f$ by the following construction.
Let $C^\infty_{\bb T^2}(M)$ be the smooth $\bb T^2$-invariant
functions of $M$.  Set
$$
 C^\infty_+(M) =
\{\phi\in C^\infty_{\bb T^2}(M)|\omega_{g_\Delta+\phi}>0\}.$$
Then for $\omega\in [\omega_o]$ there exists a function $\phi$ such that
$$
\mathsf f=\mathsf g+\phi:=(g_\bullet+\phi)
$$
and $\omega=\omega_\mathsf f$. Set
$$
\mc C^\infty(M,\omega_o)=\{\mathsf f|\mathsf f=\mathsf g+\phi, \phi\in C^\infty_+(M)\}.
$$

\begin{remark}\label{remark_1.2.1}
 Suppose that  $f_\bullet=g_\bullet+\phi.$
Consider the matrix
$$
\mathfrak{M}_ {f}=(\sum_k g^{i\bar
k} f_{j\bar k}).
$$
Though this is not a globally well defined matrix on $M$, its
eigenvalues are defined globally. Set $ \nu_f$ to be the set of eigenvalues and $
H_ {f}=\det\mathfrak{M}_ {f}\inv. $ These are global
functions on $M$.
\end{remark}
   Under   a coordinate chart with potential function $f$, the
Christoffel symbols, the curvature tensors, the Ricci curvature
and the scalar curvature of K\"{a}hler metric
$\omega_{f}$ are given by
$$ {\Gamma}^{k}_{ij}=\sum_{l=1}^{n} f^{k \b l}\frac{\p f_{i \b l}}{\p z_j}, \;\;\;
 {\Gamma}^{\b k}_{\b i\b j}=\sum_{l=1}^{n} f^{\b k  l}\frac{\p
f_{\b i l}}{\p z_{\b j}} ,  $$
$$ {R}_{i\b jk\b l}=- \frac{\p ^2 f_{i\b j}}{\p z_k \p z_{\b l}}+\sum f^{p\b q}\frac{\p f_{i \b q}}{\p z_k}
\frac{\p f_{p \b j}}{\p z_{\b l}},$$
$$
R_{i\b j}= - \frac{\partial^{2}}{\partial z_{i}\partial \bar
z_{j}} \left(\log \det\left(f_{k\bar l}\right)\right), \;\;\;\mc
S=\sum f^{i\bar j}R_{i\b j},
$$
respectively. When we use the log-affine coordinates
on $\t$,   the Ricci curvature
and the scalar curvature can be written as
$$
R_{i\b j}= - \frac{\partial^{2}}{\partial x_{i}\partial x_{j}}
\left(\log \det\left(f_{kl}\right)\right), \;\;\;\mc S=-\sum
f^{ij}\frac{\partial^{2}}{\partial x_{i}\partial x_{j}} \left(\log
\det\left(f_{kl}\right)\right).
$$
 We  treat
$\mc S$ as an operator for $f$ and denote it by
$\mc S(f)$.

 Define
\begin{equation}\label{eqn_1.1c}
\mc K\;=\;\|Ric\|_f +\|\nabla Ric\|_f
^{\frac{2}{3}}+\|\nabla^2 Ric\|_f
^{\frac{1}{2}}.
\end{equation}
  We also   denote by $\dot \Gamma_{ij}^k$, $\dot
R^m_{ki\bar{l}}$  and $\dot R_{i\b j}$ the connections, the
curvatures and the Ricci curvature of the metric $\omega_{o} $
respectively.

When  focusing on $\CHART_\Delta$ and using the log-affine
coordinate (cf.\;Remark \ref{remark_1.1.1}), we have $f(x)=g(x)+\phi(x)$. We  remark that when restricting on $\real^2\cong \real^2\times2\sqrt{-1}\{1\}$,
the Riemannian metric induced from $\omega_\mathsf{f}$ is the Calabi metric $G_f$ (cf.\;\S \ref{sect_2.1}).

  Fix a large constant $K_o>0$. We set
$$
\mc C^\infty(M,\omega_o;K_o)
=\{\mathsf{f}\in \mc C^\infty(M,\omega_o)|
|\mc S(f)|\leq K_o\}.
$$
In this paper, we  mainly study the apriori estimates for the functions in this class.

\subsection{The Legendre transformation, moment maps and potential functions}
\label{sect_1.3}
 Let $f$ be a (smooth) strictly convex function on
$\t$.
 The gradient of $f$
defines a (normal) map $\nabla^f$ from $\t$ to  $\t^\ast$:
$$
\xi=(\xi_1,\xi_2)=\nabla^f(x) =\left(\frac{\partial f}{\partial
x_1},\frac{\partial f}{\partial x_2}\right).
$$
The function $u$ on $\t^\ast$
$$
u(\xi)=x\cdot\xi - f(x).
$$
is called the Legendre transformation of $f$. We write $u=L(f)$. Note that $f=L(u)$.

Now we restrict  to $\CHART_\Delta.$ When we use the coordinate $(z_1,z_2),$
the moment map with respect to $\omega_ {f}$
 is given by
\begin{equation}\label{eqn_1.1b}
\tau_{{\mathsf{f}}}:\CHART_\Delta\xrightarrow{(\log |z_1^2|,\log|z_2^2|)}
\t\xrightarrow{\nabla^f} \Delta.
\end{equation}
Note that the first map is induced from \eqref{eqn_1.1a}. Let $u=L(f)$. It
is known that $u$ must satisfy certain behavior near boundary of
$\Delta$ by the following theorem.
\begin{theorem}[Guillemin]\label{theorem_1.3.1}
 Let $v=L(g)$, where $g=g_\Delta$ is the potential
function of the Guillemin metric. Then
$v(\xi)=\sum_{i} l_i\log l_i$, where $l_i$ is defined in \eqref{eqn_1.1}.
\end{theorem}

	For $u=L(f)$, we have $u=v+\psi$, where $\psi\in C^\infty(\bar\Delta)$. Motivated by this,
	  we set
$$
\mc C^\infty(\Delta,v)=\{u| u=v+\psi \mbox{ is strictly convex, }
\psi\in C^\infty(\bar \Delta)\},
$$
 with $v$ as a
reference point.   Note that  this space only depends on $\Delta$.

We summarize the fact we just  presented: let $\mathsf{f}\in \mc
C^\infty(M,\omega_o)$, then the moment map $\tau_{{\mathsf{f}}}$
is given by $f=f_\Delta$ via the diagram \eqref{eqn_1.1b} and $u=L(f)\in \mc
C^\infty(\Delta,v)$. Conversely, $f_\bullet$ can be constructed from $u$ as well.

Given a function $u\in \mc C^\infty(\Delta,v)$, we can get an
$\mathsf{f}\in \mc C^\infty(M,\omega_o)$ as  follows.
\begin{itemize}
\item On $\CHART_\Delta$,
 $f_\Delta=L(u)$;
\item    on $\CHART_{\vartheta_i}$, $f_{\vartheta_i}$ is constructed in the
following steps: (i), since $\vartheta_i=\ell_i\cap \ell_{i+1}$, let $B\in SL(2,\mathbb Z)$ be the transformation of $\mathfrak t^\ast$ such that
$$
B(\nu_i)=(1,0), \;\;\;B(\nu_{i+1})=(0,1).
$$
  Meanwhile, $u$ is transformed to a function in the following format
$$u'=\xi_1\log\xi_1+\xi_2\log\xi_2+\psi';$$
(ii), $f'=L(u')$ defines  a function on $\t$ and therefore is a
function on $(\cplane^\ast)^2\subset
 \CHART_{\vartheta_i}$ in terms of log-affine coordinate;   \\
 (iii),   it is known that
 $f'$ can be extended over $\CHART_{\vartheta_i}$ and we set $f_{\vartheta_i}$
to be this function; \item on $\CHART_\ell$, the construction of
$f_\ell$ is similar to $f_\vartheta$.  The reader may refer to
\S\ref{sect_1.5} for the construction.
\end{itemize}

\subsection{The Abreu equation on $\Delta$}\label{sect_1.4}
We can transform  the scalar curvature operator $\mc S(f)$ to  an
operator $\mc S(u)$ of $u$ on $\Delta$, where $u=L(f)$.   Then
$$
\mc S(u)=\mc S(f)\circ\nabla^u.
$$
  The operator $\mc S(u)$  is known
to be
$$
\mc S(u)=-\sum U^{ij}w_{ij}
$$
where $(U^{ij})$ is the cofactor matrix of the Hessian matrix
$(u_{ij})$, $w=(\det(u_{ij}))\inv$.
It is well known that $\omega_{\mathsf{f}}$ gives an extremal
metric if and only if $\mc S\circ\nabla^ u$ is a linear function
of $\Delta$. Let $K$ be a smooth function on $\bar\Delta$, the Abreu
equation is
\begin{equation}\label{eqn_1.2}
\mc S(u)=K.
\end{equation}
We set $\mc C^\infty(\Delta, v;K_o)$
to be the functions {\bf $u\in\mc C^\infty(\Delta,v)$} with $|\mc S(u)|\leq K_o$.

\begin{defn}\label{definition_1}
Let $K$ be a smooth function on $\bar\Delta$. It is called {\em edge-nonvanishing}
if it does not vanish on any edge of $\Delta$.
\end{defn}
In our papers, we will always assume that $K$ is edge-nonvanishing.

\subsection{A special case: $\cplane\times \cplane^\ast$}\label{sect_1.5}

Let $\halfplane\subset \t^\ast$ be the half plane given by
$\xi_1\geq 0$. The boundary is the $\xi_2$-axis and we  denote it by
$\t_2^\ast$.
 The corresponding fan consists of only one lattice $ \nu=(1,0)$.
 The coordinate chart is
$\CHART_{\halfplane}= \cplane\times\cplane^\ast$. Let
$Z=Z_{\t_2^\ast}=\{0\} \times \cplane^\ast$ be its divisor.

Let $v_{\halfplane}= \xi_1\log\xi_1+\xi_2^2$. Set
$$
\mc C^\infty(\halfplane,v_{\halfplane})=\{u|u=v_{\halfplane}+\psi
\mbox{ is strictly convex, } \psi\in C^\infty(\halfplane)\}
$$
and $\mc C^\infty(\halfplane, v_{\halfplane};K_o)$
be the functions whose $\mc S$ is less than $K_o$.

Take  a function $u\in \mc C^\infty(\halfplane,
v_{\halfplane})$. Then $f=L(u)$ is a function on $\t$.
Hence it defines a function on the $\cplane^\ast\times\cplane^\ast
\subset \CHART_{\halfplane}$ in terms of log-affine coordinates
$(w_1,w_2)$. Then the function $  f_\mathsf{h}(z_1,w_2):=f(\log |z_1^2|,
Re(w_2)) $ extends smoothly over $Z$,  and hence is defined on
$\CHART_{\halfplane}$. We conclude that for any $u\in \mc
C^\infty(\halfplane,v_{\halfplane})$ it yields a potential
function $ f_\mathsf{h}$ on $\CHART_{\halfplane}$.

 When we choose the coordinate $(z_1,w_2),$
the moment map with respect to $\omega_{f}$
 is given by
\begin{equation}\label{eqnc_1.1b}
\tau_{{\mathsf{f}}}:\CHART_\halfplane\xrightarrow{(\log |z_1^2|,Re(w_2))}
\t\xrightarrow{\nabla^f}\halfplane.
\end{equation}

Using $v_{\halfplane}$ and the above argument, we can get a potential function denoted by $g_\mathsf{h}$ on $\CHART_\halfplane$.


\section{Calabi geometry}\label{sect_2}

\subsection{Calabi metrics and basic affine invariants}\label{sect_2.1}
Let $f(x)$ be a smooth, strictly convex function defined on a
convex domain  $\Omega\subset\real^n\cong \t$.  As $f$ is strictly
convex,
$$
G:=G_f=\sum_{i,j} f_{ij}d x_i dx_j
$$
defines a Riemannian metric on $\Omega$. We call it the {\em Calabi} metric. We recall some fundamental facts on the Riemannian manifold $(\Omega, G)$ (cf.\cite{P}). The Levi-Civita connection is given by
$$\Gamma_{ij}^k=\frac{1}{2}\sum {f^{kl}}{f_{ijl}}.$$
The Fubini-Pick tensor is
$$A_{ijk}= -\frac{1}{2}f_{ijk}.$$
Then the  curvature tensor and the Ricci tensor are
\begin{eqnarray*}
R_{ijkl} &=& \sum f^{mh}(A_{jkm}A_{hil}-A_{ikm}A_{hjl})
\\
R_{ik}&=&\sum f^{mh}f^{jl}(A_{jkm}A_{hil}-A_{ikm}A_{hjl}).
\end{eqnarray*}
Let $u$ be
 the Legendre transformation of $f$ and $
 \Omega^\ast=\nabla^f(\Omega)\subset \t^\ast$.
Then it is known that $ \nabla^f: (\Omega, G_f)\to (\Omega^\ast,
G_u) $ is isometric.

Let $ \rho=\left[\det(f_{ij})\right]^{-\frac{1}{n+2}}, $ we
introduce the following affine invariants:
\begin{equation}\label{eqn_2.1}
\Phi=\frac{\|\nabla\rho\|^2_{G_f}}{\rho^2}=\frac{1}{(n+2)^2}\|\nabla \log\det(u_{ij})\|_{G_{u}}^2 \end{equation}
\begin{equation}\label{eqn_2.2}
4n(n-1)J=\sum f^{il}f^{jm}f^{kn}f_{ijk}f_{lmn}= \sum
u^{il}u^{jm}u^{kn}u_{ijk}u_{lmn}.
\end{equation}
  Note that $\Phi$ is called the norm of {\em the Tchebychev vector field} and
$J$ is called {\em the Pick invariant}. In this paper, we usually give estimates for $J$ and $\Phi$ simultaneously. For convenience,   Put
\begin{equation}\label{eqn_2.3}
\Theta = J + \Phi.
\end{equation}

\subsection{Affine transformation rules}\label{sect_2.2}
We study how the affine transformation affects the Calabi geometry.
\begin{defn}\label{defn_2.2.1}
By an affine transformation, we mean a transformation of  the following form:
$$
\hat A: \t^\ast\times\real\to \t^\ast\times\real; \;\;\; \hat
A(\xi, \eta)=(A\xi, \lambda \eta),
$$
where $A$ is an affine transformation on $\t^\ast$. We may write
$\hat A=(A,\lambda)$. If $\lambda=1$
we call $\hat A$ the base-affine transformation.

Let $u$ be a function on $\t^\ast$. Then
$\hat A$ induces an affine
transformation on $u$:
$$
u^\ast(\xi)= \lambda u(A\inv\xi).
$$
We denote $u^\ast$ by $\hat A(u)$.
\end{defn}

Then we have the following lemma of the {\em affine transformation
rule} for the affine invariants.
\begin{lemma}\label{lemma_2.2.1}
Let $u^\ast=\hat A(u)$ be as above, then
\begin{enumerate}
\item $\det(u^\ast_{ij})(\xi)=\lambda^n\det(A)^{-2}\det(u_{ij})(A\inv
\xi)$. \item $G_{u^\ast}(\xi)= \lambda G_u(A\inv \xi)$; \item
$\Theta_{u^\ast}(\xi)=\lambda\inv \Theta_u(A\inv \xi)$; \item $\mc
S(u^\ast)(\xi)= \lambda\inv\mc S(u)(A\inv \xi)$.
\end{enumerate}
\end{lemma}
These can be easily verified. We skip the proof. As a corollary,   we have the following invariants with respect to affine transformations
\begin{lemma}\label{lemma_2.2.2}
$G$ and $\Theta$ are invariant with respect to the base-affine
transformation. $\Theta\cdot G$ and $\mc S\cdot G$ are invariant
with respect to affine transformations.
\end{lemma}

\subsection{The differential inequality of $\Phi$}\label{sect_2.3}

 We calculate $\Delta\Phi$ on the Riemannian manifold $(\Omega, G_f)$
  and derive a
differential inequality, where $\Delta$ is the Laplacian operator with respect to the Calabi metric.. We only need  this inequality for the case
$\mc S(f)=0$ in this paper.  This is already done in \cite{L-J}. The
following results can be found in  \cite{L-J} (See also
\cite{L-J-4,L-S}).
\begin{lemma}\label{lem_5.1}
The following two equations are equivalent
\begin{eqnarray}\label{eqn_2.4}
&& \mc S(f)=0,\\
&& \Delta \rho= \frac{n+4}{2}\frac{\|\nabla\rho\|_G^2}{\rho}.
\end{eqnarray}
\end{lemma}
\begin{theorem}\label{theorem_2.3.1}
Let $f$ be a smooth strictly convex function on $\Omega$ with $\mc
S(f)=0$, then
\begin{eqnarray*} \Delta \Phi &\geq&
\frac{n}{2(n-1)} \frac{\|\nabla
\Phi\|_G^2}{\Phi}+\frac{n^{2}-4}{n-1}\langle \nabla \Phi,
\nabla\log
\rho \rangle_{G} \\
&&+\frac{(n+2)^{2}}{2}\left(\frac{1}{n-1}-\frac{n-1}{4n}\right)
\nonumber
  \Phi^2.\end{eqnarray*}\end{theorem}

\subsection{The equivalence between Calabi metrics and Euclidean metrics}
\label{sect_2.4}

 In this subsection,    we  compare a Calabi metric  $G$ on a  convex  domain $\Omega\subset   \real^n$ defined by a smooth, strictly convex  function $u$ and the  standard Euclidean  metric  on $\Omega$  under the condition that the sum of the Pick invariant and the norm of the Tchebychev vector field
$$\Theta  = J + \Phi $$  is bounded  from the above.   This comparison  is based on Lemma
\ref{lemma_2.4.1} which gives an estimate for eigenvalues of    the Hessian matrix
$\left(u_{ij}\right)$ of $u$. A similar estimate was first proved in
\cite{L-J-1}, but the lemma given here is  slightly stronger.

For $p\in \Omega$, denote by  $D_a(p)$  and $B_a(p)$  the balls of radius $a$ that are  centered at
$p$ with respect to the Euclidean metric and Calabi metric respectively. Similarly,
$d_E$ and $d_u$   denote  the
distance functions with respect to the Euclidean metric and Calabi metric respectively.

First we prove the following lemma.

\begin{lemma}\label{lemma_2.4.0}
Let $\Gamma:\xi=\xi(t),t\in [0,t_0],$  be a curve from $\xi(0)=0$ to
$p_0=\xi(t_0)$  in $\Omega$.  Suppose that $\Theta \leq N^2$ along $\Gamma.$
Then for any $p\in \Gamma,$
\begin{equation}\label{eigen_est}
\left| \frac{d \log T}{ds}\right| \leq 2n^2 {N} ,\;\;\; \left|\frac{d \log \det(u_{ij})}{ds} \right|
\leq (n+2) {N}  ,
\end{equation} where   $T=Trace(u_{ij})=\sum u_{ii} $ and $s$ denotes the arc-length parameter   for the curve $\Gamma$  with
respect to the Calabi metric $G$.
\end{lemma}
{\bf Proof.} For any $p\in \Gamma,$ by an  orthogonal coordinate transformation, we can choose a coordinate system
$\xi_1,...,\xi_n$ such that $ u_{ij}(p) = \lambda_{i}\delta_{ij}.$
Denote by $\lambda_{max}$ (resp.\;$\lambda_{min}$) the maximal (resp.\;minimal) eigenvalue  of $(u_{ij})$. Suppose that, at $p$,
$$\frac{\partial}{\partial t} = \sum a_i \frac{\partial}{\partial
\xi_i}.$$  We have
\begin{eqnarray*}\frac{1}{T^2}\left(\frac{\partial T}{\partial
t}\right)^2 &=& \frac{1}{T^2}\left(\sum_{i,j} u_{iij}a_j\right)^2
\leq \frac{n}{T^2}\sum_{j}( \sum_{i}u_{iij})^2a_j^2 \nonumber \\
&\leq& \frac{n}{T^2}\left(\sum_{j}\frac{1}{\lambda_j}(
\sum_{i}u_{iij})^2\right) \left( \sum_{j} \lambda_j
a_j^2\right)\nonumber \\ &\leq& n^2 \left({\sum
u^{im}u^{jn}u^{kl}u_{ijk}u_{mnl}}\right)\left( \sum_{i,j} u_{ij}a_ia_j\right)\\
&\leq & 4n^4 \Theta\sum_{i,j} u_{ij}a_ia_j=4n^4\Theta \cdot
\left(\frac{ds}{dt}\right)^2.
\end{eqnarray*}
Here we use the definition of $\Theta$ and $J$ (cf.\;\eqref{eqn_2.2}
and \eqref{eqn_2.3}). For the third "$\leq$",  we use the following
computation:
\begin{eqnarray*}\sum_{i,j,k}
u^{im}u^{jn}u^{kl}u_{ijk}u_{mnl}&=& \sum_{i,j,k}
u_{ijk}u_{mnl}\frac{1}{\lambda_{l}}\delta^{k}_l\frac{1}{\lambda_{i}}\delta_i^m
\frac{1}{\lambda_{j}}\delta_j^n \\ & \geq &\sum_{i,k} (u_{iik})^2
\frac{1}{\lambda_{k}\lambda_{max}^2} \geq \frac{1}{n}\frac{\sum\limits_{k}(\sum\limits_i
u_{iik})^2\frac{1}{\lambda_{k}}}{(\sum\limits_i \lambda_i)^2}.
\end{eqnarray*}
Therefore, we have the differential inequality
\begin{equation*}
\left|  \frac{d \log T}{ds} \right| \leq 2n^2\sqrt{\Theta} \leq 2n^2N
.\end{equation*} By the definition of $\Phi$ and \eqref{eqn_2.3} we have $\left| \frac{d \log \rho}{ds} \right|  \leq  N.$
Then the lemma follows. $\blacksquare$

\begin{lemma}\label{lemma_2.4.1}
 Let $u$ be a smooth, strictly convex function
defined on  $ \Omega\subset  \mathbb  R^n$   and $0\in \Omega$.  Suppose that
\begin{equation}\label{eqn_2.5}
\Theta\leq \mathsf{N}^2\quad in \quad \Omega,
\end{equation}
and the Hessian matrix $(u_{ij})$ satisfies $ u_{ij}(0)=\delta_{ij}$.
Let $\Gamma:\xi=\xi(s), s\in[0,\mff a],$ be a curve lying in $\Omega,$ starting from $\xi(0)=0$ with arc-length parameter with respect to the Calabi metric $G_{u}=\sum u_{ij} d\xi_i d\xi_j$. Let $\lambda_{\min}$
and $\lambda_{\max}$ be the minimal and maximal eigenvalues of
$(u_{ij})$ along the path $\Gamma$.
Then there exists a constant  $\mff C_1$  such that
\begin{description}
\item[(i)] $ \exp\left(-\mff C_1\mff a\right)\leq
\lambda_{\min}\leq \lambda_{\max}\leq n\exp\left(
\mff C_1\mff a\right), $
\item[(ii)] $\Gamma\subset D_{\mff a\exp\left(\frac{1}{2}\mff C_1\mff a\right)}(0)$.
\end{description}
\end{lemma}
\v\n{\bf Proof.} (i)  With the initial values  given by  $ u_{ij}(0)=\delta_{ij}$, applying the  integration  to   \eqref{eigen_est} with respect to $s$, we
have
\begin{eqnarray*} && n\exp\left\{-2n^2\mff  Na \right\} \leq
T(q)\leq n\exp\left\{2n^2\mff Na
\right\}\\&&
\exp\left\{- (n+2)\mff Na \right\}\leq \det(u_{ij})(q)
\leq\exp\left\{(n+2)\mff Na \right\}
\end{eqnarray*} for any
$q\in B_a(0).$    Then we have the bounds of eigenvalues
of $(u_{ij})$.

Note that the Euclidean length of $\Gamma$ is
less than $\mff a\exp\left(\frac{1}{2}C_1\mff a\right)$. Hence we have (ii).  $\blacksquare$

\v As a corollary, we have
\begin{corollary}\label{ball_ball}
 Let $u$ be the function given in Lemma \ref{lemma_2.4.1}.
\begin{enumerate}
\item[(1)] Suppose that $D_r(0)\subset \Omega$, then there exists a constant $\mff a_1 $ depending only on $\mff N$ and $r$ such that
$B_{\mff a_1}(0)\subset D_r(0)$;
\item[(2)] Suppose that  $B_{\mff a}(0)\subset \Omega$,   then there
exists a constant $r_1$ depending only on $\mff N$ and $\mff a$ such that
$
D_{r_1}(0)\subset B_{\mff a}(0).
$
\end{enumerate}
\end{corollary}
{\bf Proof.} \begin{enumerate}
\item[(1)]  Let $v$ be any nonzero vector at $0$ and $\Gamma_v(s)$ be a geodesic ray initiating from 0 along $v$ direction  with arc-length parameter. Then by (i) of Lemma \ref{lemma_2.4.1}
$$
d_E(0,\Gamma_v(s))\leq s\exp\left(\frac{1}{2}\mff C_1s\right).
$$
Hence $\Gamma_v(s)\subset D_r(0)$ if $s\exp\left(\frac{\mff C_1s}{2}\right)<r$. Hence we may take
$$
\mff a_1=\min\left\{1, \half r\exp\left(-\frac{1}{2}\mff C_1\right)\right\}
$$ in order to get (1).

\item[(2)]   We use the same argument. Note that the eigenvalues of $(u_{ij})$ are uniformly bounded in (i) of Lemma
\ref{lemma_2.4.1}.  Let $\Gamma'_v(t)=tv$ be the ray and assume that $|v|=1$. Hence $t$ is the arc-length parameter with respect to the Euclidean metric. Then if $\Gamma'_v\subset B_\mff a(0)$,  the length of $\Gamma'$ with respect to
Riemannian metric $G_u$ is less than $t\exp\left(\frac{1}{2}\mff C_1 \mff a\right)$. Hence we have the estimate for the Riemannian distance with respect to the metric $G_u$:
$$
d_{u}(0,\Gamma'_v(t)) \leq t\sqrt{n}\exp\left(\frac{1}{2}\mff C_1 \mff a\right).
$$
Take $r_1=\frac{1}{2\sqrt{n}} \mff a\exp\left(-\frac{1}{2}\mff C_1 \mff a\right)$, we have (2).
$\blacksquare$

\end{enumerate}

\v
With these preparations, now  we can establish   the following  convergence theorem.
\begin{theorem}\label{theorem_2.4.2}
Suppose that $\{u\indexm\}_{k=1}^\infty$ is a sequence of smooth
strictly convex functions on $\Omega \subset \real^n$ containing 0
and $\Theta_{u\indexm}\leq\mff  N^2;$
suppose that $u\indexm$ are already normalized such that
$$
u\indexm\geq u\indexm(0)=0,\;\;\;
\partial^2_{ij}(u\indexm)(0)=\delta_{ij}.$$ Then
\begin{enumerate}
\item[(a)] there exists a constant  $\mff a>0$ such that $B_{\mff a,u_k}(0)\subset \Omega $,
\item[(b)] there exists a subsequence of $u\indexm$ that locally
$C^2$-converges to a strictly convex function $u_\infty$
in $B_{\mff a,u_\infty}(0)$,
\item[(c)] moreover,  if $(\Omega ,G_{u\indexm})$ is complete,
$(\Omega ,G_{u_\infty})$ is complete.
\end{enumerate}
\end{theorem}
{\bf Proof.} By (1) in Corollary \ref{ball_ball}  we have a constant
$\mff a_1$ such that
$$
B_{\mff a_1,u_k}(0)\subset D_r(0)\subset \Omega.
$$
Now by (2) of Corollary \ref{ball_ball},  there exists a constant $r_1$ such that
$$
D_{r_1}(0)\subset B_{\mff a_1,u_k}(0)\subset \Omega.
$$
 It follows from  $\Theta\leq \mff N^2$
that  $u\indexm$
$C^2$-converges (passing to a subsequence) to a strictly convex function $u_\infty$ in $D_{r_2}(0),$ where $r_2=\half r_1$. Take $\mff a$ to be a constant such that
$$
B_{\mff a,u_k}(0)\subset D_{r_2}(0).
$$
Then we have proved (a) and (b).

Now we prove (c), i.e,  if $(\Omega ,G_{u\indexm})$ is complete, then
$(\Omega ,G_{u_\infty})$ is complete.
We already conclude
that $u\indexm$
$C^2$-converges to a strictly convex function $u_\infty$ in $D_{r_2}(0),$ and  $B_{\mff a, u_\infty}(0)\subset D_{r_2}(0)$.

For any point $q\in \partial B_{\mff a,u_\infty}(0),$
we normalize $u_k$ by the following transformation:
 $$\tilde \xi=A_k(\xi-\xi(q)),\;\;\;\; \tilde u_k(\tilde \xi)=u_k(\xi)-\nabla u_k(q)\cdot(\xi-\xi(q))-u_k(q),$$
 where  $A_k$ is the matrix such that $\partial^2_{ij}\tilde u_k (0)=\delta_{ij}.$

From the convergence of $u\indexm$ in $D_{r_2}(0)$ we have the following convergence
\begin{equation}\label{eqn_add_2.8}
\lim_{k\to\infty}A_k=A_{\infty} ,\;\;\;\;\lim_{k\to\infty}\nabla u_k(q)=\nabla u_{\infty}(q),\;\;\; \lim_{k\to\infty} u_{k}(q)=  u_{\infty}(q),
\end{equation}
where $A_{\infty}$ is the matrix such that
 $\partial^2_{ij}\tilde u_{\infty}(0)=\delta_{ij}.$

 By the same argument above we conclude that  $\tilde u\indexm$
$C^2$-converges to a strictly convex function $\tilde u_\infty$ in   $B_{\mff a,\tilde u_\infty}(0).$ Since  $ \partial B_{a,u_\infty}(0)$ is compact, by \eqref{eqn_add_2.8},
$u\indexm$
$C^2$-converges to a strictly convex function $u_\infty$ in  $B_{2\mff a,u_\infty}$.
Note that $\mff a$ is a uniform constant, hence by repeating this procedure we can prove (c).
$\blacksquare$

\v    The following theorem provides an equivalence between the Calabi metric and the Euclidean metric in  a convex domain $\Omega$  under the condition that $\Theta$ is bounded in the sense that
the Calabi metric is complete if and only if the Euclidean metric is complete.

\begin{theorem}\label{theorem_2.4.3}
Let $\Omega\subset \real^n$ be a convex domain. Let $u$ be a smooth strictly convex function on
$\Omega$ with  $\Theta\leq\mff N^2$. Then $
(\Omega,G_u)$ is
 complete if and only if  the graph of $u$ is Euclidean complete
in $\real^{n+1}$.
\end{theorem}
{\bf Proof.} "$\Longrightarrow$". Without loss of generality, we
assume that $0\in \Omega$ and
$$
u\geq u(0),\;\;\; u_{ij}(0)=\delta_{ij}.
$$
We will prove that $u|_{\partial\Omega}=+\infty$. Namely, for any Euclidean
unit vector $v$ and the ray from 0 to  $\partial\Omega$ along $v$
direction, given by
$$\ell_v: [0,t)\to \Omega, \;\;\; \ell_v(s)=sv,\;\;\; \lim_{s\to t}\ell_v(s)\in \partial \Omega,
$$
we have $u(\ell_v(s))\to \infty$ as $s\to t$.

On $[0,t)$, we choose a sequence of points
$$
0=s_0<s_1<\cdots<s_k<\cdots
$$
such that $d_u(\ell_v(s_i),\ell_v(s_{i+1}))=1$. Since $\Omega$ is
complete, the sequence is an infinite sequence. Set
$u_i=u(\ell_v(s_i)).$ We now claim that \v\n {\em Claim. There is a
constant $\delta>0$ such that $u_{i+1}-u_i\geq \delta$.} \v {\em
Proof of claim.} We show this for $i=0$. By Lemma \ref{lemma_2.4.1},
we conclude that $B_1(0)$ contains an Euclidean ball $D_{a_1}(0)$.
Since the eigenvalue of $(u_{ij})$ are bounded from below and
$\nabla u(0)=0$, we conclude that there exists a constant $\delta$
such that $u(\xi)-u(0)\geq \delta$ for any $\xi\in B_1(0)-
D_{a_1(0)}$. In particular, $u_1-u_0\geq \delta$. Now consider any
$i$. Let $p=\ell_v(s_i)$ and $q=\ell_v(s_{i+1})$. Since the claim is
invariant with respect to the base-affine transformations, we first
apply a base-affine transformation such that
$u_{ij}(p)=\delta_{ij}$. Furthermore, we normalize $u$ to a new
function $\tilde u$ such that $p$ is the minimal point:
$$
\tilde u(\xi)=u(\xi)-u(p)-\nabla u(p)(\xi-p).
$$
Then by the same argument as $i=0$ case, $\tilde u(q)-\tilde u(
p)\geq \delta$. Moreover
$$
u(q)-u(p)=\tilde u(q)-\tilde u(p) + \nabla u(p)(q-p)\geq \tilde
u(q)-\tilde u(p)\geq \delta.
$$
Here we use the convexity of $u$ for the first inequality.  This completes the
proof of the claim. \v By the claim, $u_k\geq k\delta$ and goes to
$\infty$ as $k\to \infty$.

\v\n "$\Longleftarrow$". We assume that $ u \geq u (0)=0 .$ Since
the graph of $u$ is Euclidean complete, for any $C>0,$ the section $S_u(0,C)$
is compact (the definition of section is given in \S\ref{sect_3.1}). Consider the function
\[ F=\exp\left\{-\frac{2C}{C-u}\right\}\frac{\|\nabla u\|_G^2}{(1+u)^2}\]
defined in $S_u(0,C)$. $F$ attains its maximum at some interior
point $p^*$. We may assume $\|\nabla u\|_G(p^*)>0.$ Choose a frame
field of the Calabi metric $G$  around $p^*$ such that
 $$u_{ij}(p^*)=\delta_{ij},\;\;\|\nabla u\|_G(p^*)=u_1(p^*)
,\;\;  u_i(p^*)=0, \; i\geq 2 .$$    Then at $p^*$, $F_i=0$
implies that
\begin{eqnarray*}
   \left(-\frac{2u_{,1}}{1+u}-  \frac{2Cu_{,1}}{(C-u)^2}\right)
   u_{,1}^2 + 2  u_{,1}u_{,11}= 0.
\end{eqnarray*}
     Using  $u_1(p^*)>0$ and $ \frac{2C}{(C-u)^2}>0,$ we
     conclude  that
\begin{equation}\label{eqn_2.6}
\frac{u_{,1}^2}{(1+u)^2}
   \leq  \frac{|u_{,11}|}{1+u}.\end{equation}
   Note that
\begin{equation*}|u_{,11}| = \left|  \Gamma^1_{11}u_{,1} + u_{11}\right| = |1+u_{111}u_{,1}|\leq 1+
|2n \sqrt{J}|u_{,1}\leq 1+2n\mff Nu_{,1} .\end{equation*} Substituting
this into \eqref{eqn_2.6} and applying the Schwarz's inequality   we
have
\begin{equation*}\frac{u_{,1}^2}{(1+u)^2}\leq 2+8n^2\mff N^2.
  \end{equation*}
Then $F(p^\ast)\leq e^{-2}(2+8n^2\mff N^2).$ Since $F$ is
maximum at $p^*$, we obtain that $F\leq C_1$.  Hence in
$S_{u}(0,\frac{C}{2}),$
\begin{equation}\label{eqn_2.7}
 \|\nabla \log (1+ u)\|_G\leq C_2, \;\;\;
\end{equation}
for some constant $C_2$ that is  independent of $C.$ Letting $C\to
\infty$, we have that \eqref{eqn_2.7} holds everywhere.

Take a point  $p_1\in \p S_{u}(0,{C})$ such that $d(0,p_1)=d(0,\p
S_{u}(0,{C})).$ Let $l$ be the shortest geodesic from $0$ to
$p_1.$ Following from \eqref{eqn_2.7} we have $$C_2  \geq \|\nabla \log
(1+ u)\|_G  \geq  \left|\frac{d\log (1+u)}{ds}\right| ,$$ where $s$ denotes the arc-length parameter with
respect to the metric $G$. Applying this and a direct
integration we obtain
\[d(0,p_1)= \int_l ds \geq C_2\inv\log(1+u)|_{p_1}=
 C_2\inv \log\left(1+C\right) . \]
As $C\to \infty$, we obtain $d(0,\p S_u(0,C))\to +\infty.$
Hence
 $(\Omega,G_u)$ is complete.    $\blacksquare$


\section{Convergence theorems and   Bernstein properties  for
  Abreu equations}\label{sect_3}

\subsection{Basic definitions}\label{sect_3.1}
 In this subsection we
review some basic concepts of convex domains and (strictly)
convex functions.

Let $\Omega \subset \mathbb{R}^2$ be a bounded convex domain. It
is well-known (see \cite{G}, p.27) that there exists a unique
ellipsoid $E$, which attains the minimum volume among all the
ellipsoids that contain $\Omega $ and that is centered at the
center of mass of $\Omega $, such that
$$2^{-\frac{3}{2}} E \subset \Omega  \subset E,$$ where
$2^{-\frac{3}{2}} E$ means the $2^{-\frac{3}{2}}$ -dilation of $E$
with respect to its center. Let $T$ be an affine transformation
such that $T(E)=D_1(0)$, the unit disk. Put $\tilde{\Omega}=
T(\Omega)$. Then
\begin{equation}\label{eqn_3.1}
2^{-\frac{3}{2}}D_1(0) \subset \tilde{\Omega} \subset
D_1(0).\end{equation} We call $T$ the normalizing transformation of
$\Omega$.
\begin{defn}\label{definition_3.1.1}
A convex domain $\Omega $ is called normalized when its center of mass is 0 and $2^{-\frac{3}{2}}D_1(0) \subset  \Omega  \subset
D_1(0).$
\end{defn}
Let $A:\real^n\to\real^n$ be an affine transformation
given by $A(\xi)=A_0(\xi)+a_0$, where $A_0$ is a linear transformation
and $a_0\in \real^n$. For any Euclidean vector $a,$ denote by $|a|$ the Euclidean norm of the vector $a.$  If there is
a constant $L>0$ such that $|a_0|\leq L$ and for any Euclidean  unit vector $v$
$$
L\inv\leq |A_0v|\leq L,
$$
we say that  $A$ is  $L$-bounded.
\begin{defn}\label{definition_3.1.2}
A convex domain $\Omega$ is called $L$-normalized if its
normalizing transformation is  $L$-bounded.
\end{defn}
The following lemma is useful to measure the normalization of a
domain.
\begin{lemma}\label{lemma_3.1.3}
 Let $\Omega \subset \real^2$ be a  convex domain. Suppose
 that there exists a pair of constants $R>r>0$ such that
 $$
D_r(0)\subset \Omega\subset D_R(0),
 $$
then $\Omega$ is $L$-normalized, where $L$ depends only on $r$
and $R$.
\end{lemma}
{\bf Proof.} Suppose the normalizing transformation is $T(\xi)=A(\xi-p),$ where $A$ is a linear
transformation and $p\in \Omega$. Then $$T(D_R(0))\supset
2^{-\frac{3}{2}}D_1(0)  \Rightarrow A(D_R(0)-p)\supset 2^{-\frac{3}{2}}D_1(0)
 \Rightarrow A(D_{2R}(0))\supset 2^{-\frac{3}{2}}D_1(0).
$$
Here, we use the fact $p\in D_R(0)$. Therefore for any Euclidean  unit vector $v$
$$
|Av|\geq 2^{-\frac{5}{2}}R\inv.
$$

On the other hand,
$$
T(D_r(0))\subset D_1(0)\Rightarrow A(D_r(0))-A p\subset D_1(0).$$
In particular,   we have
$-Ap\in D_1(0).$ Then
$$A (D_r(0))\subset   D_2(0).$$
This says that for any Euclidean  unit vector $v$
$$
|Av|\leq \frac{2}{r}.
$$
We get  upper and lower bounds of eigenvalues of $A$. Furthermore
$
|Ap|\leq 1.
$ So  $T$ is $L$-bounded for some constant $L$.
$\blacksquare$

\v Let $\Omega$ be a convex domain in $\real^2$. We  define its {\em uniform $\xi_i$-width} $wd_i(\Omega)$, $i=1,2,$ by the following:
$$
wd_1(\Omega)=\max_{t\in \real} \max_{\{\xi,\xi'\in
\Omega|\xi_2=\xi'_2=t\}} |\xi_1-\xi'_1|,$$
$$
wd_2(\Omega)=\max_{t\in \real} \max_{\{\xi,\xi'\in
\Omega|\xi_1=\xi'_1 =t\}} |\xi_2-\xi'_2|.
$$
Then
\begin{lemma}\label{lemma_3.1.4}
Let $\Omega$ be a bounded convex domain in $\real^2$ with
$\Omega \subset D_b(0)$. Suppose that $ wd_i(\Omega) \geq a>0$.
Then $\Omega$ is $L$-normalized, where $L$ depends only on $a$ and
$b$.
\end{lemma}
{\bf Proof.} Since $wd_i$ has lower bounds, there exists a disk of
radius $r$ in $\Omega$, where $r$ depends only on $a$. Let $p_0$
be the center of disk. Then by $\Omega\subset D_{b}(0)$ it is
obvious that $\Omega\subset D_{2b}(p_0)$. Then this lemma is a
corollary of Lemma \ref{lemma_3.1.3}. $\blacksquare$

\v

Let $u$ be a convex function on $\Omega$. Let $p\in \Omega$.
Consider the set
$$
\{\xi\in \Omega| u(\xi)\leq u(p)+\nabla
u(p)\cdot(\xi-p)+\sigma\}.
$$
If it is compact in $\Omega$, we call it {\em a section of $u$ at
$p$ with height $\sigma$} and denote it by $S_u(p,\sigma)$.

In this section, we need an important result from the classical
convex body theory (see \cite{BU}). Let $M$ be a convex
hypersurface in $\mathbb R^{n+1}$ and $e$ be a subset of $M$. We
denote by $\psi_M(e)$ the spherical image of $e$. Denote by
$\sigma_M(e)$ the area (measure) of the spherical image $\psi_M(e)$,
denote by $A(e)$ the area of the set $e$ on $M$. The ratio
$\sigma_M(e)/A(e)$ is called the specific curvature of $e$. In the
case $n=2$, the following theorem holds  (see (5.5) in Page 35 of \cite{BU}):
 \begin{theorem}[Alexandrov-Pogorelov]\label{theorem_3.1.5}A convex
surface whose specific curvature is  bounded   away from zero is
strictly convex.   \end{theorem}

\subsection{The convergence theorem for Abreu equations}\label{sect_3.2}

In this subsection we  prove a very useful convergence theorem for the Abreu
equations (cf.\;Theorem \ref{theorem_3.2.7}). Essentially, there are two types of convergence
theorems in this paper for the Abreu equations. One is to get the
convergence by controlling  $\Theta$  (cf.\;Theorem \ref{theorem_2.4.2}), the other is to
get the convergence by controlling the sections. Here, we
explain the latter one.

Denote by $\mathcal{F}(\Omega,C)$ the class of smooth convex functions
defined on $\Omega$ such that
$$ \inf_{\Omega} u  = 0,\;\;\;
u= C>0\;\;on\;\;\partial \Omega.$$ Let 
$$
\mc F(\Omega,C;K_o)=\{u\in \mc F(\Omega,C)| |\mc S(u)|\leq K_o\},
$$
where $\mc S(u)=-\sum U^{ij}w_{ij}.$ Recall that $(U^{ij})$ is the cofactor matrix of the Hessian matrix
$(u_{ij})$ and $w=(\det(u_{ij}))\inv$.
We will assume that  $\Omega\subset \mathbb R^2$ is normalized in this subsection.

The main result of this subsection is the following convergence
theorem.
\begin{theorem}\label{theorem_3.2.7}
Let $\Omega\subset \mathbb R^2$ be a  normalized domain. Let $u\indexm\in \mc
F(\Omega,1; K_o)$ be a sequence of functions and $p^o\indexm$ be
the minimal point of $u\indexm$. Then there exists a subsequence
of functions, without loss of generality, still denoted by
$u\indexm$,  locally uniformly  converging to a function $u_\infty$ in $\Omega$ and
$p^o\indexm$ converges to $p^o_\infty$ satisfying:
\begin{enumerate}
\item[(i)] there exist constants $\mff s$ and $\mff C_2$ such that $d_E(p\indexm^o,\partial\Omega)>2\mff s$, and in
$D_{\mff s}(p^o_\infty)$
$$
\|u\indexm\|_{C^{3,\alpha}}\leq \mff C_2
$$
for any $\alpha\in (0,1)$;
in particular,
 $u\indexm$ $C^{3,\alpha}$-converges
to $u_\infty$ in $D_\mff s(p^o_\infty)$.\item[(ii)] there exists a
constant $ \delta\in (0,1),$ such that $S_{u\indexm}(p^o\indexm,
\delta)\subset D_\mff s (p^o_\infty). $ \item[(iii)] there
exists a constant $\mff b>0$ such that
$S_{u\indexm}(p^o\indexm,\delta)\subset B_\mff b(p\indexm^o)\subset \Omega.$
\end{enumerate}
Furthermore, if $u\indexm$ is smooth and $C^k$-norms
of $\mc S(u\indexm)$ are uniformly bounded, then $u\indexm$
$C^{k+3,\alpha}$-converges to $u_\infty$ in
$D_\mff s(p^o_\infty)$.
\end{theorem}

\begin{remark}\label{remarka_3.2.7}
Let  $\Omega_k$ be a sequence of  L-normalized convex  domains and $u\indexm\in \mc
F(\Omega_k,C_k; K_o)$ be a sequence of functions. If the sequence of constants $C_k$ satisfies
$C\inv\leq C_k\leq C$ for some constant $C>0,$     then Theorem \ref{theorem_3.2.7} still holds.
\end{remark}
\begin{remark}\label{remark_?}
This theorem with $K_o=0$ has been proved in \cite{L-J-4}, the proof
is  purely analytical and consists of a very long calculation. The authors of \cite{L-J-4}  pointed out that one can use the
Alexandrov-Pogorelov Theorem (Theorem 3.5) to give a simpler proof.
Here we follow this suggestion to give the proof for $K_o\ne 0$.
\end{remark}

We prove the following proposition, which is equivalent to Theorem
\ref{theorem_3.2.7}.
\begin{prop}\label{proposition_3.2.6}
Let $\Omega\subset \mathbb R^2$ be a  normalized domain. Let $u\in \mc F(\Omega,1;
K_o)$ and $p^o$ be its minimal point. Then
\begin{enumerate}
\item[(i)] there exist constants $\mff s$ and $\mff C_2$ such that $d_E(p^o,\partial\Omega)>\mff s$ and in
$D_\mff s(p^o)$
$$
\|u\|_{C^{3,\alpha}}\leq\mff  C_2
$$
for any $\alpha\in (0,1)$;
 \item[(ii)] there exists a constant
$0< \delta <1$ such that $S_u(p^o,\delta)\subset D_\mff s
(p^o). $ \item[(iii)] there exists a constant $\mff b>0$ such that
$S_u(p^o,\delta)\subset B_\mff b(p^o)\subset \Omega.$
\end{enumerate}
\end{prop}
In the statement, all the constants  only depend on $K_o$.

Furthermore, if $u$ is smooth, then for any $k\in \mathbb Z^{+},$
$$
\|u\|_{C^{k+3,\alpha}(D_\epsilon(p^o))}\leq \mff C_3
$$
where $\mff C_3$ depends on the $C^k$-norm of $\mc S(u)$.

To prove Proposition \ref{proposition_3.2.6},
we
need some useful  properties for functions in $\mc F(\Omega,1;K_o)$ that are stated in
Lemma \ref{proposition_3.2.1}-Lemma \ref{lemma_3.2.5}.

\begin{lemma}\label{proposition_3.2.1}
Let $u\in \mc F(\Omega, 1;K_o)$. If there is a constant $C_1>0$ such
that in $  \Omega$ \begin{equation}\label{eqn_A}C_1\inv\leq
\det(u_{ij})\leq C_1,\end{equation}
 then for any $\Omega^*\subset\subset \Omega $, $p>2$, we have the
estimate \begin{equation}\|u\|_{W^{4, p}(\Omega^*)}\leq C ,\;\;\;\;\|u\|_{C^{3,\alpha}(\Omega^*)}\leq C,
\end{equation} where $C$ depends on $n, p,C_1, K_o$, $d_E(\Omega^*,
\partial \Omega )$.
\end{lemma}
{\bf Proof.} The argument is similar to \cite{TW}. In \cite{CG} Caffarelli-Gutierrez  proved a H\"older estimate of
$\det(u_{ij})$ for homogeneous linearized Monge-Amp\`ere equations
assuming that the Monge-Amp\`ere measure $\mu[u]$ satisfies some
condition, which is guaranteed by \eqref{eqn_A}.  Consider the Abreu equation
$$ \sum U^{ij}w_{ij}=-K
,\;\;\; w=[\det(u_{kl})]\inv,$$ where $K\in
L^\infty(\Omega).$ By the same  argument in \cite{CG}  one can obtain the   H\"older
continuity of $\det(u_{ij})$. Then Caffarelli's $C^{2,\alpha}$
estimates for Monge-Amp\`ere   equations \cite{C1} give  us
$$\|u\|_{C^{2,\alpha}(\Omega^*)}\leq C_2.$$
Then $U^{ij}\in C^{\alpha}(\Omega^*).$ Following from the standard elliptic regularity theory we have
$\|u\|_{W^{4, p}(\Omega^*)}\leq C $. By the Sobolev embedding theorem
\[\|u\|_{C^{3,\alpha}(\Omega^*)}\leq C_2  \|u\|_{W^{4, p}(\Omega^*)}.
\]
Then  the proposition follows.
$\blacksquare$

\begin{lemma} \label{lemma_3.2.2}
Let $u\in \mc F(\Omega, 1;K_o)$. Let $\Omega^\ast=\nabla^u(\Omega)$
and $f=L(u)$. For any $\Omega_1^\ast\subset\subset \Omega^\ast$
there exists a constant $\mff d>0$ such that  $\det D^2u>\mff d$ on
$\Omega_1=\nabla^f (\Omega^\ast_1)$, where $\mff d$ depends on $\Omega,
d_E(\Omega_1^\ast,\partial\Omega^\ast),diam(\Omega_1^\ast)$ and
$K_o$.
\end{lemma}
This is proved by Li-Jia in \cite{L-J-4}.

\begin{lemma} \label{lemma_3.2.3}
Let $u\in \mc F(\Omega,1;K_o)$.  There exist two  constants $\mff d_1,
\mff  d_2>0$
such that
\begin{equation}\label{eqn_3.2}
\exp\left\{-\frac{4}{1-u}\right\}\frac{\det
(u_{ij})}{(\mff d_1+f)^{4}}\leq\mff d_2.
\end{equation}
\end{lemma}
This is proved by Chen-Li-Sheng in \cite{CLS1}.

\begin{lemma}\label{lemma_3.2.4}
For $u\in \mc F(\Omega,1;K_o)$ with $u(p_o)=0$. Let $R$ be a constant such
that
\begin{equation}\label{eqn_3.3}
D_R(p_o)\supset \Omega\supset S_u(p_o,\frac{1}{2}).
\end{equation}
Then $\nabla^u(S_u(p_o,\frac{1}{2}))$ contains the disk $D_r(0),
r=(2R)\inv$. Moreover, in $D_r(0)$
\begin{enumerate}
\item[(i)] $-1\leq f\leq r$;
\item [(ii)]
$\det(f_{ij})\geq d'$ for some constant $d'$ depending only on
$K_o$ and $R$.
\end{enumerate}
\end{lemma}
{\bf Proof.} The first statement follows from \eqref{eqn_3.3}.
 At $\nabla^u(0),$
$$f_1(\nabla^u(0))=f_2(\nabla^u(0))=0,\;\;\;\;u(0)+f(\nabla^u(0))=0.$$ By the convexity of $f$, we have that the minimum of $f$ is $f(\nabla^u(0)).$ Then
$$f\geq f(\nabla^u(0))=-u(0)\geq -1.$$
 Since $ f(x)+u(\xi)=x\cdot\xi, $\;$u\geq 0$ and
$\Omega\subset D_1(0)$, we have
$$
f(x)\leq x\cdot\xi\leq r\cdot 1=r.
$$
From \eqref{eqn_3.2} and (i), we know that $ \det(u_{ij}) $ is bounded from
above in $D_r$. Therefore $\det(f_{ij})=[\det(u_{ij})]\inv$ is
bounded from below.   $\blacksquare$

\begin{lemma}\label{lemma_3.2.5} Use the notations in Lemma
\ref{lemma_3.2.4}. Then for any $\alpha\in (0,1)$, there exists a constant $C_2$ such
that in $D_{\frac{r}{2}}(0)$
\begin{equation}\label{eqn_3.4}
\|f\|_{C^{3,\alpha}}\leq C_2,
\end{equation}
and for any  eigenvalue  $\lambda_f$ of $(f_{ij})$
\begin{equation}\label{eqn_3.5}
C_2\inv\leq \lambda_f\leq C_2.
\end{equation}
\end{lemma}
{\bf Proof.} If not, we have a sequence of $u_k$ such that in
$D_{\frac{r}{2}}(0)$
\begin{equation*}
\|f_k\|_{C^{3,\alpha}}\to\infty.
\end{equation*}
By (i) of Lemma \ref{lemma_3.2.4}, we know that $f_k$ locally uniformly converges to a convex function $f_\infty$ on $D_r(0)$. The Gauss curvature $K_{Gauss}$ of the convex surface $(x_1,x_2,f(x))$ can be written as
$$K_{Gauss}(x)=\frac{\det(f_{ij})}{(1+|\nabla f|^2)^2}.$$
Since $\nabla f (D_r(0))\subset \Omega\subset D_1(0),$ we have $|\nabla f|\leq 1,$ in $D_r(0).$
By (ii) of Lemma \ref{lemma_3.2.4}, $$K_{Gauss}\geq \frac{d'}{4}.$$  Then the specific curvature of convex surface $(x_1,x_2,f_{\infty}(x))$ is bounded away from zero.
By Theorem \ref{theorem_3.1.5}, we conclude that $f_\infty$ is {\em strictly
convex} on $D_{\frac{r}{2}}(0)$. Hence there is a constant $b_0$
such that $S_{f_\infty}(0,b_0)$ is compact in
$D_{\frac{r}{2}}(0)$.  Then $S_{f\indexm}(0,b_0)\subset \subset D_{\frac{r}{2}}(0),$ when $k$ is large enough. By $|\nabla f|\leq 1 $ in $D_{r}(0),$ we have $$d_E(S_{f_{k}}(0,b_{0}/2),\partial S_{f_{k}}(0,b_{0}))\geq \tau$$ for some positive constant $\tau$ independent of $k.$ Moreover, in the section $S_{f_{k}}(0,b_{0}/2)$,
$\det(D^2f_k)$ is bounded from above (by Lemma \ref{lemma_3.2.2}) and below (by (ii) of Lemma \ref{lemma_3.2.4}). Hence by Lemma \ref{proposition_3.2.1}, we conclude that the $C^{3,\alpha}$-norms of $f_k$ in the section $S_{f_k}(0,\frac{1}{2}b_0)$ are uniformly bounded by a constant $N$ that is independent of $k$. Instead of considering $0$, we repeat the argument for any $x\in D_{\frac{r}{2}}(0)$ we conclude that the $C^{3,\alpha}$-norms of $f_k$ are uniformly bounded in some section $S_{f_k}(x,b_x)$ for some $b_x$. By the compactness of $D_{\frac{r}{2}}(0)$, we conclude that the $C^{3,\alpha}$-norms of $f_k$ in $D_{\frac{r}{2}}(0)$ are uniformly bounded. This contradicts  the assumption. $\blacksquare$

\v
\n
{\bf Proof of Proposition \ref{proposition_3.2.6}.} Due to the lower bound of eigenvalues in
Lemma \ref{lemma_3.2.5}, there exists a constant $b$ such that for
any $f$
$$
S_f(0,b)\subset\subset D_{\frac{r}{2}}(0).
$$
Let $V_f= \nabla^f(S_f(0,b))$. Since $f$ is strictly convex and
$C^{3,\alpha}$-bounded, $V_f$ must contain a disk
$D_{\epsilon}(p_o)$, where $\epsilon$ is independent of $f$. The
regularity of $u$ in $V_f$ follows from Proposition \ref{proposition_3.2.1}. (iii) is a
direct consequence of (i) and (ii). $\blacksquare$

\subsection{Bernstein properties for Abreu equations}\label{sect_3.3}
When $\mathcal  S(u)=0$, we expect that $u$ must be a quadratic
polynomial under certain completeness assumptions. This is called
the Bernstein property.

We state the results here. The proof of these theorems can be found in \cite{L-J}.
\begin{theorem}\label{theorem_3.3.1}[Jia-Li] Let $\Omega\subset \real^n$ and
$f:\Omega\to \real$ be a smooth strictly convex function.
Suppose that $(\Omega,G_f)$ is complete and
$$
\mc S(f)=-\sum f^{ij}(\log\det( f_{kl}))_{ij}=0
$$  If $n\leq 5$,  then the graph of $f$ must be an elliptic paraboloid.
\end{theorem}

\begin{remark}\label{remark_3.3.1}   If we only assume that $f\in C^{5}$, Theorem \ref{theorem_3.3.1} remains  true.
\end{remark}
Now we use the convergence theorem \ref{theorem_3.2.7} to prove the
following theorem.
 \begin{theorem}\label{theorem_3.3.3}[Li-Jia] Let $u(\xi_1,\xi_2)$ be a $C^\infty$
 strictly convex function defined in a convex
domain $ \Omega \subset \real^2$. If
 $$
\mc S(u)=0,\;\;\;
u|_{\partial\Omega}=+\infty,$$
then the graph of $u$ must be an elliptic paraboloid.
\end{theorem}

\vskip 0.1in \noindent {\bf Proof.}  First we show  $\Phi \equiv 0.$ If
this is not true, then there is a point $p \in M$ such that
$\Phi(p)>0.$ By subtracting a linear function we may suppose that
$u(\xi)\geq u(p)=0.$

Choose a sequence $\{C_k\}$ of positive numbers such that
$C_k\rightarrow\infty$ as $k\rightarrow\infty$. For any $C_k>0$ the
section $S_u(p,C_k)$ is a bounded convex domain. Let
$$u_k(\xi)=C_k\inv{u(\xi)},\quad k=1,2,\dots$$
Then
\begin{equation}\label{eqn_3.12a}
\Phi_{u_k}(p)=C_k\Phi(p)\to \infty,\;\;\; k\to\infty.
\end{equation}
Let $T_k$ be the normalizing transformation of $S_u(p,C_k)=S_{u_k}(p,1)$.
We obtain a sequence of convex functions
$$\tilde{u}_k({\xi}):=
u_k(T_k\inv(\xi)).$$
Let $\tilde p_k=T_k(p).$
By Theorem \ref{theorem_3.2.7}
 we conclude that $\tilde u^{(k)}$ converges with all
derivatives to a smooth and strictly convex function $u_\infty$
in a neighborhood of $\tilde p_\infty$, which is the limit of $\tilde p_k$.
Therefore $
\Phi_{u_k}(p)=\Phi_{\tilde u_k}(\tilde p_k)
$ (cf.\;Lemma \ref{lemma_2.2.1}) is uniformly bounded, which contradicts
to \eqref{eqn_3.12a}.

  {
  It follows from  $\Phi\equiv 0$ that
   $$\det(u_{ij})=Const.$$  Since $\lim\limits_{p\to \partial \Omega} u(p)=+\infty,$ the graph of $u$ is Euclidean complete. By a   theorem of Cheng-Yau (p.118 Theorem 2.3.17 \cite{L-S-Z}) the graph of $u$ is also an affine complete parabolic affine hypersphere. Using a theorem of Calabi (p.128 Theorem 2.6.6 \cite{L-S-Z}), $u$ must be a quadratic polynomial.
}
 $\blacksquare$


\section{Affine blow-up analysis and interior estimates of $\Theta$}
\label{sect_4}
When considering the convergence of a sequence of functions, we may encounter the case that the sequence blows-up  somewhere and the convergence fails.
Then we may apply blow-up analysis to  renormalize the domain and  the functions
such that the new sequence converges nicely after normalization.  In this paper, we apply the
affine transformations (cf.\;Definition \ref{defn_2.2.1}) for normalization.  This is a
 very powerful tool to estimate {\em affine
invariants} in various circumstances. Hence, we call the technique
the ``affine blow-up analysis".

\subsection{Affine blow-up analysis}\label{sect_4.1}

Let $\mc C$ be a class of  smooth strictly convex functions on
$\Omega$. For simplicity, we suppose that $(\Omega, G_u)$ is complete for any $u\in
\mc C$.

Given a function $u\in \mc C$, let $Q_u: \Omega\to \real$ be a non-negative
function which is invariant  under base-affine transformations (cf Definition \ref{defn_2.2.1}). Suppose that we want to prove that there exists
a constant $N$, such that for any $u\in \mc C$
\begin{equation}\label{eqn_4.1}
Q_u\leq N <\infty.
\end{equation}
We usually argue by assuming that it is not true. Then
 there will be a sequence of functions $u\indexm$ and
 a sequence of points $p\indexm$ such that
\begin{equation}\label{eqn_4.2}
Q_{u\indexm}(p\indexm)\to \infty.
\end{equation}
Usually, we perform affine transformation to $u\indexm$ and get a
sequence of normalized function $\tilde u\indexm$ to get
contradictions (this is already used in the proof of Theorem \ref{theorem_3.3.3}). However, sometimes, we
need a more tricky argument to get $\tilde u\indexm$. This is
explained in the following remark.
\begin{remark}\label{remark_4.1.1}
We explain a standard blow-up argument used in our papers. Suppose
that \eqref{eqn_4.2} happens.
\begin{enumerate}
\item[(i)] Consider the function $F\indexm:=Q_{u\indexm}(p)
d^{\alpha}_{u\indexm}(p,\partial B_{1}(p\indexm))$ in the unit geodesic
ball $B_1({p\indexm})$ for some $\alpha>0$. \v  Suppose that
 $Qd^\alpha$ is invariant with respect
to affine transformations in Definition \ref{defn_2.2.1}.
 \item[(ii)] Suppose that
$F\indexm$
attains its maximal at $p^{\ast}\indexm$. Set
$$d\indexm=\frac{1}{2}d_{u\indexm}(p^{\ast}\indexm, \partial B_1({p\indexm})).$$ Then
we conclude that
\begin{itemize}
\item $Q_{u\indexm}(p^{\ast}\indexm)d^\alpha\indexm\to \infty$. \item  $
Q_{u\indexm}\leq 2^\alpha Q_{u\indexm}(p^{\ast}\indexm)$ in
$B_{d_k}({p^{\ast}\indexm})$.
\end{itemize}
We now restrict ourself on $B_{d_k}({p^{\ast}\indexm})$; 
\item[(iii)]
Without loss of generality, by translation, we set $p_k^\ast$ to
0. Now we re-normalize  functions $u\indexm$ by  proper affine
transformations such that for  the new functions, denoted by $\tilde
u\indexm$, satisfy
$$
Q_{\tilde u\indexm}(0)=1.
$$
Then the ball  $B_{d_k}(0)$ is normalized to a ball
$B_{\tilde d_k}(0)$, where $\tilde d_{k}\to\infty$.
\end{enumerate}
To summarize, we have a sequence of normalized functions $\tilde
u\indexm$ and we will study the convergence of the sequence to get
contradiction. Usually, we conclude that $\tilde u\indexm$
uniformly converges to $\tilde u_\infty$ in a neighborhood of $0$
and $Q_{\tilde u_\infty}\equiv0$ which contradicts  the fact that
$$
\lim_{k\to\infty}Q_{\tilde u\indexm}(0)=1.
$$
In this final step, we need the aid of convergence theorems.
\end{remark}

\subsection{Interior estimate of $\Theta$}\label{sect_4.2}
The purpose of this section is to estimate  $\Theta$ near the boundary $\partial \Delta.$
We use the blow-up analysis explained in Remark \ref{remark_4.1.1}
 to prove the following result.

\begin{prop}\label{theorem_4.2.1} Let $u$ be a smooth strictly convex function on
a bounded convex domain $\Omega\subset \mathbb R^2$ with $\|\mc S(u)\|_{C^{2}(\Omega)}<K_o$, where  $\|\cdot\|_{C^2}$ denotes the Euclidean $C^2$-norm. Suppose that
for any $p \in \Omega$, \begin{equation}\label{eqnc_4.1} d_u(p,\partial{\Omega})<\mff N
 \end{equation}for some constant $\mff N>0.$  Then
there exists a constant $\mff C_5>0$, depending only on $\Omega,\mff N$
and $K_o$, such that
\begin{equation}\label{eqn_4.4}
\Theta(p) d^2_u(p,\partial {\Omega})\leq \mff C_5,\;\;\; \forall p\in \Omega. \end{equation}
Here $d_u(p,\partial \Omega)$ is the distance from $p$ to $\partial\Omega$
with respect to the Calabi metric $G_u$.
\end{prop}
{\bf Proof.} If not, then there
exists a sequence of functions $u\indexm$ and a sequence of points
$p_k$ such that
$$
\Theta_{u\indexm}(p_k)d^2_{u\indexm}(p_k,
\partial{\Omega})\to \infty.$$
Let $B\indexn$ be the $\frac{1}{2}d_{u\indexm}(p_k,
\partial{\Omega})$-ball centered at $p_k$ and
consider the {\em affine transformationally  invariant function} (cf.\;Lemma \ref{lemma_2.2.2})
$$
F_k(p)=\Theta_{u\indexm}(p)d^2_{u\indexm}(p,
\partial{B}\indexn).
$$
$F_k$ attains its maximum at $p_k^*$. Put
$$d_k=\frac{1}{2}d_{u\indexm}(p^*_k,\partial{B}\indexn).$$
By adding linear functions we assume that
\[u_k(p_k^*)=0,\;\;\; \nabla u_k(p_k^*)=0.\;\;\; \]
 By taking a proper coordinate  translation  we may assume that the coordinate of $p^*_k$ is $0$. Then (cf.\;(ii) in
Remark \ref{remark_4.1.1})
\begin{itemize}
\item $\Theta_{u\indexm}(0)d^2_{k}\to \infty$. \item
$\Theta_{u\indexm}\leq 4 \Theta_{u\indexm}(0)$ in
$B\indexn_{d_k}(0)$.
\end{itemize}
By \eqref{eqnc_4.1}  we have
\begin{equation}\label{eqnc_4.2}\lim\limits_{k\to\infty}
\Theta_{u\indexm}(0)=+\infty.\end{equation}
We take an affine transformation   on $u_k$:
$$\xi^{\star}=A_k\xi,\;\;u^{\star}\indexm(\xi^{\star}):= \lambda_k
u_k\left(A_k\inv \xi^{\star}\right),$$ where $\lambda_k = \Theta_{u\indexm}(0).$ Choose $A_k$ such that  $\partial^2_{ij}u^{\star}_{k}(0)=\delta_{ij}.$ Denote $A_k\inv=(b^k_{ij})$. Then by the affine transformation
rule we know that
$\Theta_{u^{\star}_k}(0) = 1,$
and for any fixed large $R$, when $k$ large enough
$$\Theta_{u^{\star}_k}\leq 4
\;\;\;\;in \;\;\;B\indexn_{R}(0).$$ Moreover,
$$
\mc S(u^{\star}\indexm)= \frac{\mc S(u\indexm)}{\lambda_k}\to 0.
$$
Assume that $d\indexm$ changes to $d\indexm^{\star}$ with respect to the affine transformation. Then $$\lim_{k\to\infty}d\indexm^{\star}=+\infty.$$
By Theorem \ref{theorem_2.4.2}, one concludes that $u^{\star}\indexm$ locally  $C^{2}$-converges
to a function $u^\star_\infty,$ and the graph of $u_\infty^{\star}$ is complete with respect to the metric $G_{u_\infty}.$  By Lemma \ref{lemma_2.4.1}, we have in $B_{R,u^\star_{\infty}}(0),$
$$\exp\left(-\mff C_1 R\right)\leq
\lambda_{\min}(u^\star\indexm)\leq \lambda_{\max}(u^\star\indexm)\leq n\exp\left(
\mff C_1 R\right),$$
where $\lambda_{\min}(u^\star\indexm)$ and $\lambda_{\max}(u^\star\indexm)$ denotes the minimal and maximal eigenvlues of $(u^\star_{ij})$ in $B_{R,u^\star_{\infty}}(0).$ Then there is a constant $r_1>0$ independent of $k$ such that $S_{u^\star_{k}}(q,r_1)\subset B_{R,u^\star_{\infty}}(0)$ for any $q\in B_{R/2,u^\star_{\infty}}(0).$ It follows from Lemma \ref{proposition_3.2.1} that
$u^{\star}\indexm$ locally  $C^{3,\alpha}$-converges to a function $u^\star_\infty.$ Since $\Theta_k$ is bounded,
by  Lemma \ref{lemma_2.4.1} we know  there exists an Euclidean ball $D_\epsilon(0)$ such that  $D_\epsilon(0)\subset A_k(\Omega)$ when $k$ large enough. Therefore
 $A_k\inv D_\epsilon(0)\subset \Omega.$ It follows that for any $1 \leq i,j\leq 2$
\begin{equation}\label{eqna_4.1}\left|\frac{\partial \xi_i}{\partial \xi^\star_j}\right|=|b^k_{ij}|  \leq  {2diam(\Omega)}{\epsilon\inv}, \end{equation}
  when $k$ large.
By a direct calculation we have
$$\left|\frac{\partial \mc S(u^{\star}\indexm)}{\partial \xi_i^{\star}}\right|=
 \lambda_k\inv   \left|\sum_{j}b^k_{ji}\frac{\partial \mc S(u \indexm)}{\partial \xi_j}\right|  \leq 8diam(\Omega){\epsilon\inv}\lambda_k\inv K_o .$$Similarly, we have
$$\|\mc S(u\indexm^{\star})\|_{C^2}\leq 128[(diam(\Omega))^2\epsilon^{-2}+1]\lambda_k\inv K_o  .$$
Then  by the standard elliptic equation technique we obtain that $ u^{\star}\indexm $ locally   $C^5$-converges to $ u^{\star}_\infty$ with
$$
\mc S(u^{\star}_\infty)=0,\;\;\;\;\;\; \Theta_{u^{\star}_\infty}(0)=1.
$$  Hence by Theorem \ref{theorem_3.3.1} and Remark \ref{remark_3.3.1},
$u^{\star}_\infty$ must
be quadratic and $\Theta\equiv 0$. We get a contradiction. $\blacksquare$

\begin{theorem}\label{corollaryc_4.2.1} Let $u$ be a smooth strictly convex function on
a bounded convex domain $\Omega\subset \mathbb R^2$ with $\|\mc S(u)\|_{C^{2}(\Omega)}<K_o$, where  $\|\cdot\|_{C^2}$ denotes the Euclidean $C^2$-norm. Assume that
for any $p \in \Omega$, \begin{equation} \label{eqnc_4.3} d_u(p,\partial{\Omega})<\infty.\end{equation}
Suppose that the height of $u$ is bounded, i.e,
\begin{equation}\label{eqn_4.3}
\max_{\Omega} u-\min_{\Omega} u\leq \mff C_4
\end{equation}
 for some constant $\mff C_4>0.$ Then
there exists a constant $\mff C_5>0$, depending only on $\Omega,\mff C_4$
and $K_o$, such that
\begin{equation}\label{eqn_4.4}
\Theta(p) d^2_u(p,\partial {\Omega})\leq \mff C_5,\;\;\; \forall p\in \Omega. \end{equation}
Here $d_u(p,\partial \Omega)$ is the distance from $p$ to $\partial\Omega$
with respect to the Calabi metric $G_u$.
\end{theorem}
{\bf Proof.} We use affine blow-up analysis as we did in Theorem \ref{theorem_4.2.1}.
We repeat the arguments before \eqref{eqnc_4.2}.
 Then by \eqref{eqn_4.3} we have
\begin{equation} \label{eqn_4.a.1}
 \inf_{\partial \Omega} u_k - u_k(p^*_k)\leq \mff C_4.
\end{equation}
If \eqref{eqnc_4.2} holds, then the proof is identical as that of Proposition \ref{theorem_4.2.1}.
Now we assume that  $\Theta_{u\indexm}(0)$ is bounded above by some
constant $ C_1.$
 Then
$d_k\to\infty$.
By a base-affine transformation we can assume that $\partial^2_{ij}u_{k}(0)=\delta_{ij}.$ Hence we
conclude that $u{\indexm}$ locally $C^2$-converges to a
function $u_\infty$ and its domain $\Omega_\infty$ is complete with
respect to $G_{u_\infty}$. Then by Theorem \ref{theorem_2.4.3}, the
graph of $u_\infty$ is Euclidean complete. This contradicts the
 equation  \eqref{eqn_4.a.1}. $\blacksquare$

\begin{remark}\label{remarkc_4.2.2}
The condition \eqref{eqnc_4.3} in Theorem \ref{corollaryc_4.2.1} is not necessary. In fact, if  $(\Omega,G_{u})$ is complete; by using  the blow-up analysis  we can prove that in  $\Omega,$
$$\Theta\leq C$$ for some constant $C>0;$ by $\Theta\leq C$ and Theorem \ref{theorem_2.4.3}
we conclude that the graph of $u$ is Euclidean complete in $\mathbb R^3.$
This contradicts \eqref{eqn_4.3}.
\end{remark}
Similarly, we can prove that
\begin{corollary}\label{corollaryc_4.2.2} Let $u$ be as  in Theorem \ref{corollaryc_4.2.1}, with one modification:
$$\|\mc S(u)\|_{C^{3}(\Omega)}<K_o.$$
Then
there is a constant $\mff C_5>0$, depending only on $\Omega,\mff C_4$
and $K_o$, such that
\begin{equation}\label{eqn_4.4b}
\mc K(p) d^2_u(p,\partial {\Omega})\leq \mff C_5,\;\;\; \forall p\in \Omega \end{equation}
where $\mc K$ is the function defined in \eqref{eqn_1.1c}.
\end{corollary}
{\bf Proof.} We replace $\Theta$ by $\Theta+\mc K$ for the proof of Proposition \ref{theorem_4.2.1}
and Theorem \ref{corollaryc_4.2.1}. Follow word by word until the very last step. We need $|S(u)|_{C^3}$
to conclude the converges of $\mc K$, then we get $(\Theta+\mc K)_{u^\ast_\infty}(0)\cong 1$. This contradicts  the
Bernstein property as before. $\blacksquare$
\v
In the proof, we use the fact that
$\Theta $ and $\mc K$  share the same affine transformation rule. Though $\|\nabla\log |\mc S(u)|\|^2$ share the same affine transformation
rule as well, the proof can not go through since
$\|\nabla\log |\mc S(u)|\|^2$ is not well defined when $\mc S(u)=0$. We need more arguments.
\begin{corollary}\label{corollaryc_4.2.3} Let $u$ be as that in Theorem \ref{corollaryc_4.2.1}, with one extra condition:
$$|\mc S(u)|\geq \delta>0.$$
Then
there is a constant $\mff C_5>0$, depending only on $\Omega,\mff C_4$
and $K_o$, such that
\begin{equation}\label{eqn_4.4b}
\|\nabla\log|\mc S(u)|\|^2(p) d^2_u(p,\partial {\Omega})\leq \mff C_5,\;\;\; \forall p\in \Omega. \end{equation}
\end{corollary}
{\bf Proof.}
 Let $F(p):=\|\nabla\log |\mc S(u)|\|^2(p)d^2(p,\partial \Omega)$ and $p^\star$ be a point such that $F(p^\star)=\max\limits_{\Omega} F.$
    Since $F$ is an affine invariant function, by  a sequence of affine transformations
    $$\xi^{\star}=A\xi,\;\;u^{\star}(\xi^{\star}):= \lambda
u\left(A\inv \xi^{\star}\right)+l(\xi^{\star}),$$
    we can assume that  $$d_{u^\star}(p^\star,\partial \Omega)=1,\;\; \partial^2_{ij}u^\star(p^\star)=\delta_{ij},\;\;u^\star\geq u^\star(p^\star)=0.$$ Here $l(\xi^{\star})$ is a linear function. 
 Using Theorem \ref{corollaryc_4.2.1}, we have $\Theta(u^\star)\leq 4\mff C_5 \; in \; B_{\half}(p^\star).$ By (2) of Corollary \ref{ball_ball}, we have $D_{r_1}(p^\star)\subset B_{\half}(p^\star).$  Then as in \eqref{eqna_4.1} we have $\left|\frac{\partial \xi_k}{\partial   \xi^\star_i} \right|\leq 2diam(\Omega)r_1^{-1}.$
  Since \begin{equation*}\|\nabla \log |\mc S(u^\star)|\|^2_{u^\star}=\|\nabla \log |\mc S(u)|\|_{u^\star}^2=\sum {u^\star}^{ij}\frac{\partial \xi_k}{\partial   \xi^\star_i} \frac{\partial \xi_l}{\partial   \xi^\star_j} \frac{\partial \log |\mc S(u)|}{\partial \xi_k}    \frac{\partial \log |\mc S(u)|}{\partial \xi_l}, \end{equation*}
 by   $\|\mc S(u)\|_{C^2(\Omega)}\leq K_o$ and $|\mc S(u)|\geq \delta$ we conclude that $F$ has a  uniform upper bound.
$\blacksquare$


\section{Complex differential inequalities}\label{sect_5} In this section
we extend the affine techniques to real functions defined on
  domains in $\mathbb C^n$. Let $\Omega \subset \mathbb C^n$, denote
$$\mc R^\infty(\Omega):=\{f\in C^{\infty}(\Omega)\;|\; f\mbox{ is a
real function and }(f_{i \bar j})>0 \},
$$ where $(f_{i\bar j})=\left(\frac{\partial^2  f}
{\partial z_i\partial \bar{z}_j}\right).$ For $ f\in \mc R^\infty(\Omega),$  $(\Omega, \omega_f)$ is
 a  K\"{a}hler  manifold.

 For the sake of notations, we set
\begin{equation}\label{eqn_5.2}
W=\detf,\;\;\; V=\log\detf, \;\;\;
\end{equation}
 Introduce the functions
\begin{equation}\label{eqn_5.3}
\Psi=\|\nabla V\|_f^2,\;\;\; P=\exp\left(\kappa
W^{\alpha}\right)\sqrt{W}{\Psi},
\end{equation}
Note that $\Psi$ is a complex version of $\Phi$ in Calabi
geometry. Constants $\kappa$ and $\alpha$ in $P$ are to be
determined in \eqref{eqn_5.3.2} or \eqref{eqn_6.4}. Denote
$$\|V_{,i\bar j}\|_f^2=\sum f^{i\bar j}f^{k\bar l}V_{i\bar l}V_{k\bar j},\;\;
\|V_{,i  j}\|_f^2=\sum f^{i\bar j}f^{k\bar l}V_{,ik}V_{,\bar l\bar j}. $$
 Denote by $\square=\sum f^{i\bar j}\frac{\partial^2}
 {\partial z_i\partial\bar z_j}$ the Laplacian operator.
 Recall that
\begin{equation}\label{eqn_5.1}
\mc S(f)=-\square V=-\sum   f^{i\bar j}V_ {i\bar j}
\end{equation}
is the scalar curvature of $\omega_f.$

\subsection{Differential inequality-I}\label{sect_5.1}
   We are
going to calculate  $\df P$ and derive a differential inequality.
This inequality is similar to that in Theorem \ref{theorem_2.3.1}.
\begin{lemma}[Inequality-I]\label{lemma_5.1.1}
\begin{equation}\label{eqn_5.4}
\frac{\df  P}{P}  \geq   \frac{\|V_{,i\bar{j}}\|_f^2}{2\Psi}
 +\alpha^2\kappa(1 -2 \kappa W^{\alpha}  )W^{\alpha} \Psi
-   \frac{2|\langle \nabla \mathcal S, \nabla V\rangle|}{\Psi}- \left({\alpha} \kappa W^{\alpha}
+\tfrac{1}{2}\right)\mathcal  S,
\end{equation}
where $\langle, \rangle$ denotes the inner product with respect to the metric $\omega_f.$
\end{lemma}
{\bf Proof.} By definition,
$${\Psi}_{,k} = \sum f^{i\bar{j}}\left(V_{,i}V_{,\bar{j}k} +
V_{,ik}V_{,\bar{j}}\right),$$\begin{equation*} \df {\Psi} =\sum
f^{i\bar{j}}f^{k\bar{l}} \left(V_{,i} V_{,\bar{j}k\bar{l}} +
V_{,ik\bar{l}}V_{,\bar{j}} + V_{,ik}V_{,\bar{j}\bar{l}} +
V_{,i\bar{l}}V_{,\bar{j}k}\right).
 \end{equation*}
  By the Ricci
identities
\begin{equation*}
V_{,\bar{j}k\bar{l}} = V_{,k\bar{l}\bar{j}},\;\;\; V_{,ik\bar{l}}
= V_{,k\bar{l}i}+\sum f^{m\bar{h}} V_{,m} R_{k\bar{h}i\bar{l}},
\end{equation*}
we have
\begin{equation}\label{eqn_5.5}
\df {\Psi}  =  \sum
f^{i\bar{j}}f^{k\bar{l}}\left(V_{,ik}V_{,\bar{j}\bar{l}}
 + V_{,i\bar{l}}V_{,\bar{j}k}- V_{i\bar{l}}V_{,k}V_{,\bar{j}}
 \right) -
  2Re (\sum f^{i\bar{j}}V_{,i}\mc S_{,\bar{j}}),
   \end{equation}
  where we use the facts
$ R_{i\bar{j}}= - V_{i\bar{j}}$ and  $\df V=-\mathcal  S . $
Denote $\Pi=\alpha \kappa W ^{\alpha} + \tfrac{1}{2}.
$ Then
\begin{eqnarray}\label{eqn_5.6}
 P_{,i}&=&  P\left(\frac{{\Psi}_{,i}}{{\Psi}} +
\Pi V_{,i}\right) =:   P\Lambda_i, \nonumber\\
\df   P&=&  P \left[\sum f^{i\bar j}\Lambda_i\Lambda_{\bar j}
+\frac{\df{\Psi}}{{\Psi}} -\frac{\|\nabla {\Psi}\|_f^2}{\Psi^2}
+\Pi\df V + \alpha^2\kappa   W ^{\alpha} {\Psi} \right].\;\;\;\;\;\;
\end{eqnarray}
 Choose a new complex  coordinate system      such
that, at $p$,
\begin{equation*}
f_{i\bar{j}}= c\delta_{ij},\;\; V_1 = V_{\bar{1}},\;\;\;V_i =
V_{\bar{i}}=0\;\;\forall \;i>1,\end{equation*} where
$c=[W(p)]^{\frac{1}{n}}.$ Then \eqref{eqn_5.6} can be re-written as
\begin{eqnarray}\label{eqn_5.7}
 \df   P&=&  \frac{P}{\Psi}\left[  \df{\Psi}   + \alpha^2\kappa
  W ^{\alpha}c^{-2}V_1^2V_{\bar 1}^2   + \Pi^2c^{-2}V_1^2V_{\bar 1}^2 -
\Pi\mathcal  S\Psi  \right.\nonumber \\ &&\left.
 +    {2\Pi}{c^{- 2}} \left( Re(V_{,11}V_{,\bar{1}}^2)+
V_{,1\bar{1}}{V_1V_{\bar 1}}\right)\right],
\end{eqnarray}
here and later we denote  the real part by $Re$.
 From
\eqref{eqn_5.5} we have
\begin{eqnarray}\label{eqn_5.8}
\df{\Psi}&=&\sum c^{-2}\left(V_{ik}V_{\bar i\bar k}+   V_{i\bar
k}V_{\bar i k}\right) -c^{-2}V_{1\bar 1}V_1V_{\bar 1}- 2Re (\sum f^{i\bar{j}}V_{,i}\mc S_{,\bar{j}}) \nonumber\\&\geq &
c^{-2}V_{11}V_{\bar 1\bar 1}+c^{-2}V_{1\bar 1}^2+\sum_{i>1}
c^{-2}V_{i\bar k}V_{\bar i k}-c^{-2}V_{1\bar 1}{V_1V_{\bar 1}}
\nonumber\\
&& - 2Re (\sum f^{i\bar{j}}V_{,i}\mc S_{,\bar{j}}).\end{eqnarray} Substituting
\eqref{eqn_5.8} into \eqref{eqn_5.7}, we have
\begin{eqnarray}
 \df P&\geq &\frac{P}{c^2\Psi} \left[\Pi^2(V_1V_{\bar 1})^2 +2\Pi
 Re(V_{,11}V_{\bar 1}^2)+   V_{11}V_{\bar 1\bar 1}\right] \nonumber\\
 && -
\frac{P}{\Psi} \left[2 \left|  \sum f^{i\bar{j}}V_{,i}\mc S_{,\bar{j}} \right|  + \Pi\mc
S\Psi\right]+\frac{P}{c^2\Psi}\sum_{i>1} V_{i\bar j}V_{j\bar i} \nonumber\\
 &&+\frac{
P}{c^2\Psi}[V_{1\bar1}^2+\alpha^2\kappa  W ^{\alpha}(V_1V_{\bar
1})^2 +2({\alpha \kappa  W ^{\alpha}}) V_{1\bar 1}{V_1V_{\bar
1}}]
  \nonumber\\
  &=&\;\frac{P}{c^2\Psi}\left[  I_1+ I_2+\sum_{i>1} V_{i\bar
j}V_{j\bar i}\right]- \frac{P}{\Psi} \left[ \Pi\mc S{\Psi}+2\left|  \sum f^{i\bar{j}}V_{,i}\mc S_{,\bar{j}} \right| \right],\nonumber
\end{eqnarray}
where \begin{eqnarray*} I_1&:=& \left[\Pi^2(V_1V_{\bar 1})^2 +2\Pi
 Re(V_{,11}V_{\bar 1}^2)+   V_{11}V_{\bar 1\bar 1}\right],
 \\
I_2&:=&\left[V_{1\bar1}^2+\alpha^2\kappa  W ^{\alpha}(V_1V_{\bar
1})^2 +2\left({\alpha \kappa  W ^{\alpha}}\right) V_{1\bar
1}{V_1V_{\bar 1}}\right]. \nonumber
\end{eqnarray*}
 It is easy to check that
\begin{eqnarray*}
I_1&=&\left|\Pi V_{1}^2 +V_{11} \right|^2,\\
I_2&=& \alpha^2\kappa ( 1-2\kappa W ^{\alpha}
 ) W^{\alpha}(V_1V_{\bar 1})^2
 +\frac{V^2_{1\bar1}}{2}+ \left|  \frac{1}{\sqrt 2}  {V_{1\bar 1}}+\sqrt{2} {\alpha \kappa
 W ^{\alpha}}V_1V_{\bar 1} \right|^2 .\;\;\;\;\;\;\;\;\;\;
\end{eqnarray*}
The lemma then follows. $\blacksquare$ \v

\subsection{Differential inequality-II}\label{sect_5.2}
Fix a function $ g\in \mc R^\infty(\Omega)$.   Denote by $
\dot{R}_{i{\bar{j}}k\bar{l}}$ and $ \dot{R}_{i\bar{j}}$ the
curvature tensor and  the Ricci curvature of $(\Omega,\omega_g)$,
respectively. Put
$$ \mc {\dot
{R}}:= \sqrt{\sum g^{m\bar{n}}g^{k\bar{l}}g^{i\bar{j}}g^{s\bar{t}}
 \dot{R}_{{m}\bar{l}i\bar{t}} \dot{R}_{k{\bar{n}}s\bar{j}}}.$$
Let $f  \in \mc R^\infty(\Omega)$ and $\phi=f-g.$
  Now we prove a
differential inequality of $ n-\square \phi.$ Obviously $
n-\square \phi=\sum f^{i\bar j}g_{i\bar j}.$
\begin{lemma}[Inequality-II]\label{lemma_5.2.1}
\begin{equation}\label{eqn_5.9}
 \df   \log \left(n-\square \phi\right)  \geq  -\| Ric\|_f - \mathcal {\dot
R}(n-\square \phi).
\end{equation}
\end{lemma}

  To prove this lemma we   calculate $\df\left(n-\square \phi\right)$ firstly.
In the following calculation,
 ";" denotes the covariant derivatives with
respect to the metric $\omega_g$.

\begin{lemma}\label{lemma_5.2.2}
\begin{eqnarray*}
\df \left(n-\square \phi\right)&\geq &\sum
f^{i\bar{j}}f^{n\bar{l}}f^{k\bar{h}}f^{m\bar{p}}g_{k\bar{l}}
\phi_{;n\bar{p}\bar{j}}\phi_{;m\bar{h}i}  - \sum
f^{m\bar{l}}f^{k\bar{h}}g_{k\bar{l}}  V _{m\bar{h}} \\&&- \mc
{\dot {R}}  \left(n-\square \phi\right)^2.
\end{eqnarray*}
\end{lemma}  {\bf Proof.}
A direct calculation gives us $$f_{;i\bar j k}=(g+\phi)_{;i\bar j
k}=\phi_{;i\bar j k},\;\;\;f_{;i\bar j \bar k}=\phi_{;i\bar j \bar
k},\;\;\;f_{;i\bar j s\bar k}=\phi_{;i\bar j s\bar k}.$$ Note that
$$\sum f_{i\bar j}f^{k\bar j}=\delta_{i}^k.$$ Then we have
\begin{equation}\label{eqn_5.10}
f^{k\bar{l}}_{;i}= -\sum
f^{m\bar{l}}f^{k\bar{h}}\phi_{;m\bar{h}i},\;\;\;\;
f^{k\bar{l}}_{;\bar{j}}= -\sum
f^{m\bar{l}}f^{k\bar{h}}\phi_{;m\bar{h}\bar{j}}.
\end{equation}
Then \begin{eqnarray}\label{eqn_5.11}
(n-\df \phi)_{;i}&=&  -\sum
f^{m\bar{l}}f^{k\bar{h}}\phi_{;m\bar{h}i}g_{k\bar{l}} ,\nonumber\\
 \df (n-\square \phi)
&=& \sum f^{i\bar{j}}g_{k\bar{l}}\left(
f^{n\bar{l}}f^{k\bar{h}}f^{m\bar{p}}\phi_{;n\bar{p}\bar{j}}
\phi_{;m\bar{h}i}+
f^{m\bar{l}}f^{n\bar{h}}f^{k\bar{p}}\phi_{;n\bar{p}\bar{j}}
\phi_{;m\bar{h}i} \right.\nonumber \\ && -\left.
f^{m\bar{l}}f^{k\bar{h}}\phi_{;m\bar{h}i\bar{j}}
 \right).
\end{eqnarray}
Now we calculate $\phi_{;m\bar{h}i\bar{j}}$.  Differentiating the  equation $V=\log\det\left( g_{i\bar{j}} +   \phi _{i  \bar j}\right)$  twice
we have
\begin{equation}\label{eqn_5.12}
\sum f^{i\bar{j}}\phi_{;i\bar{j}k\bar{l}} = V_{k  \bar l}+\dot R_{k\bar l}
+ \sum f^{i\bar{j}}f^{n\bar{m}}\phi_{;n\bar{j}\bar{l}}
\phi_{;\bar{m}ik}.
\end{equation}
 By the Ricci identities we
have
\begin{eqnarray}\label{eqn_5.13}
\sum f^{i\bar j}\phi_{;m\bar{h}i\bar{j}} &=&
\sum f^{i\bar j}f_{;i\bar{h}m\bar{j}}\nonumber\\
&=& \sum f^{i\bar j}\left(f_{;i\bar{h}\bar{j}m} + f_{n\bar{h}}
\dot R^n_{im\bar{j}} -
f_{i\bar{n}}  \dot R^{\bar{n}}_{\bar{h}\bar{j}m}\right)\nonumber\\
&=&\sum f^{i\bar j}\left(f_{;i\bar{j}\bar{h}m} + f_{n\bar{h}} \dot
R^n_{im\bar{j}} -
f_{i\bar{n}} \dot R^{\bar{n}}_{\bar{h}\bar{j}m}\right)\nonumber\\
&=& \sum f^{i\bar j}\left(\phi_{;i\bar{j}m\bar{h}}  - f_{n\bar{j}}
\dot R^n_{im\bar{h}}
 + f_{n\bar{h}}\dot R^n_{im\bar{j}}\right)\nonumber\\
 &=& V_{; m\bar{h}} + \sum  f^{i\bar j}f_{n\bar{h}}
 \dot R^n_{im\bar{j}}+\sum f^{i\bar{j}}f^{a\bar{b}}
 \phi_{;a\bar{j}\bar{h}}
\phi_{;\bar{b}im},\;\;\;\;\;
\end{eqnarray}
where we use \eqref{eqn_5.12} in the last step.
 Inserting \eqref{eqn_5.13} into \eqref{eqn_5.11},   we   obtain
  \begin{eqnarray}\label{eqna_5.14}
\df \left(n-\square \phi\right)&=&\sum
f^{i\bar{j}}f^{n\bar{l}}f^{k\bar{h}}f^{m\bar{p}}g_{k\bar{l}}
\phi_{;n\bar{p}\bar{j}}\phi_{;m\bar{h}i}  - \sum
f^{m\bar{l}}f^{k\bar{h}}g_{k\bar{l}}  V _{m\bar{h}} \nonumber\\&&-\sum f^{i\bar j}f^{m\bar l} \dot{R}_{i\bar{l} m \bar {j}} .
\end{eqnarray}
To obtain the lemma we only need to prove
  \begin{equation}\label{eqna_5.15}\sum f^{i\bar j}f^{m\bar l} \dot{R}_{i\bar{l} m \bar {j}}\leq \mc
{\dot {R}}  \left(n-\square \phi\right)^2.
\end{equation}  Note that this inequality is invariant under coordinate transformations. We choose another coordinate system so that
$g_{k\bar{l}}=\delta_{k {l}},\;f_{k\bar{l}}
=\delta_{k {l}}f_{k {k}}. $ Then $$n-\df \phi=\sum
f^{k\bar{k}},\;\;\; |\dot{R}_{i\bar{l} m \bar {j}}|\leq \mc {\dot {R}}.$$
A direct calculation gives us \eqref{eqna_5.15}.
    \hfill$\blacksquare$ \v\n {\bf Proof of Lemma
\ref{lemma_5.2.1}.} We choose another coordinate system so that
$g_{k\bar{l}}=\delta_{k {l}},\;\;\;f_{k\bar{l}}
=\delta_{k {l}}(1+\phi_{k\bar{k}}). $ Then $n-\df \phi=\sum
(1+\phi_{i\bar{i}})^{-1}.$  By the Cauchy inequality one can show that
\begin{equation}\label{eqn_5.14}
\frac{\left|\sum
f^{m\bar{l}}f^{k\bar{h}}g_{k\bar{l}}V_{m\bar{h}}\right|}{n-\square
\phi}=\frac{\sum (1+\phi_{k\bar k})^{-2}V_{k\bar k}}{\sum (1+\phi_{i\bar{i}})^{-1}} \leq \|V_{m\bar{l}}\|_f.
\end{equation}
Now, by using a similar calculation as in \cite{Yau}, we
prove the following inequality:
\begin{equation}\label{eqn_5.15}
\sum f^{i\bar{j}}f^{n\bar{l}}f^{k\bar{h}}f^{m\bar{p}}
g_{k\bar{l}}\phi_{;n\bar{p}\bar{j}} \phi_{;m\bar{h}i}\geq
\frac{\sum f^{i\bar{j}}(\df \phi)_{i}(\df \phi)_{\bar{j}}}{n-\df
\phi}.\end{equation} A direct calculation gives us
\begin{equation*}
(\df \phi)_{i}=(n-\sum f^{k\bar l}g_{k\bar l})_{;i}=\sum f^{k\bar
a}f^{b\bar l}g_{k\bar l}\phi_{;\bar a bi}=\sum (1+\phi_{k\bar
k})^{-2}\phi_{; k\bar ki}\end{equation*} and
\begin{eqnarray*}
&&\frac{\sum f^{i\bar{j}}(\df \phi)_{i}(\df \phi)_{\bar{j}}}{n-\df
\phi}=\nonumber\\
&=& (n-\df \phi)^{-1}\sum_i (1+\phi_{i\bar{i}})^{-1}\left|\sum_k
(1+\phi_{k\bar{k}})^{-2}\phi_{;k\bar{k}i}\right|^2\nonumber\\
&=&(n-\df \phi)^{-1}\sum_i (1+\phi_{i\bar{i}})^{-1}\left|\sum_k
(1+\phi_{k\bar{k}})^{-\frac{3}{2}}\phi_{;k\bar{k}i}
(1+\phi_{k\bar{k}})^{-\frac{1}{2}}\right|^2\nonumber\\
&\leq& (n-\df \phi)^{-1}\left(\sum_{i,k} (1+\phi_{i\bar{i}})^{-1}
(1+\phi_{k\bar{k}})^{-3}\phi_{;k\bar{k}i}\phi_{;k\bar{k}\bar{i}}\right)
\left(\sum_k(1+\phi_{k\bar{k}})^{-1}\right)\nonumber\\
&=&\sum_{i,k} (1+\phi_{i\bar{i}})^{-1}
(1+\phi_{k\bar{k}})^{-3}\phi_{;k\bar{k}i}\phi_{;k\bar{k}\bar{i}}\nonumber\\
&\leq& \sum_{i,k,l} (1+\phi_{i\bar{i}})^{-1}
(1+\phi_{k\bar{k}})^{-1}
(1+\phi_{l\bar{l}})^{-2}\phi_{;k\bar{l}i}\phi_{;l\bar{k}\bar{i}}.\nonumber\\
&=&\sum f^{i\bar{j}}f^{n\bar{l}}f^{k\bar{h}}f^{m\bar{p}}
g_{k\bar{l}}\phi_{;n\bar{p}\bar{j}} \phi_{;m\bar{h}i}\nonumber
\end{eqnarray*}
Thus \eqref{eqn_5.15} is proved. Note that
\begin{equation*}
\df \log (n-\df \phi)=\frac{\df (n-\df \phi)}{n-\df
\phi}-\frac{\sum f^{i\bar j} (n-\df \phi)_{i}(n-\df \phi)_{\bar
j}}{(n-\df \phi)^2}.
\end{equation*}
 Using  \eqref{eqn_5.14}, \eqref{eqn_5.15} and Lemma \ref{lemma_5.2.2},
  we obtain Lemma \ref{lemma_5.2.1}.   $\blacksquare$\\

\subsection{Differential inequality-III}\label{sect_5.3}

In this subsection, we apply  Differential Inequality-I on the toric
surfaces. We use notations introduced in \S\ref{sect_1.1} and
\S\ref{sect_1.2}. Let $0\in \Omega\subset \cplane^2_{\vartheta}$,
$f=f_\vartheta$ be the potential
 function
on $\Omega$. Let $T=\sum f^{i\bar i}$.
 Put
\begin{equation*}
  Q=e^{N_1(|z|^2-A)}\sqrt{W}T,\;\;\;\;
\end{equation*} where $A,N_1$ are constants.

\begin{lemma}\label{lemma_5.3.1}
Let $K= \mc S(f)\circ \tau_f\inv$ be the scalar curvature function
on $\t$. Suppose that
  \begin{equation}\label{eqn_5.3.1}
  \max_{\tau_{f}(\bar \Omega)}\left(|  K|+ \sum\left|\frac{\p
K}{\p \xi_i}\right|\right)\leq \mff N_2,\;\;\;\;\;\max_{\bar \Omega} W\leq\mff
N_2,\;\;\;\;\;\; \max_{\bar \Omega}|z|\leq \mff N_2\end{equation}  for some constant
$\mff N_2>0.$ Then we may choose
\begin{equation}\label{eqn_5.3.2}
A=\mff N_2^2+1, N_1=100,\alpha=\frac{1}{3},\kappa=[4\mff N_2^\frac{1}{3}]\inv
\end{equation}
 such that
\begin{equation}\label{eqn_5.3.3}
\df (P+Q+ \mff C_7f ) \geq \mff C_{8}(P+Q)^2>0
\end{equation} for some positive constants $\mff C_7$ and $\mff C_8$ that
depend only on $\mff N_2$ and $n$.
\end{lemma}
{\bf Proof.} Applying  Lemma \ref{lemma_5.1.1} and the choice of $\alpha$ and $\kappa$, in particular,
$\kappa W^\alpha\leq 1/4$, we have
\def \al{\frac{1}{3}}
\begin{equation}\label{eqn_5.3.4}
\frac{\Psi\df  P}{P} \geq \left(\frac{1}{2}\|V_{,i\bar{j}}\|_f^2
 +\frac{1}{18} \kappa W^{\al}\Psi^2\right)
- \left(2 |\langle \nabla \mathcal S, \nabla V\rangle| + \Psi
|\mathcal  S| \right).
\end{equation}
Treatment for $\langle \nabla \mathcal S, \nabla \log W\rangle$:
using log-affine coordinates we have
\begin{equation*}
|\langle \nabla \mathcal S, \nabla V\rangle| =
\left|\sum f^{ij}\frac{\partial \mc S}{\partial x_i}\frac{\partial
V}{\partial x_j} \right| =\left|\sum f^{ij}\frac{\partial K}{\partial
\xi_k}\frac{\partial \xi_k} {\partial x_i}\frac{\partial V}{\partial
x_j} \right| \leq \mff N_2\sum_{j}\left|\frac{\partial V}{\partial x_j}
\right|,
\end{equation*}
where we use the fact $\frac{\partial \xi_k} {\partial x_i}=f_{ki}$;
if we use the complex coordinates $z_i$, we have
$$
\left|\frac{\partial V}{\partial x_j} \right|=   \left|z_j \frac{\partial V}{\partial z_j}\right|.
$$
Since $|z|$ is bounded, we conclude that
\begin{equation*}
|\langle \nabla \mathcal S, \nabla V\rangle| \leq C\sum_j|V_{\bar
j}|\leq C\sqrt{2WT\Psi}.
\end{equation*}
We explain the last step: suppose that $0<\nu_1\leq\nu_2$ are the
eigenvalues of $(f_{i\b j})$, then
$$
\Psi=\sum f^{i\bar j}V_iV_{\bar j} \geq \nu_2\inv
(|V_1|^2+|V_2|^2)\geq(WT)\inv(|V_1|^2+|V_2|^2).
$$
Note that $$(eW)^{-\half}\leq \Psi/P=(\exp(\kappa W^\alpha)
W^\half)\inv \leq W^{-\half},$$ \eqref{eqn_5.3.4} is then
transformed to be
\begin{equation}\label{eqn_5.3.5}
{\df  P} \geq W^\half\left(\frac{1}{2}\|V_{,i\bar{j}}\|_f^2
 +\frac{1}{18} \kappa W^{\al}\Psi^2\right)
- C' W^\half\left(\sqrt{WT\Psi} + \Psi |\mathcal  S| \right).
\end{equation}
Applying the Young inequality and the Schwartz inequality to  terms
in \eqref{eqn_5.3.5},
  we  have that
\begin{equation}\label{eqn_5.3.7}
\df  P \geq \half W^\half \|V_{,i\bar j}\|_f^2 + C_1W^{-\frac{1}{6}}
P^2  -\epsilon QT -C_2(\epsilon).
\end{equation}
For example,
$$
 W (T\Psi)^{\half}\leq
\delta (W^{\frac{1}{8}}T^\half)^4+\delta
(W^{\frac{5}{24}}\Psi^\half)^4 +C_\delta W^{\frac{4}{3}}\leq C_2'\delta
(QT+ W^{-\frac{1}{6}}P^2) +C_\delta'.
$$
We skip the proof of \eqref{eqn_5.3.7}.

 By a direct calculation we have
\begin{equation*}\df Q \geq Q\left( N_1 T+\half \square V +
 \df \log T
\right) \geq   Q \left(N_1T-\half\mc
S-\left\|V_{,i\bar{j}}\right\|_f\right),
\end{equation*}
where we use the formula \eqref{eqn_5.9}   with $g_{i\bar j}=\delta_{ij}$ to calculate  $\square
\log T$. Here $\left\|V_{,i\bar{j}}\right\|_f=\|Ric\|_{f}.$
 Using the explicit value $A,N_1$ and the bounds of $W$ and $\mc S$,
 applying the Schwartz inequality properly,
we can get
\begin{equation}\label{eqn_5.3.6}
\df Q\geq   -\frac{1}{4} W^\half\|V_{,i\bar j}\|_f^2 + \frac{N_1}{3}
QT -C_3(N_1,\mff N_2).\end{equation}
Combining \eqref{eqn_5.3.7} and \eqref{eqn_5.3.6}, and choosing $\epsilon=\frac{1}{100}$, we have
\begin{equation*}
 {\df (  Q+  P)} \geq  {C_1W^{-\six} P^2+\frac{N_1}{4}
QT}  - {C_4} .
\end{equation*}
 Note that
$$
T= e^{-N_1(|z|^2-A)}W^{-\half}  Q\geq
e^{-N_1(|z|^2-A)}\mff N_2^{-\al}W^{-\six} Q \geq C_5W^{-\six}  Q,
$$
 we get
\begin{equation*}
{\df (  Q+  P)}\geq C_6W^{-\six}(  Q+  P)^2 -C_7
\end{equation*}
where  $C_6,C_7$ are constants depending only on $C_1,N_1$ and $\mff
N_2$. Our lemma follows from $\Box  f=n$ and $|W|\leq \mff N_2$.
$\;\;\;\; \blacksquare$


\section{Complex interior estimates and regularities}\label{sect_6}

We apply the differential inequalities to derive interior
estimates in terms of the norm $\mc K$ of Ricci tensor  (cf.\;\eqref{eqn_1.1c}).

\subsection{Interior estimate of $\Psi$}\label{sect_6.1}

We use Lemma \ref{lemma_5.1.1} to
derive the interior  estimate of   $\Psi$ in a geodesic ball.
\begin{lemma}\label{lemma_7.16a} Let $ f  \in \mc R^\infty(\Omega)$
 and $B_{a}(o)\subset \Omega$ be a closed
geodesic ball of radius $a$ centered at o. Set $  W_\diamond:=\max\limits_{B_{a}(o)} W.$  Suppose that
\begin{equation} \label{eqn_7.31a}
\min\limits_{B_{a}(o)}|\mathcal S|\neq0,
\mbox{ and, }  W^{\half} (\mc K+\|\nabla \log |\mc S|\|^2_f+\Psi)\leq 4,
\end{equation} in $B_{a}(o)$. Then the following estimate  holds in
$\; B_{a/2}(o)\;$
\begin{eqnarray}\label{eqn_7.32a} {W^{\half}}  \Psi
 \leq \mff C_{6}\left[   W_\diamond ^{\half}\max_{B_{a}(o)}|\mc S| +   W_\diamond ^{\frac{1}{3}}\max_{B_{a}(o)}|\mc S|^{\frac{2}{3}}  +
a\inv W_\diamond^{\frac{1}{4}}+ a^{-2} W_\diamond ^{\half} \right],\;\;
   \;
 \end{eqnarray}    where $\mff C_6$ is a constant  depending
only on $n$.
 \end{lemma}
\v \n {\bf Proof.} Consider the function
$$F:= ( a^2-r^2)^2 {P} $$
defined in   $B_{a}(o, \omega_f),$ where $r$ denotes the geodesic distance from $o$ to $z$ with
respect to the metric $ \omega_f.$  $F$ attains its
supremum at some interior point $q^\ast$.
 Then, at $q^\ast$,
\begin{eqnarray}\label{eqn_7.33a}
&& \frac{P_{,i}}{P} -\frac{2(r^2)_{,i}}{{a}^2-r^2}=0,
\\\label{eqn_7.34a}
 &&\frac{\df  P}{P}- \frac{\|\nabla  {P} \|_f^2}{P^2} -
\frac{2   \|\nabla( r^2)  \|_f^2}{(a^2-r^2)^{2}}
 -\frac{4(1+r\df r)}{a^2-r^2}\leq 0.
\end{eqnarray}
By choosing
\begin{equation}\label{eqn_6.4}
\kappa= \frac{1}{4 {W_\diamond}^{\frac{1}{4}}},\;\;\;\;\alpha=\frac{1}{4}
\end{equation}
 in Inequality-I (Lemma \ref{lemma_5.1.1}), we have

\begin{equation*}
\frac{\df  P}{P}  \geq \epsilon_o\left(\frac{W}{W_\diamond}\right)^{\tfrac{1}{4}}\Psi
-   \frac{2|\langle \nabla \mathcal S, \nabla V\rangle|}{\Psi}- \left[\frac{1}{16}
\left(\frac{W}{W_\diamond}\right)^{\tfrac{1}{4}}+\frac{1}{2}\right]|\mathcal  S|,
\end{equation*}
 where $\epsilon_{o}=\frac{1}{128}.$
  Inserting this and \eqref{eqn_7.33a} into \eqref{eqn_7.34a} , we have, at $q^\ast$
\begin{equation}\label{eqn_7.35a}
 \epsilon_{o}\left(\frac{W}{W_\diamond}\right)^{\tfrac{1}{4}} {\Psi}
-   \frac{2\|\nabla  {\mathcal S}\|_{f}}{\sqrt{\Psi}}- \frac{9}{16} {|\mathcal  S|}
 -\frac{24 {a}^2}{(a^2-r^2)^{2}}
 -\frac{4(1+r\df r)}{a^2-r^2}\leq 0.\end{equation}
Now we estimate $r\df r$. Denote by $\Gamma$ the geodesic from $o$ to $q^\ast$. We consider two cases:
\v\n
{\bf Case 1.} Along $\Gamma$ the following estimate holds
$$W\geq \frac{1}{100}W(q^\ast).$$
Then   we have
$$\| {Ric}\|(q)\leq \mathcal K (q)\leq  \frac{4}{\sqrt{W(q)}} \leq \frac{40}{\sqrt{W(q^\ast)}},\;\;\;\;\forall q\in \Gamma.$$
By the Laplacian
comparison theorem (see Remark \ref{remark_3.3.2}) we have
\begin{equation}\label{eqn_7.36a}
r\df r\leq C(n)\left[1+\frac{\sqrt{40} {a}}{(W(q^\ast))^{\frac{1}{4}}}\right].\end{equation}
Inserting \eqref{eqn_7.36a} into \eqref{eqn_7.35a}, we have, at $q^\ast,$
 $$\epsilon_{o}\left(\frac{W}{W_\diamond}\right)^{\tfrac{1}{4}} {\Psi}
-   \frac{2\|\nabla  {\mathcal S}\|}{\sqrt{\Psi}}-\frac{9}{16}  {|\mathcal  S|}-\frac{24 a^2}{(a^2-r^2)^{2}}
-\frac{4+4C(n)}{a^2-r^2}
 -\frac{4C(n)\sqrt{40} {a}}{(a^2-r^2) W ^{\frac{1}{4}}}\leq 0. $$
We rewrite it as
\begin{equation}
 {\Psi}\leq b_1A + \frac{b_2B}{\sqrt {\Psi}},\end{equation}
where $b_1$, $b_2$ are constants depending only on $n$ and
$$A:=\left(\frac{W_\diamond}{W}\right)^{\tfrac{1}{4}}\left[\frac{a^2}{(a^2-r^2)^{2}}+ \frac{1}{a^2-r^2}+\frac{a}{(a^2-r^2)W^{\frac{1}{4}}}+ {|\mathcal  S|}\right](q^\ast),$$
$$B:=\|\nabla  {\mathcal S}\|_f\left(\frac{W_\diamond}{W}\right)^{\tfrac{1}{4}}(q^\ast).$$
Applying Young's inequality and multiplying $\sqrt{W} $ on the both sides, we get
\begin{equation}
\sqrt{W} {\Psi}\leq b_3\sqrt{W}A + b_4 \sqrt{W}B^{2/3}.\end{equation}
 Note that, at $q^\ast$
$$B{W} ^{\tfrac{3}{4}}=\|\nabla  {\mathcal S}\|_f\left(\frac{W_\diamond}{W}\right)^{\tfrac{1}{4}}{W} ^{\tfrac{3}{4}}=\|\nabla \log {|\mathcal S|}\|_f W^{\tfrac{1}{4}} {|\mathcal S|} ( {W_\diamond}{W})^{\tfrac{1}{4}} \\
  \leq    2  {|\mathcal S|}{W_\diamond}^{\tfrac{1}{2}}
$$
It follows that, at $q^\ast$,
$$
 F(q^\ast) \leq b_5\left[(W_\diamond)^{\frac{1}{2}} {a}^2 + (W_\diamond)^{\frac{1}{4}} {a}^3
 +  \max_{B_ {a}}
\left(| {\mathcal S}|(W_\diamond)^{\frac{1}{2}} +   |{\mathcal S}|^{\tfrac{2}{3}}(W_\diamond)^{\frac{1}{3}}\right){a}^4 \right] $$
for some constants $b_5>0$  depends only on $n$.
Since $F\leq F(q^\ast),$ we conclude that \eqref{eqn_7.32a} holds in $B_{a/2}(o)$.

\v\n
{\bf Case 2.} There is a point $q\in \Gamma$ such that $W(q)<\frac{1}{100}W(q^\ast)$. Then there is a point $q_1\in \Gamma$ such that
$$W(q_1)=\frac{1}{100}W(q^\ast);$$
$$W(q)\geq \frac{1}{100}W(q^\ast),\forall q\in [q_1,q^\ast].$$
As $\sqrt{W} {\Psi}\leq   4$£¬ we have
$$\left| \frac{dW^{\frac{1}{4}}}{ds} \right|\leq \frac{1}{2}.$$
Denote by $r_{q_1}(z)$ the geodesic distance function from $q_1$ to $z$ with respect to the metric $\omega_{f}$.  It follows that the geodesic distance between $q_1$ and $q^\ast$ satisfies
\begin{equation}r_{q_1}(q^\ast)\geq 2(W^{\frac{1}{4}}(q^\ast)-W^{\frac{1}{4}}(q_1))= 2W^{\frac{1}{4}}(q^\ast)\left(1-\frac{1}{\sqrt{10}}\right)>\frac{4}{3}W^{\frac{1}{4}}(q^\ast).
\end{equation}
We consider the "support function" $\tilde  F$ of $F$ (see \cite{Y-S}),
$$\tilde F  (z)=\left[ {a}^2-(r(o,q_1)+ r_{q_1}(z))^2\right]^2 {P}.$$
Set $\tilde r(q) =r(o,q_1)+ r_{q_1}(q)$. {We claim} that
\begin{enumerate}\item[(1)] $\tilde F(q^\ast)=F(q^\ast)$,
\item[(2)] in $  B_{a-r(o,q_1)}(q_1),$ $\tilde F$ attains its maximum at $q^\ast$.
\end{enumerate}
{\em Proof of Claim}. (1) is obvious. We prove (2). By the  triangle inequality
 $$r(q)\leq r(o,q_1)+ r_{q_1}(q)\leq a,\;\;\;\;\forall   q\in B_{a-r(o,q_1)}(q_1),$$
 we have, for any point $q\in B_{a-r(o,q_1)}(q_1)$,
$$\left[ {a}^2-(r(o,q_1)+ r_{q_1}(q))^2\right]^2 {P}(q)\leq ( {a}^2-r^2)^2 {P}(q)\leq F(q^\ast)=\tilde F (q^\ast). \;\;$$
The claim follows.

When restricting on $[q_1,q^\ast]$ of $\Gamma$ we have
$$W\geq \frac{1}{100}W(q^\ast)$$
and $\tilde r(q)=r(q).$
Obviously, $q^\ast$ is an interior point in $  B_{a-r(o,q_1)}(q_1).$ Using the maximum principle,  by the same calculation as above, we get, at $q^\ast$,
\begin{equation}
\epsilon_{o}\left(\frac{W}{W_\diamond}\right)^{\tfrac{1}{4}} {\Psi}
-   \frac{2\|\nabla  {\mathcal S}\|}{\sqrt{\Psi}}- \frac{9}{16} {|\mathcal  S|}
 -\frac{24 {a}^2}{( {a}^2-\tilde r^2)^{2}}
 -\frac{4(1+\tilde r \df \tilde r )}{a^2-\tilde r^2}\leq 0.\end{equation}
Note that $\tilde r \df \tilde r =(r(o,q_1)+ r_{q_1}(q))\df r_{q_1}(q)$ and $\tilde r \leq a.$ Applying the Laplacian comparison theorem to $ \df r_{q_1}(q)$ we obtain that
$$\tilde r \df \tilde r=\tilde r\df r_{q_1}(q^\ast) \leq a(n-1)\left[ \frac{1}{r_{q_1}(q^\ast)}+\frac{C}{W^{\frac{1}{4}}(q^\ast)}\right]\leq \frac{C'(n) a}{W^{\frac{1}{4}}(q^\ast)} $$
for some constant $C'(n)>0$, where we used $r_{q_1}(q^\ast)\geq \frac{4}{3}W^{\frac{1}{4}}(q^\ast).$
By the same argument as that in {\bf Case 1} we get
$$ \tilde F  (q^\ast)\leq b'_5\left[(W_\diamond)^{\frac{1}{2}} {a}^2 + (W_\diamond)^{\frac{1}{4}} {a}^3+| {\mathcal S}|(W_\diamond)^{\frac{1}{2}}{a}^4 + \left(\|  {\mathcal S}|(W_\diamond)^{\frac{1}{2}}\right)^{2/3}{a}^4\right]$$
Hence $F\leq F(q^\ast)=\tilde F  (q^\ast).$

  Note that
$$
(a^2-r^2)^2P=(a^2-r^2)^2\exp(\kappa W^\alpha) W^{\half}\Psi
\geq (a^2-r^2)^2W^\half\Psi.
$$
For both cases we have estimates \eqref{eqn_7.32a} in $B_{a/2}(o)$. $\blacksquare$

\begin{remark}
[The Laplacian comparison theorem]\label{remark_3.3.2} Let $M^n$ be a
complete Riemannian manifold. Assume that $\Gamma:[0,a]\to M$ is a minimal geodesic parametrized by arc length with $\Gamma(0)=p$. Let $r(q)$ be the distance function from $p$ to $q.$   Suppose that  $$Ric(\Gamma(t))\geq -(n-1)C_0^2,\;\;\;\;\;\;\mbox{ for any } t\in [0,a],$$ where
$C_0>0 $ is a constant. Then for any $q\in \Gamma,$
\[ r\df  r(q)\leq C(n)(1+ {C_0}r(q)).\]
\end{remark}
The proof is similar to the standard Laplacian comparison theorem by using the Jacobi fields.

\v By the same  argument  as in Lemma \ref{lemma_7.16a}, applying Laplacian comparison theorem and  Young's inequality  to \eqref{eqn_7.35a}, we have
\begin{lemma}\label{lemma_6.1.2} Let $ f  \in \mc R^\infty(\Omega)$ and $B_{a}(o)\subset \Omega$ be the geodesic ball of radius $a$ centered at $o$. Suppose that there are constants $N_1,N_2>0$ such that
\begin{equation*}
 \mc K\leq
N_1,\;\;\;\;\;  W\leq N_2,
\end{equation*} in $B_{a}(o)$. Then
in $ \;   B_{{a}/{2}}(o)$
\begin{eqnarray*} W^{\half}   \Psi
&\leq&\mff C_9 N_2^\half\left[\max_{B_{a}(o)}\left(|\mc S| + \|\nabla \mc S\|_f
^{\frac{2}{3}}\right) + a\inv + a^{-2} \right].
 \end{eqnarray*}
    where $\mff C_9$ is a constant depending
only on  $n$ and $N_1.$
 \end{lemma}

\begin{remark}\label{remark_6.1.3}
Recall that $\Psi=\|\nabla\log W\|_f^2$. Hence
\begin{equation*}
W^{\half}\Psi=16 \|\nabla W^{\frac{1}{4}}\|_f^2.
\end{equation*}
Therefore, the result in Lemma \ref{lemma_7.16a}  and Lemma \ref{lemma_6.1.2} can be treated as a
bound of  $\nabla W^{\frac{1}{4}}$.
\end{remark}

\subsection{Interior  estimate  of  $\|\nabla f\|_f $}\label{sect_6.2}
\begin{lemma}\label{lemma_6.2.1}
Let $ f\in \mc R^\infty(\Omega) $  with
$f(p_0)=\inf_{\Omega} f=0$. Suppose that for $a>1,$
 \begin{equation}\label{eqn_6.110}
\;\; \mc K\leq N_0,
 \;\;in \;\; B_a(p_0 ).
\end{equation} Then in $  B_{{a}/{2}}(p_0 )$
\begin{equation}\label{eqn_6.10}
 \frac{\|\nabla f\|_f^2}{(1+ f)^2}\leq \mff C_{10}
\end{equation}
where  $\mff C_{10}>0$ is a constant depending only on $ n$ and
$N_0.$ Then,  for any $q\in B_{a/2}(p_0),$
\begin{equation}\label{eqn_6.11}
   f(  q)-  f(  p_0) \leq \exp(\sqrt {\mff C_{10}}a).
   \end{equation}
\end{lemma}
{\bf Proof.}  Consider the function
$$
F=(a^2-r^2)^2\frac{\sum f^{i\b j}  f_{i}  f_{\b j}}{(1+f)^2},$$ in
$B_a(p_0)$.  $F$ attains its supremum at some interior point
$p^\ast$.
Then, at $p^\ast$,
\begin{eqnarray} \label{eqn_6.12a} &&
-\left[\frac{2(r^2)_{,k}}{a^2-r^2}+2\frac{f_{,k}}{1+f}\right]\sum
f^{i\b j}f_{i}f_{\b j} +\sum f^{i\b j}
f_{,i  k}f_{\b j}+f_k=0, \\
&& -\left[\frac{2\|\nabla r^2\|_f^2}{(a^2-r^2)^2}+\frac{2\df
(r^2)}{a^2-r^2}+\frac{2\df f}{1+f} \right]\|\nabla f\|_f^2 +n\nonumber\\
&&+\sum
f^{i\bar j}f^{k \bar l}f_{,ki}f_{,\b l\b j}+\sum f^{i\bar j}f^{k
\bar l}f_{,ik\b l}f_{\b j} \nonumber\\
\label{eqn_6.13a} &&
+\frac{2\|\nabla f\|_f^4}{(1+f)^2}
-\sum\left[\frac{2(r^2)_{,k}}{a^2-r^2}+2\frac{f_{,k}}{1+f}\right]
f^{k \bar l}f^{i \bar j}( f_{i \b l}f_{\b j}+
   f_{i}f_{,\b l\b j}) \leq 0 .\end{eqnarray}
 Choose the coordinate system such that at $p^\ast,$ we
have
$$f_{i\b j}=\delta_{ij},\;\;  f_{,1}=f_{,\bar 1}=\|\nabla f\|_f,\;\; f_{,i}=f_{,\bar i}=0,\;\;\mbox{ for } i\geq2.$$
Then \eqref{eqn_6.12a} and \eqref{eqn_6.13a} can be read as
\begin{eqnarray} \label{eqn_6.112} && \;\;\;\;\;-\left[\frac{2(r^2)_{,k}}{a^2-r^2}+2\frac{f_{,k}}{1+f}\right]
 f_{1}f_{\b 1} +f_{,1  k}f_{\b 1}+f_1\delta_{1k}=0, \\
&& -\left[\frac{2\|\nabla r^2\|_f^2}{(a^2-r^2)^2}+\frac{2\df
(r^2)}{a^2-r^2}+\frac{2\df f}{1+f} \right]f_{1}f_{\b 1} +n+\sum
f_{,lk}f_{,\b l\b k}+\sum f_{,1k\b k}f_{\b 1}
\nonumber\\\label{eqn_6.113} && \;\;\;\;\;+\frac{2( f_{1}f_{\b
1})^2}{(1+f)^2}
-\sum\left[\frac{2(r^2)_{,k}}{a^2-r^2}+2\frac{f_{,k}}{1+f}\right](
\delta_{1  k}f_{\b 1}+
   f_{1}f_{,\b 1\b k}) \leq 0 .\end{eqnarray}
 Applying the Ricci
identities, the fact $\sum r_kr_{\b k}=\frac{1}{4},$  and the Laplacian
comparison theorem  to \eqref{eqn_6.113},  and  inserting
 \eqref{eqn_6.110}, then we have
\begin{eqnarray} &&-\left[\frac{3a^2}{(a^2-r^2)^2}+\frac{
C_1(n)(1+\sqrt{N_0}a)}{a^2-r^2}+\frac{2n}{1+f}-N_0 \right]f_{1}f_{\b
1} +n\nonumber\\&&
+\frac{2( f_{1}f_{\b 1})^2}{(1+f)^2}
 +\sum  f_{,lk}f_{,\b l\b k}-\left|\sum\left[\frac{2(r^2)_{,k}}{a^2-r^2}
+\frac{2f_{,k}}{1+f}\right]
  f_{1}f_{,\b 1\b k}\right|
 \nonumber\\&&
\label{eqna_6.114}  -\left| \frac{2(r^2)_{,1}f_{\b 1}}{a^2-r^2}
+\frac{2f_{,1}f_{\b 1}}{1+f} \right|  \leq  0.
\end{eqnarray}
Note that
\begin{eqnarray*}2\left|\sum \left[\frac{(r^2)_{,k}}{a^2-r^2}+\frac{f_{,1}\delta_{1k}}{1+f}\right]
  f_{1}f_{,\b 1\b k}\right|  \leq
   \sum \left|\frac{(r^2)_{,k}}{a^2-r^2}+\frac{f_{,1}\delta_{1k}}{1+f}\right|^2f_{,1}f_{,\bar 1}+\sum f_{,\b 1\b k}f_{, 1 k} &&\\
    =  \frac{\|\nabla r^2\|^2_ff_{,1}f_{,\bar 1}}{(a^2-r^2)^2}+2 Re\left( \frac{(r^2)_{,1}}{a^2-r^2}\frac{f_{,1}^2f_{,\bar 1}}{1+f}\right)
   +\frac{(f_{,1}f_{,\bar 1})^2}{(1+f)^2} +\sum f_{,1k}f_{,\b1\b k}. && \end{eqnarray*}
  Inserting the above inequality into \eqref{eqna_6.114}  we have
\begin{eqnarray}-\left[\frac{C_2a^2}{(a^2-r^2)^2}+\frac{
{C_3}a}{a^2-r^2}+\frac{2n}{1+f} \right]f_{1}f_{\b 1}+\frac{(
f_{1}f_{\b 1})^2}{(1+f)^2} \nonumber\\\label{eqn_6.114}
-N_0f_{1}f_{\b 1}-2 \left|\frac{(r^2)_1f_{\b
1}}{a^2-r^2}\right|-\frac{2f_1 f_{\b
1}}{1+f}- 2Re\left[\frac{(r^2)_{,1}}{a^2-r^2}\frac{f_{,1}^2f_{,\bar 1}}{1+f}\right]\leq 0.\end{eqnarray}
 Multiplying $\frac{(a^2-r^2)^4}{(1+f)^2}$ on the both sides of  \eqref{eqn_6.114}
   and applying Schwarz's inequality,
   we obtain the following inequality
$$F\leq C_6a^4+C_7(a^2+a^3).\;\;\; $$
Then in $B_{\frac{a}{2}}(p_0)$ we obtain \eqref{eqn_6.10}.\;\;\;
$\blacksquare$

\subsection{Interior estimate of eigenvalues of $(f_{i\bar j})$}\label{sect_6.3}
For simplicity, we assume that $g=\sum z_i z_{\bar i}.$ Denote $T=\sum f^{i\bar i}.$ Then Inequality-II (Lemma \ref{lemma_5.2.1}) can be re-written as
\begin{equation}\label{eqn_6.12}
\df \log T   \geq -\|Ric\|_f  .
\end{equation}
We apply this to prove the following lemma.
\begin{lemma}\label{lemma_6.3.1}Let $ f  \in \mc R^\infty(\Omega)$
and $B_a(p_0)\subset \Omega.$  Suppose
\begin{equation*}W\leq \mff N_4,\;\;\;  \mc K \leq \mff N_4,
\;\;\; |z|\leq
\mff N_4.\end{equation*}in $B_a(p_0)$, for some constant $\mff N_4>0$. Then
there exists a constant $\mff C_{11}>1$ such that
$$\mff C_{11}\inv\leq \lambda_1\leq \cdots
\leq\lambda_n\leq \mff C_{11},\;\;\forall\;q\in B_{a/2}(p_0).$$
 where $\lambda_1,\cdots,\lambda_n$ are eigenvalues of the matrix $(f_{i\b
 j}),$ $\mff C_{11}$
 is a positive constant depending on $n,a$ and $\mff N_4.$
  \end{lemma}
{\bf Proof.} Consider the function
$$F:=(a^2-r^2)^2e^{|z|^2}T$$ in $B_a(p_0 ).$
  $F$ attains its supremum at some interior point $p^*$.   Then, at $p^*$, we have
\begin{eqnarray}\label{eqn_6.13}
&& z_{\bar i}+   (\log T)_{,i}    - \frac{2(r^2)_{,i}}{a^2-r^2}=0 ,\nonumber\\
&& T+ \square \log T   - \frac{2\|\nabla
r^2\|_f^2}{(a^2-r^2)^2}-\frac{2\square (r^2)}{a^2-r^2}\leq
0.
\end{eqnarray}   Note that  $\mc K \leq \mff N_4$.
Thus, by \eqref{eqn_6.12} and
the Laplacian comparison theorem, we have
$$T \leq
\frac{b_1a^2}{(a^2-r^2)^2}+\frac{b_2a}{a^2-r^2}+\mff N_4,$$ for some
positive constants $b_1,$ $b_2$ depending only on $n$ and $\mff N_4$.
Then
\begin{equation*} F\leq C_3a^3+C_4a^4,\end{equation*}
for some positive constants $C_3$ and $C_4$  that  depends only on
$n$ and $\mff N_4$. Then there is a constant $C_5>0$ such that
$$\lambda_{1}\inv\leq
\sum f^{i\bar i}\leq C_5\;\;\;\; in \;\;\;
B_{\frac{a}{2}}(p_0,\omega_f).$$
Since $W$ is bounded and $\lambda_1$ is bounded below, hence $\lambda_n$ is bounded above.
$\blacksquare$

\subsection{Bootstrapping}\label{sect_6.4}
In the subsection, for  the sake of  convenience, we state and prove our results only on the toric surfaces, although the arguments work for higher dimensional toric manifolds.
Suppose that $ \mc S(u)=K $  for some $K\in C^\infty(\bar\Delta)$.
Then in terms of the complex coordinates, the equation becomes
$\mc S(f)=\tilde K$. For example, on the complex
torus in terms of log-affine coordinates:
$$
\tilde K= K\circ \nabla^f.
$$
The bootstrapping argument says that  the regularity
of $f$ can be obtained via $\|\tilde K\|_{C^k}$ if
$\|f\|_{C^2}$ is bounded. However, it is not so obvious that
this is still true in terms of $\|K\|_{C^k}$.
In this subsection, we verify this by direct computations.
The main theorem in this subsection is
\begin{theorem}\label{theorem_6.4.1}
Let $\CHART$ to be one
of the coordinate charts
$\CHART_\Delta,\CHART_\ell$ and $\CHART_\vartheta$.
Suppose that $U\subset \CHART$ is a bounded open set  and
 \begin{equation}\label{eqna_6.4}
 C_1\inv \leq \chi_1 \leq   \chi_2 \leq
C_1\end{equation} for some constant $C_1>0,$
where
$\chi_1,\chi_2$ are the
eigenvalues of the matrix $(\sum g^{i \b j}f_{k\b j}).$
  Then for any $U'\subset U$ and any $k\geq 0$
$$
\|f\|_{C^{k+3,\alpha}(U')}\leq\mff  C_k,
$$
where $\mff C_k$ depends on $\|K\|_{C^k}$,
$d_E(U',\partial U) $ and the bound of $U$.  Recall that   $g$ is the potential
function of the Guillemin metric.
\end{theorem}
The rest of  this subsection is devoted to the proof of the theorem.

\v\n
{\em Case 1, bootstrapping on $\CHART_\Delta$.} The interior regularity of $u$ in $\Delta$
is equivalent to that of $f$ in the complex
torus. Hence we may apply the bootstrapping argument
for $u$ to conclude the regularity of $f$.

\v\n
{\em Case 2, bootstrapping on $\CHART_\ell$.} Without loss of generality, we assume that
$\CHART_\ell=\CHART_{\halfplane}$.
One of the main  issues is to study the
relation between the derivatives
of $\mc S(f)$ and $\mc S(u)$.
Recall that the complex  coordinate is $(z_1,w_2)$.
By an explicit computation, we have
$$
\tilde K(z_1,w_2)=K(\xi_1,\xi_2),
$$
where
$$
\xi_1= {z_1} \frac{\partial f}{\partial z_1},\;\;\;
\xi_2=2\frac{\partial f}{\partial w_2}.
$$
Then by a direct computation, we have
\begin{lemma}\label{lemma_6.4.2}
For any $k$ there exists a constant $C'_k$ such that
$$
\|\tilde K\|_{C^k(U)}\leq C_k',
$$
where $C'_k$ depends only on $\|f\|_{C^{k+1}(U)}$, $|K|_{C^k}$ and
the bound of $U$.
\end{lemma}
{\bf Proof.} When $k=1$, then
$$
\frac{\partial\tilde K}{\partial z_1} = \frac{\partial
K}{\partial \xi_1} \left(\frac{\partial f}{\partial z_1}+
z_1\frac{\partial^2 f}{\partial^2 z_1}\right) +2\frac{\partial
K}{\partial\xi_2}\frac{\partial^2 f} {\partial z_1\partial w_2}.
$$
Other derivatives $\partial \tilde K/\partial
\bar z_1$, $\partial \tilde K/\partial
 w_2$ and $\partial \tilde K/\partial
\bar w_2$ can be computed similarly. They justify the
assertion for $k=1$. It is easy to see that the
assertion is true for any $k$. $\blacksquare$

\v When $\|f\|_{C^2(U)}\leq C$ for some constant $C>0,$ applying Lemma \ref{lemma_6.4.2} with $k=1$, we have $\|\tilde{K}\|_{C^1(U)}\leq C'.$
Since $C_0\inv\leq (g_{i\bar j})\leq C_0$ in $U$ for some constant $C_0>0,$ by \eqref{eqna_6.4} we get
$C_1\inv\leq   (f_{i\bar j}) \leq C_1  $ for some constant $C_1>0.$
Then applying Krylov-Safonov's estimate to the equation \begin{equation}\label{complex_scalar_curvature_equations}-\sum f^{i \bar j} V_{i \bar j}=\tilde K,
\end{equation}
 we have $V\in C^{\alpha}(U).$  To obtain
a $C^{2,\alpha}$ estimate we need the following lemma (see \cite{D-Z-Z}).
\begin{lemma}[Dinew-Zhang-Zhang] \label{Din-Z-Z}
 Let $f \in \mathcal R^{\infty} (D_1(0))$  be a function that solves the equation
$$det(f_{i\bar j})= W \;\;\;\;\; in \;\;\;\;D_1(0).$$
Assume that $W\in C^{\alpha}(D_1(0))$  has the property
$\inf_{D_1(0)} W\geq \lambda > 0, $
and
 $$\|f\|_{C^2(D_1(0))}\leq C_0 $$
for some $C_0$. Then  $f\in C^{2,\alpha}(D_1(0))$ and there exists a constant $C_1$ depending only on $n,\alpha,C_0$ and $\|W\|_{C^{\alpha}(D_1(0))}$
such that $$\|f\|_{C^{2,\alpha}(D_{\half}(0))}\leq C_1.$$
\end{lemma}
Note that $W=e^{V}\in C^{\alpha}(U).$ It follows from Lemma \ref{Din-Z-Z} that $f\in C^{2,\alpha}(U). $ Then $f^{i \bar j}\in C^{\alpha}(U).$ Hence by Schauder estimates of the equation \eqref{complex_scalar_curvature_equations} we have
 $V\in C^{2,\alpha}(U).$
 Applying a  bootstrapping argument to the equation $\det(f_{k\bar l})=e^{V},$ we have for any compact set $U_1\subset\subset U$ $$\|f\|_{C^{4,\alpha}(U_1)}\leq C_1'.$$ Hence, it is crucial to get the interior estimate of
$\|f\|_{C^2}$. We need the following lemma.

\begin{lemma}\label{lemma_6.4.3}
 Let $\chi_1,\chi_2$ be the
eigenvalues of the matrix $(\sum g^{i \b j}f_{k\b j}).$ Suppose
that in $U$
\begin{equation}\label{eqn_7.3}
  C_1\inv \leq \chi_1 \leq   \chi_2 \leq
C_1,\;\;\;\; \|D^2g\|_{C^0(U)}\leq C_1
\end{equation}
  for some   constant $C_1>0$. Then there exists a constant $C_2>0$  such
  that
$$\|D^2f\|_{C^0(U)} \leq C_2.$$
\end{lemma}
{\bf Proof.}  With log-coordinate $(w_1, w_2)$ we have
$$\frac{\partial^2 f}{\partial w_i \partial w_j} =\frac{\partial^2 f}{\partial w_i \partial {\bar w_j}}=\frac{\partial^2 f}{\partial {\bar w_i} \partial {\bar w_j}},\;\;  $$
and
\begin{equation}\label{eqn_7.4}\lim_{x_1\to - \infty}\frac{\p g}{\p x_1} (x)=\lim_{x_1\to -
\infty}\frac{\p f}{\p x_1} (x)=0.\end{equation} It follows from
\eqref{eqn_7.3} and \eqref{eqn_7.4} that
\begin{equation}\label{eqn_7.5}0<\frac{\p}{\p x_1}f\leq C_1 \frac{\p}{\p x_1}g.
\end{equation}
Obviously, \eqref{eqn_7.3} gives us
\begin{equation}\label{eqna_7.5}
\left|\frac{\partial^2 f}{\partial w_2  \partial {\bar w_2}}\right| + \left|\frac{\partial^2 f}{\partial z_1 \partial {\bar z_1}}\right|
+\left|\frac{\partial^2 f}{\partial z_1 \partial {\bar w_2}}\right|+\left|\frac{\partial^2 f}{\partial w_2 \partial {\bar z_1}}\right|
\leq C_3\end{equation} for some constant $C_3>0.$
Since $\frac{\partial f}{\partial {\bar w_2}}=\frac{\partial f}{\partial w_2},$ we have
$$ \left|\frac{\partial^2 f}{\partial w_2  \partial {w_2}}\right|+ \left|\frac{\partial^2 f}{\partial z_1 \partial {w_2}}\right|\leq C_3.$$
 Note  that
$$\frac{\p^2f}{\p z_1^2} =  \frac{\p^2f}{\p w_1^2}\left(\frac{\p w_1}{\p z_1}\right)^2
+ \frac{\p^2 w_1}{\p z_1^2} \frac{\p f}{\p w_1}.$$
Using \eqref{eqna_7.5},   we have
\begin{eqnarray*}\left|\frac{\p^2f}{\p w_1 ^2}\left(\frac{\p w_1}{\p z_1}\right)^2\right| =  \left|\frac{\p^2f}{\p w_1 \p {\b
w_1}}\frac{\p w_1}{\p z_1}\frac{\p {\b w_1}}{\p {\b z_1}}\right|
 = \left|\frac{\partial ^2 f}{\partial z_1\partial \bar{z}_1}\right|\leq C_3,
\end{eqnarray*} and by \eqref{eqn_7.5}, we have
\begin{eqnarray*} \left|\frac{\p^2 w_1}{\p z_1^2} \frac{\p f}{\p w_1}\right|
&\leq& C_2\left|\frac{\p^2 w_1}{\p z_1^2} \frac{\p g}{\p w_1}\right| \\
&=&C_2\left|\frac{\p^2 w_1}{\p z_1^2} \frac{\p g}{\p w_1}+\sum
\frac{\p^2g}{\p w_1 \p w_1}\frac{\p w_1}{\p z_1}\frac{\p w_1}{\p
z_1} -\sum \frac{\p^2g}{\p w_1 \p w_1}\frac{\p w_1}{\p
z_1}\frac{\p
w_1}{\p z_1}\right|  \\
&\leq& C_2\left|\frac{\p^2 g}{\p
z_1^2}\right|+C_2\left|\frac{\p^2 g}{\p z_1 \p{\b z_1}}\right|.
\end{eqnarray*}
Then
$\|D^2 g\|_{C^2(U)}\leq C $ implies that
$$\|D^2f\|_{C^2(U)}\leq C.$$ The lemma is proved.
$\blacksquare$

 \v\n
{\em Case 3, bootstrapping on $\CHART_\vartheta$.} The proof
is same as that of Case 2.

\v We complete the proof of Theorem \ref{theorem_6.4.1}.

 \v  By subtracting   a linear function from $f$ we may
normalize $f$ such that $f(z)\geq f(z_o)=0$.  This can be easily done
if $\CHART=\CHART_\Delta$.    When $\CHART=\CHART_\ell$, this is
explained in the beginning of \S\ref{sect_7.1}. Finally, when
$\CHART=\CHART_\vartheta$, this can be done similarly.  We leave the
verification to readers.
\v
 Combining Lemma \ref{lemma_6.3.1} and Theorem \ref{theorem_6.4.1}
  we have the following
theorem.

\begin{theorem}\label{corollary_6.4.4}
 Let $z_o\in\CHART_\bullet$, where $\bullet$ can be $\Delta, \ell_i$ or $\vartheta_i$, and
 $B_a(z_o)$ be a geodesic ball in $\CHART_\bullet$.
 Suppose that there is a
 constant $C_1>0$ such that
$f(z_o)=0,\;\; \nabla f (z_o)=0,\;\;$ and $$ \mc K(f) \leq C_1,\;\;\;\;    W \leq   C_1,\;\;\; |z|\leq C_1 $$
in  $B_a(z_o)$.
 Then there is a constant $a_1>0$, depending on $a$ and $C_1$,
such that  $D_{{2a_1}}(z_o)\subset B_{\frac{a}{2}}(z_o),$ and for any $k\geq 0,$
$$
\|f\|_{C^{k+3,\alpha}(D_{a_1}(z_o ))}\leq C(a,C_1, \|\mc S(f)\|_{C^k(D_{a_1}(z_o ))}).
$$
\end{theorem}
{\bf Proof.}  Let $\lambda_1,\lambda_2$  (resp. $\mu_1,\mu_2 $)  be the eigenvalue of the matrix  $(f_{i \bar j})$ (resp. the matrix $(g_{i\bar j})$). By Lemma \ref{lemma_6.3.1},
 we conclude that  $$C_2\inv\leq \lambda_1\leq\lambda_2\leq C_2,\;\;\; in\;\;\; B_{\frac{a}{2}}(z_o)$$ for some constant $C_2>0$.
 For  the Guillemin metric
 $$C_3\inv\leq  \mu_1\leq \mu_2\leq C_3,\;\;\; in\;\;\; B_{\frac{a}{2}}(z_o),\;\;\mbox{ and }\;\;\;\;\;\; \|g\|_{C^2(D_{C_1}(0))} \leq C_3$$ for some constant $C_3>0$.
 Hence,
by Lemma \ref{lemma_6.4.3}, we have bounds on $|D^2f|$. Then there is a constant
$a_1 >0$  such that
 $
D_{{2a_1}}(z_o)\subset B_{\frac{a}{2}}(z_o).$  By integrating we have  an  upper
bound of the $C^1$-norm of $ f$  in $D_{2a_1}(z_o)).$
Combining these together, we have  a bound on the $C^2$-norm of $f$ in
$D_{2a_1}(z_o)$. Then
  the theorem follows from Theorem \ref{theorem_6.4.1}. $\blacksquare$

\begin{remark}
The results in this subsection can be easily generalized to  higher dimensional cases.
\end{remark}


\section{Estimates of $\mc K$ near divisors}\label{sect_7}
\def \fkz{\mathfrak{z}}
Let $\Delta\subset \mathbb R^2.$ In this section, we prove the following theorem.
\begin{theorem}\label{theorem_7.0.1}
Let $u\in \mc C^\infty(\Delta,v)$. Let $\fkz_o$ be a point on a
divisor $Z_\ell$ for some $\ell$. Choose a coordinate system $(\xi_1,\xi_2)$ such that $\ell=\{\xi|\xi_1=0\}.$
Let $p\in \ell$ and $D_b(p)\cap \bar \Delta$ be an Euclidean half-disk such that its intersection  with $\partial\Delta$ lies in the
interior of $\ell$. Let $B_a(\fkz_o)$ be a closed geodesic
ball satisfying  $\tau_f(B_a(\fkz_o))\subset D_b(p)$. Suppose that
\begin{equation}\label{eqn_7.1}
 \min_{D_b(p)\cap \bar \Delta}|\mc S(u)|\geq \delta>0,\;\;\;
 \|\mc S(u)\|_{C^3(D_b(p)\cap \bar \Delta)}\leq\mff N_5,  \;\;\;
 h_{22}|_{D_b(p)\cap\ell}\geq \mff N_5\inv
\end{equation}
 for some constant $N_5>0,$ where $h=u|_{\ell}$ and $\|.\|_{C^3(\Delta)}$ denotes the Euclidean  $C^{3}$-norm.
 Then there is a
 constant $\mff C_{12}>0$, depending only on $a,b,\delta$ and $\mff N_5$,
  such that
\begin{equation}\label{eqn_7.2}
 \frac{W^{\frac{1}{2}}(\fkz)}{\max\limits_
 {B_a(\fkz_o )} {W}^{\frac{1}{2}}}
\left( \mc K(\fkz)+\|\nabla \log |\mc S|\|^2(\fkz)+\Psi(\fkz)\right)a^2 \leq \mff C_{12},\;\;\;\forall \fkz\in B_{a/2}(\fkz_o )
\end{equation}
where $W=\detf$.
\end{theorem}

\subsection{Affine transformation rules on $\cplane\times\cplane^\ast$}\label{sect_7.1}
Recall that the coordinate chart for $\ell$ is $\CHART_\ell\cong
\cplane\times \cplane^\ast$. By the assumption, $B_a(\fkz_o)$
is inside this chart. The situation is  the same for
$\CHART_{\halfplane}$. Hence for the sake of simplicity of
notations, we assume that we are working on $\CHART_{\halfplane}.$
  Let $u\in \mc C^\infty(\halfplane,
v_{\halfplane})$. Then $f=L(u)$ is a function on $\t$.
Hence it defines a function on  $\cplane^\ast\times\cplane^\ast
\subset \CHART_{\halfplane}$ in terms of log-affine coordinate
$(w_1,w_2)$. Then the function $  f_\mathsf{h}(z_1,w_2):=f(\log |z_1^2|,
Re(w_2)) $ extends smoothly over $Z$, hence  is defined on
$\CHART_{\halfplane}$.  We denote $f_\mathsf{h} $ by $f$ to simplify notations.
Then $$\lim_{x_1\to-\infty}\frac{\partial f}{\partial x_1}=0,\;\;\;\frac{\partial f}{\partial x_1}>0.$$

Fix a point $\fkz_o\in Z$, we claim that by substracting
  a linear function $l$ on $\t$, for $\hat f=f-l$,
  $\hat f(\fkz_o)=\min \hat f$. In fact, let
  $(0,a)=\tau_f(\fkz_o)$. Then
$$
\hat f(x_1,x_2)=f(x_1,x_2)- ax_2.
$$
In terms of complex coordinates,
$$
\hat f(z_1,w_2)=f(z_1,w_2)-\frac{1}{2}a(w_2+\bar w_2).
$$
We find that $\tau_{\hat f}(\fkz_o)=(0,0)$. Hence it is easy to
show that $\hat f$ achieves its  minimum at $\fkz_o$.

We will explain how the affine transformation
affects the invariants of functions.

Let $u\in \mc C^\infty(\halfplane,v_{\halfplane})$. We consider the following
affine transformation on $u$:
\begin{equation}\label{eqn_7.1a}
\tilde u(\xi)=\lambda u(A\inv(\xi))+\eta\xi_1+b\xi_2+c,
\end{equation}
where $ A(\xi_1,\xi_2)=(\alpha\xi_1,\beta \xi_2+\gamma ). $ Let $\tilde
f=L(\tilde u)$.
  In this section, we always take $\lambda=\alpha$. Then  $\tilde f$
  is still a potential function of $\CHART_\ell$. $\tilde f$ can be computed directly.
The above transformation induces a transformation on $\t$:
$$
B(x_1,x_2)=( x_1+\eta,\frac{\lambda}{\beta}x_2+b).
$$
Then by a direct computation, we have
\begin{lemma}\label{lemma_7.1.1}
$\tilde f(x) =\lambda f\circ B\inv(x)+\gamma x_2-c-\gamma b$.
\end{lemma}
Now we explain the coordinate change  (in the complex sense) that covers the transformation $B$. Let
\begin{center}
\begin{tabular}{rll}
$ \mathbb R^1 \times \mathbb R^1$ & $\longrightarrow$  &  $  \mathbb R^1 \times \mathbb S^1$\\
$(x_2,\breve y_2)$ & $\xrightarrow{(I ,p_r)} $  &  $(x_2,y_2).$
\end{tabular}
\end{center}
be the universal holomorphic covering of $\cplane^\ast$. It induces a covering map
 \begin{center}
\begin{tabular}{rll}
$\mathbb C \times \mathbb C $ & $\longrightarrow$  &  $\mathbb C  \times \mathbb C^*$\\
$ (z_1,\breve w_2)$ & $\xrightarrow{\breve p_r}  $  &  $(z_1,w_2),$
\end{tabular}
\end{center}
where  $\breve w_2=x_2+\sqrt{-1}\breve y_2,w_2=x_2+\sqrt{-1}y_2,$ and  $ y_2= p_r(\breve y_2).$

 Since $f$ is $\mathbb T^2$-invariant function, i.e., $f$ is independent of $y$, so $f$ can be naturally treated as the function defined in the coordinate $(z_1,\breve w_2)$.
As in \eqref{eqn_5.2},  let $W $ be the   determinant
of the Hessian of   $f$. Then
\begin{lemma}\label{lemma_7.1.2}
For any $
u\in \mc C^\infty(\halfplane,v_{\halfplane}),
$  consider the following
affine transformation on $u$:
\begin{equation}\label{eqnc_7.1}
\tilde u(\xi)=\lambda u(A\inv(\xi))+\eta\xi_1+b\xi_2+c,\;\;\;\;\;
\end{equation}
where $ A(\xi_1,\xi_2)=(\alpha\xi_1,\beta \xi_2+\gamma )$ and $\lambda=\alpha.$
Then it  induces an affine transformation  in  complex coordinate $(z_1,\breve w_2):$
 \begin{equation}\label{eqnc_7.2}
 B_\cplane(z_1,\breve w_2)=(e^{\frac{\eta}{2}}z_1,\frac{\alpha}{\beta}\breve w_2 +b).\;\;\;\;
\end{equation}
Moreover,
\begin{itemize}
\item $\tilde f(z)=\alpha f (B_\cplane\inv z)+\gamma x_2-\gamma b-c;$
\item $\tilde W(z)=\beta^2e^{-\eta}W( B\inv_{\cplane}z)$;
\item  $\tilde \Psi(z)=\alpha\inv \Psi(B\inv_{\cplane}z)$;
\item  $\mc {\tilde {K}}(z)=\alpha\inv \mc K(B\inv_{\cplane}z)$;
\item  $ \|\nabla\log |\mc S(\tilde f)|\|^2_{\tilde f}(z)=\alpha\inv \|\nabla\log |\mc S( f)|\|^2_{f}(B\inv_{\cplane}z)$,
\end{itemize}
 where $z=(z_1,\breve w_2)$.
\end{lemma}
Note that all these functions are $\mathbb T^2$-invariant, hence the formulae can be  pushed forward from $\cplane\times \cplane$
 down to $\cplane\times\cplane^\ast $.
The   lemma follows from a direct calculation. We omit it.

\subsection{Uniform control of sections}\label{sect_7.3}

In this subsection, we consider functions
$
u\in \mc C^\infty(\halfplane,v_{\halfplane};K_o)
$ (i.e., $u=\xi_1\log\xi_1 + \xi_2^2+ \psi$ is strictly convex and $|\mc S(u)|\leq K_o$  for some $\psi\in C^\infty(\halfplane)$)
with the property
\begin{equation}\label{eqn_7.6}
\Theta_u(p) d_u^2(p, \t^\ast_2)\leq \mff C_5.
\end{equation}

 Let $\fkz^\circ \in
\CHART_{\halfplane}$ be any point such that $d(\fkz^\circ,Z)=1$ and
$\fkz^\ast\in B_1(\fkz^\circ)\cap Z$.
Without loss of generality, we assume that $\fkz^\circ$ is a representative point
of its orbit and assume that it is on $\t$ (cf.\;Remark \ref{remark_1.1.3}).
Let $p^\circ,p^\ast$ be
their images of moment map $\tau_f$.
\begin{remark}
In this section, when we consider a point $z\in \CHART_{\halfplane}$,
without loss of generality, we assume $z$ is the representative
of its $\bb T^2$-orbit (cf Remark \ref{remark_1.1.3}).
Hence when $p=\tau_f(z)$, we assume that $z=\nabla^u (p)\in
\t$. If $z$ is on $Z$, we may assume that it is the point in
$\t_2$. For such $z$ we write $\dot\tau_f\inv(p)$.
\end{remark}

 By adding a linear function we normalize $u$ such that $p^\circ$ is the minimal point of $u$; i.e.,
 \begin{equation}\label{eqnc_7.4}
  u(p^\circ)=\inf u.
 \end{equation}
Let $\check p$ be the minimal point of $u$ on $\t_2^\ast$ which is the boundary of $\halfplane$.
By adding some constant to $u$, we may require that
\begin{equation}\label{eqnc_7.5}
u(\check p)=0,
\end{equation}
and by a coordinate translation we may assume that \begin{equation}\label{eqnc_7.6} \xi(\check p)=0.\end{equation}
(see Figure 1). Denote $z^\ast=\nabla ^u (p^\ast).$ Set $$S_0:=\left\{z\in \t_2 \;|\; \left|\int^{x_2(z^\ast)}_{x_2(z)} \sqrt{f_{22}}dx_2\right|\leq 1 \right\}.$$ By a coordinate transformation \begin{equation}\label{eqnc_7.3}A(\xi_1,\xi_2)=(\xi_1,\beta\xi_2)\end{equation} we can normalize $u$ such that
 \begin{equation}\label{eqn_7.7}
\left|S_0\right|= 10.
  \end{equation}
In fact,  \eqref{eqnc_7.3} induces a coordinate transformation in $(x_1,x_2)$ as following $$A(x_1,x_2)=(x_1,\beta\inv x_2),$$
then by choosing the proper $\beta$ we have \eqref{eqn_7.7}.
The above transformations are affine transformations in \eqref{eqnc_7.1} with $\lambda=\alpha=1.$ By Lemma \ref{lemma_7.1.2} it is easy to see that  $\mathcal K,\Psi,\|\nabla\log |\mathcal S(f)|\|_f^2$ and  $W(z)/W(z')$ for any $z,z'\in \CHART_{\halfplane}$  are invariant under these transformations.

We say
$(u,p^\circ, \check p)$ is a {\em minimal-normalized-triple}
if $u$ satisfies \eqref{eqnc_7.4}, \eqref{eqnc_7.5}, \eqref{eqnc_7.6},  \eqref{eqn_7.7} and
$$d(p^\circ, \t_2^\ast)=1.$$

Note that $\check p$
is determined by $u$ and $p^\circ$ already.

 \begin{figure}
 \center
\includegraphics[height=2in]{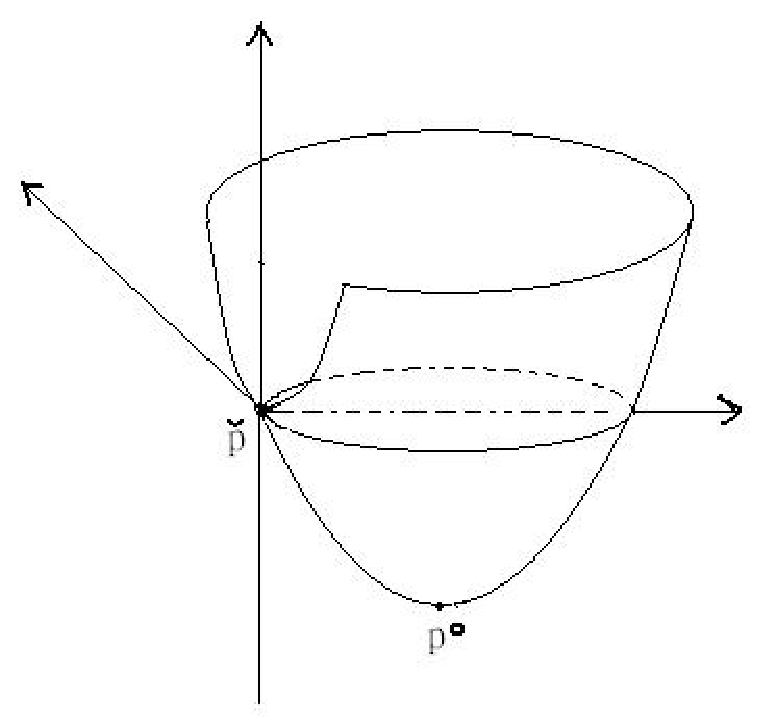}
 \caption{\label{fig1}}
 \end{figure}



\begin{lemma}\label{lemma_7.3.1}
Let $(u, p^\circ, \check p)$ be a minimal-normalized triple.
 Then there exists a constant  $ C_1>1 $
such that
$$ u(\check p)-
u(p^\circ)\geq C\inv_1.$$
\end{lemma}
{\bf Proof.} By \eqref{eqn_7.6} we have $$\Theta\leq 4\mff C_5,\;\;\;\;\; in\;
B_\half(p^\circ) .$$  Then by  an
argument similar to that in the proof of the claim in Theorem \ref{theorem_2.4.3}, we conclude that
$u(p)-u(p^\circ)\geq \delta$ for any $p$ not in
$B_\half(p^\circ)$. In particular, $\check p$ is such a point.
Hence, we prove the lemma. $\blacksquare$

\v
\begin{lemma}\label{lemma_7.3.2}
Let $(u\indexm, p\indexm^\circ, \check p\indexm)$
 be a sequence of minimal-normalized triples  with
\begin{eqnarray*}\lim_{k\to \infty} \max|\mc S(u\indexm)|
=0,\;\;\;\;\;\;\;\;
\end{eqnarray*}
 Then there exists a constant  $ C_1>1 $
such that
$$C_1\inv\leq  u\indexm(\check p\indexm)-
u(p^\circ\indexm)\leq C_1$$
when $k$ large enough.
\end{lemma}
{\bf Proof.} The lower bound is proved in Lemma \ref{lemma_7.3.1}.

By adding a constant to $u\indexm$, we assume that
$u\indexm(p^\circ\indexm)=0$.
 Now suppose that  $ u\indexm(\check p \indexm) $ has no upper
bound. Then we can choose a sequence  of constants $N\indexm\to
\infty$ such that
$$
 0<N\indexm<u\indexm(\check p \indexm),\;\;
\lim_{k\to \infty} N\indexm \max|\mc S(u\indexm)|=0.  $$ For each
$u\indexm$ we take an affine transformation $\hat A\indexm:=(A\indexm,
(N\indexm)\inv)$ (cf.\;Definition \ref{defn_2.2.1}) to get a new function $\tilde u\indexm$, i.e,
$$
\tilde u\indexm= (N\indexm)\inv u\indexm\circ (A\indexm)\inv.
$$
Then original section $S_{u\indexm}(p^\circ,N\indexm)$ is
transformed to be $S_{\tilde u\indexm}(A\indexm p\indexm^\circ, 1)$. We
choose $A\indexm$ such that the latter one is normalized (cf.\;\S\ref{sect_3.1}).
Then, by Lemma \ref{lemma_2.2.1},
$$ \lim_{k\to \infty}\mc S(\tilde u\indexm)=
\lim_{k\to \infty}N\indexm\mc S(u\indexm)= 0,$$
$$\lim_{k\to\infty} d_{\tilde u\indexm}(A\indexm p^\circ\indexm,
A\indexm\t^\ast_2) =\lim_{k\to\infty}\frac{1}{\sqrt{N\indexm}}
 d_{u\indexm}(p^\circ,\t_2^\ast)=0. $$
On the other hand,  by Theorem \ref{theorem_3.2.7} we conclude that
 $\tilde u\indexm$ locally  $C^{3,\alpha}$-converges to a
  strictly convex function $\tilde u_\infty$  in a neighborhood of $p^\circ_\infty$,
which is  the minimal point of $u_\infty$ and the limit of $A\indexm p^\circ\indexm$. In
  particular, there is a constant $C_3>0$ such that
  $$
d_{\tilde u\indexm}(A\indexm p^\circ\indexm, A\indexm\t_2^\ast)\geq
C_3.  $$
 We get a contradiction.
 $\blacksquare$

\begin{defn}\label{definition_7.3.3}
Let $(u,p^\circ,\check p)$ be a minimal-normalized-triple. Let
$$N=\max(200, 100bC_1,100\exp(100\sqrt {\mff C_{10}})),$$
where $b$ is the
constant in Proposition \ref{proposition_3.2.6},  $C_1$ is the constant in Lemma \ref{lemma_7.3.2}  and $\mff C_{10}$ is the constant in Lemma \ref{lemma_6.2.1} with $N_0=4,n=2$.
If the following inequalities hold
\begin{eqnarray}
&& \mc K(z)\leq 4, \;\;\;\forall z\in \tau_f\inv(B_N(p^\circ));\label{eqn_7.8}\\
&& \frac{1}{4}\leq \frac{W(z)}{W(z')} \leq 4,\;\;\;\;
\forall z,z'\in \tau_f\inv(B_N(p^\circ)),\label{eqn_7.9}
\end{eqnarray}
we say that $(u,p^\circ, \check p)$ is a bounded-normalized triple.
\end{defn}

\v Let $u\indexm$ be a sequence
   given in Lemma \ref{lemma_7.3.2} and  $(u\indexm,p^\circ\indexm, \check p\indexm)$ be  a sequence of bounded-normalized triples. In this subsection, we mainly   prove that   the sections $
S_{u\indexm}(p^\circ\indexm, \sigma\indexm) $ are $L$-normalized  for some constant $L$ that is independent of
$k$,
where $\sigma\indexm=\delta |\min u\indexm|$ and
$\delta$ is the constant in Proposition \ref{proposition_3.2.6}.
 By Lemma \ref{lemma_7.3.2}, we already know that the
$\sigma\indexm$ are uniformly bounded.
By Lemma \ref{lemma_3.1.4} we  only need to prove the following theorem.

\begin{theorem}\label{theorem_7.3.4}
Let $(u\indexm,p^\circ\indexm,\check p\indexm)$ be a sequence of
bounded-normalized triples. Suppose that $\lim_{k\to
\infty} \mc |S(u\indexm)|= 0$. Let
$$\Omega\indexm=S_{u\indexm}(p^\circ\indexm,\sigma\indexm).$$ Then
there exist constants $c_{in}$ and $c_{out}$, independent of $k$,
such that
$$
c_{in}\leq wd_i(\Omega\indexm), \;\;\; i=1,2.
$$ and
$$|\xi_i(\Omega\indexm)|\leq c_{out},i=1,2.$$
\end{theorem}
The proof consists of the following lemmas: Lemma \ref{lemma_7.3.5},
\ref{lemma_7.3.6} and \ref{lemma_7.3.9}.

\begin{lemma}\label{lemma_7.3.5}
$wd_1(\Omega\indexm)$ has a uniform lower bound.
\end{lemma}
{\bf Proof.}  Let $$ E\indexm=\{p\in S_{u\indexm}(p^\circ\indexm, \sigma\indexm)\;|\xi_2(p)=\xi_2(p^\circ\indexm)\}$$
and $p^\pm\indexm$ be
its right and left ends. Then we claim that \v\n {\em Claim. There
is a constant $C_2>0$ such that $ |p^\circ\indexm-p^\pm\indexm|\geq
C_2\inv.$}  \v\n {\em Proof of claim.} We prove the claim for
$p^-$. The proof for $p^+$ is identical.

If the claim is not true, then $
|\xi_1(p^\circ\indexm)-\xi_1(p^-\indexm)|\to 0. $
 We omit the
index $k$.  Take an affine  transformation on $u$
$$ \tilde u(\xi)=
u \left( A \inv\xi\right),$$ where  $A $ is the normalizing transformation of
  $S_{u}(p^\circ ,\delta\inv\sigma)$.
  By Theorem \ref{theorem_3.2.7} we conclude that $\tilde u_k$  $C^{3,\alpha}$-converges to a
smooth and strictly convex function $\tilde u_\infty$ in
$S_{\tilde u_\infty}(\tilde p_\infty^\circ,\sigma)$. Note that
 the geodesic distance and the ratio ${\det(u
_{ij})(\xi)}/{\det(u _{ij})(p^\circ)}$ are base-affine
 transformation invariants, we have that for any $p\in S_{u}(p^\circ,\sigma)$
 \begin{eqnarray}&&
C_3\inv\leq d(p^\circ, \partial S_{u}(p^\circ,\sigma) )\leq C_3,
\label{eqn_7.10}\\ && C_3\inv \leq {\det(u _{ij})(p)}/{\det(u
_{ij})( p^\circ )}\leq C_3, \label{eqn_7.11} \end{eqnarray}
 for some constant $C_3>1 $ independent of $k.$

On the other hand, by the convexity of $u$, we have
\begin{equation}\label{eqn_7.12}
| \partial_1u (p^-)|\geq \frac{\sigma}{|p^\circ -p^-|} \geq
\frac{\delta}{C_1|p^\circ-p^-|}\to \infty.
\end{equation}
By the coordinate relation $z_1=e^{\frac{w_1}{2}},z_2=w_2$, we  calculate  $W(z)$ in term of $\xi$-coordinates,
\begin{eqnarray*}
W(z)=\frac{1}{4}\det\left(\frac{\p^2 f}{\p x_i \p x_j}\right)e^{-x_1}( z)=\frac{1}{4}\left[\det\left(\frac{\p^2 u}{\p \xi_i \p \xi_j}\right)\right]\inv e^{-u_1}(\tau_f(z)).\end{eqnarray*}
  Let $z^-\in {\tau_f}(p^-)\cap B_N(\fkz^\circ)$.
As $ \partial_1u (p^-)<0$ and $ \partial_1u (p^\circ)=0$,  from \eqref{eqn_7.11} and \eqref{eqn_7.12} we
conclude that
$${W(z^-)}/{W(z^\circ)}=\exp({-\partial_1u(p^-)})\frac{\det(u_{ij})(p^\circ)}{\det(u_{ij})(p^-)}
\to +\infty,
$$
 which contradicts the  assumption \eqref{eqn_7.9}. This completes the proof of the claim. \v \n Since the width is at least $|p^+-p^-|$, we get its lower bound. $\blacksquare$
\v Recall  $\fkz^\ast\in Z$ is the point such that $d(\fkz^\circ,\fkz^\ast)=1$. Let ${p}^\ast=\tau_f(\fkz^\ast)$.
\begin{lemma}\label{lemma_7.3.6}
$wd_2(\Omega\indexm)$ has a uniform lower bound.
\end{lemma}
{\bf Proof.} We divide the proof  into three steps.
\v\n{\em Step
1. There is  a constant $C_4=\exp(100\sqrt {\mff C_{10}})$ such that for any $p\in B_{20}(p^\circ)\cap \t_2^\ast $
\begin{equation}\label{eqn_7.a.8a}u(p)-u(p^\circ)\leq C_4,\end{equation}
where   $\mff C_{10}$  is the constant in Lemma \ref{lemma_6.2.1} with $N_0=4,a=100,n=2.$
In particular,
\begin{equation}\label{eqn_7.13a}
u( p^\ast )-u( p^\circ )\leq C_4.\end{equation}}
  \n Proof of Step 1. Let $p^\bullet$ be the point such that
  $$u(p^\bullet)=\max_{B_{20}(p^\circ)\cap \t_2^\ast}u.$$ We only need to prove $u(p^\bullet)$ satisfies \eqref{eqn_7.a.8a}.   By  a coordinate  translation $\tilde \xi_1 =\xi_1,\tilde \xi_2 =\xi_2-\alpha$ we can assume that     $\tilde  \xi (p^\bullet)=(0,0).$ Obviously
$u(p^\bullet)-u(p^\circ)$ is invariant under the coordinate transformation.

\v Let
$
 \fkz^\circ\in {\tau}_{f}\inv (  p^\circ)$ and $
 \fkz^\bullet\in  {\tau}_{f}\inv ( p^\bullet)$ with $d(\fkz^\circ,\fkz^\bullet)\leq 20.$ Let $\tilde  f$ be the function associated with $u$ with respect to $\tilde \xi$-coordinates, i.e., $\tilde  f=\sum \frac{\p u}{\p \tilde \xi_{i}}\tilde \xi_{i}-u$.  Obviously $\tilde f (\fkz^\bullet)=\inf\limits_{\mathbb C\times \mathbb C^*}\tilde  f .$ Note that
$u=\tilde \xi_1\log\tilde \xi_1+\psi $ with $\psi \in C^\infty(\halfplane).$
Following  from  $$\lim\limits_{\tilde \xi\to p^\bullet}\tilde \xi_1(\log\tilde \xi_1+\psi_1)=0,\;\;\;\;\;\tilde \xi( p^\bullet)=0,\;\;\;\;\frac{\p  u}{\p \tilde \xi_1}( p^\circ)= \frac{\p  u}{\p\tilde  \xi_2}( p^\circ)=0$$ we have
\begin{equation}\label{eqn_7.a.8}
\sum \frac{\p u}{\p\tilde  \xi_i}\tilde \xi_i(p^\circ)=0,\;\;\;\lim_{\tilde \xi\to  p^\bullet}\sum \frac{\p  u}{\p \tilde \xi_i} \tilde \xi_i=0.\end{equation}
Then
$$\tilde  f ( \fkz^\bullet)+ u (  p^\bullet) =0,
\;\;\tilde   f (  \fkz^\circ)+ u (  p^\circ) =0.$$ It follows that
$|  u(  p^\circ ) - u(  p^\bullet)|= | \tilde  f (
\fkz^\circ)- \tilde  f ( \fkz^\bullet)|.$  Applying Lemma \ref{lemma_6.2.1} to $\tilde   f$ in $B_{100}( \fkz^\bullet)$ with $N_0=4,n=2,$  we have $$|\tilde   f (
\fkz^\circ)-\tilde   f ( \fkz^\bullet)|\leq C_4$$ for $C_4=\exp(100\sqrt {\mff C_{10}}).$
 Then the claim follows.

\v\n{\em Step 2.} Let $q_1=(-\infty,m_1)$, $q_2=(-\infty,m_2)$ be two ends of $S_0$, with $m_1<m_2$.
\v\n Denote $p_i=(0,\frac{\partial f}{\partial x_2}(m_i)).$ We have
\begin{eqnarray*}1&=&\left| \int_{m_i}^{x_2(z^\ast)}\sqrt{f_{22}}dx_2 \right|^2\;=\;\left|\int_{\xi_2(p_i)}^{\xi_2(p^*)}\sqrt{u_{22}}d\xi_2\right|^2\\
&\leq&\left|\xi_2(p_i)-\xi_2(p^*)\right| \left|\int_{\nabla^f(S_0)} u_{22} d\xi_2\right|\\
&\leq& 10|\xi_2(p_i)-\xi_2(p^* )|.
\end{eqnarray*}
It follows that
\begin{equation}\label{eqn_7.13}
|\xi_2(p_i)-\xi_2(p^*)| \geq \frac{1}{10}.\end{equation}

\v\n {\em Step 3.} For  $C_5=C_4+1,$ by the definition of $S_0$ and the result of Step 1,  we have $$  \nabla^f( S_0)\subset B_{20}(p^\circ)\cap \t_2^\ast\subset
S_u(p^\circ ,C_{5})\cap \t_2^\ast.$$
   Let  $p_3,p_4$ be the boundary of $S_u(p^\circ ,C_{5})\cap \t_2^\ast.$ Then by the result of Step 2 we have
\begin{equation}\label{eqn_7.15}
|p_3-p_4|\geq |p_1-p_2|\geq\frac{1}{5}. \end{equation}
 Take the triangle $\triangle p^\circ p_3p_4\subset
\t^*$. On the other hand, let
$$
P^\circ=(p^\circ, u(p^\circ)), P_3=(p_3,u(p_3)), P_4=(p_4,u(p_4)).
$$
They form a triangle $\triangle P^\circ P_3P_4$ in 3-dimensional
space $\t^\ast\times \real$. It is above the graph of $u$ over
$\triangle p^\circ p_3p_4$.

Now we cut the triangle $\triangle P^\circ P_3P_4$ by the
horizontal plane $\xi_3= u(p^\circ)+\sigma$, i.e:
$$
\triangle^{\sigma}=\{(\xi_1,\xi_2,\xi_3)\in \triangle P^\circ
P_3P_4| \xi_3\leq u(p^\circ)+\sigma\}.
$$
Let $\triangle_\sigma$ be the projection of $\triangle^\sigma$
onto $\t^\ast$. It is easy to see that $\triangle_\sigma\subset
S_u(p^\circ,\sigma)$.  By the similarity between triangles $\triangle^\sigma$ and
$\triangle P^\circ P_3P_4$, we conclude that the upper edge has
length at least $\sigma(10C_5)\inv$.
It follows that
$$wd_2(S_u(p^\circ,\sigma))\geq wd_2(\triangle_\sigma)\geq \sigma (10C_5)\inv. $$
This completes the proof of the Lemma.
 $\blacksquare$
\v From the first two steps in the proof, we have the following corollary.
\begin{corollary}\label{corollary_7.3.7}
There exists a  constant  $C_{6}>0 $, depending only on $\mff C_{10}$, such that
\begin{equation*}
|u_2(p^\ast)-u_2(p^\circ)|=|u_2(p^\ast)|\leq C_6 . \end{equation*}
\end{corollary}
{\bf Proof.} Without loss of generality we assume that
$\xi_2(p_1)<\xi_2(p^\ast)<\xi_2(p_2),$ where $p_1$ and $p_2$ are introduced in the Step 2 of the proof of Lemma \ref{lemma_7.3.6}. By the convexity of $u$ we conclude that
$$\frac{u(p_1)-u(p^\ast)}{\xi_2(p_1)-\xi_2(p^\ast)}\leq u_2(p^\ast)\leq\frac{u(p_2)-u(p^\ast)}{\xi_2(p_2)-\xi_2(p^\ast)}.$$
Then Corollary follows from \eqref{eqn_7.a.8a} and \eqref{eqn_7.13}.
$\blacksquare$

We resume  the proof of Theorem \ref{theorem_7.3.4}.
\begin{lemma}\label{lemma_7.3.8}
The  Euclidean volume
 of $\Omega\indexm=
 S_{u\indexm}(p^\circ\indexm,\sigma\indexm)$ has a uniform upper bound.
\end{lemma}
{\bf Proof.} First we claim that there exists a constant $C_9<N$ such that
$$
\tau_f\inv(\Omega\indexm)\subset B_{C_9}(\mathfrak z^\circ).
$$
Note that $B_{C_9}(\mathfrak z^\circ)\subset \mathbb C\times\mathbb C^*$ with the coordinate system  $(z_1,w_2).$
We omit  $k$ again. Note that for points in $\nabla^u(\Omega),$ by Proposition \ref{proposition_3.2.6}, this is
true. It remains to consider the distance for points on the fiber
over points in $\nabla^u(\Omega)$.
 Our  discussion is taking place on an
affine-log coordinate chart.

\v Suppose that
$$
\mathfrak z^\circ=(x_1^\circ,
x_2^\circ,0,0).
$$
On the one hand, by the definition of $S_0,$ we have $d_{u}(p,p^\circ)\leq 2$ for any $p\in \nabla^{f}(S_{0}).$ Then by \eqref{eqn_7.7} we can find  a point
 $ z^\alpha=(-\infty,x_2^\alpha,0,0)\in  B_2(\fkz^\circ)$ with $(-\infty,x_2^\alpha)\in S_0$ such that
\begin{equation}\label{eqn_7.17}\;\; f_{22}(z^\alpha)\leq \frac{1}{10}.
\end{equation}
This says that the metric along the torus of $y_2$-direction over $z^\alpha$ is bounded.
In fact,  if this is not true, we have $f_{22}|_{S_0}>\frac{1}{10}.$ By a direct calculation we have
$$\int_{S_0}\sqrt{f_{22}}dx_2> \int_{S_0}\frac{1}{4}dx_2=\frac{5}{2}$$ where we used \eqref{eqn_7.7}.
This contradicts the definition of $S_{0}$.
On the other hand, let $p^{\alpha}=\tau_f(z^\alpha),$ then $p^{\alpha}\in \t^\ast_2$ and
 $$\lim\limits_{\xi \to p^{\alpha}}u^{11}=\lim\limits_{\xi \to p^{\alpha}}\frac{u_{22}}{u_{11}u_{22}-u_{12}^2}=\lim\limits_{\xi
 \to p^{\alpha}}\frac{\xi_1u_{22}}{(1+\xi_1\psi_{11})u_{22}-\xi_1u_{12}^2} =0.$$ Hence
\begin{equation}\label{eqn_7.16}f_{11}(z^\alpha)=0.\end{equation}
This says that the metric along the torus of $y_1$-direction over $z^\alpha$ is bounded.


Now consider any point
$$\tilde z=(x_1,x_2,y_1,y_2)$$
over $z=(x_1,x_2,0,0)\in \nabla^u(\Omega)$.
We construct a bounded path from $\tilde z$ to $\fkz^\circ$ as  follows: \\
{\em (i) keeping $y$ coordinates  fixed,
move "parallelly" from $\tilde z$ to the point $ (-\infty,x^\alpha_2,y_1,y_2)$  over $z^\alpha$,   denote
$\hat z^\alpha=(-\infty,x^\alpha_2,y_1,y_2)$; \\
(ii) move along the fiber over $z^\alpha$ along $y_1$-direction from $\hat z^\alpha$ to the point  $(-\infty,x^\alpha_2,0,y_2),$  denote
$\tilde z^\alpha=(-\infty,x^\alpha_2,0,y_2)$;\\
(iii) move along the fiber over $z^\alpha$ along $y_2$-direction from $\tilde z^\alpha$ to $z^\alpha$; \\
(iv) keeping $y$ coordinates fixed, move from $z^\alpha$ to $\fkz^\circ$.}\\
By a direct calculation we have
\begin{eqnarray*}d(\tilde z,\fkz^\circ)&\leq& d(\tilde z, \hat z^\alpha)+
d(\hat z^\alpha, \tilde z^\alpha)+d( \tilde z^\alpha,  z^\alpha ) +d(   z^\alpha, \fkz^o )\\
&\leq& 2+C_1b+4\pi+ 2:=C_9,
\end{eqnarray*}
where we used \eqref{eqn_7.17}  and \eqref{eqn_7.16} in the second inequality.  Here $b$ is the constant in Proposition \ref{proposition_3.2.6} and $C_1$ is the constant in Lemma \ref{lemma_7.3.2}. The claim is proved.

 Since the Ricci curvature tensor is bounded, by the Bishop Volume Comparison Theorem we have
$$\mathrm{Vol}(\tau_f\inv(\Omega))\leq \mathrm{Vol}(B_{C_9}(\fkz^\circ))\leq C_{10}$$
for some constant $C_{10}>0$. By a direct calculation we have
\begin{eqnarray*} 16\pi^2 \mathrm{Vol}_E(\Omega)&=&16\pi^2\int_{\Omega}d\xi_1d\xi_2\nonumber\\
&=&\int_{\mathbb T^2} \int_{\nabla^u(\Omega)}\det\left(\frac{\p^2 f}{\p x_i \p x_j}\right)dx_1dx_2dy_1dy_2\nonumber\\
 &=&\mathrm{Vol}(\tau_f\inv(\Omega))\leq  C_{10},\end{eqnarray*}
 where   $\mathrm{Vol}_E(\Omega)$ denotes the Euclidean  volume of the set $\Omega.$
$\blacksquare$ \v We have the following consequence,

\begin{lemma}\label{lemma_7.3.9}
There is a Euclidean ball $D_R(0)$ such that for any k\;
$$S_{u\indexm}(p\indexm^\circ, C_5)\subset D_{R}(0).$$
In particular, $$\Omega\indexm\subset D_R(0)$$ and
\begin{equation}\label{eqn_7.a9}|\xi(p\indexm^\ast)|\leq R,\;\; |\xi(p\indexm^\circ)|\leq R.\end{equation}
\end{lemma}
{\bf Proof.} Since $wd_i(\Omega\indexm)$ is bounded below, let $c_{in}$
be the bound.
Let
$$
\tilde\Omega\indexm=S_{u\indexm}(p^\circ\indexm, |u\indexm(p^\circ\indexm)|+C_5).
$$
Then by the convexity of $u\indexm$
 and the similarity we have
\begin{eqnarray*}
\mathrm{Vol}_E(\tilde\Omega\indexm)&\leq& \left(\frac{|u\indexm(p^\circ\indexm)|+C_5}
{\sigma\indexm}
\right)^2\mathrm{Vol}_E(\Omega\indexm)\\&\leq&  \frac{(1+C_1 C_5)^2\mathrm{Vol_E}(\Omega\indexm)}{\delta^2}\leq
\frac{(1+C_1 C_5)^2C_{10}}{\pi^2\delta^{2}}.
\end{eqnarray*}
Note that $\check{p}=(0,0)\in \tilde\Omega\indexm$ and
$$
wd_i(\tilde\Omega\indexm)\geq wd_i(\Omega\indexm)\geq c_{in}.
$$
Then there is a constant $R>0$ independent of $k$ such that for any $p \in \tilde\Omega\indexm$
$|\xi_{i}(p)|\leq \frac{R}{2},\mbox{ for } i=1,2$. $\blacksquare$

\begin{lemma}\label{lemma_7.a.9}
There exists a constant $C_{11}>0$ independent of k, such that
for any  $p\in B_8(\fkz^o),$
$$|z_1|(p)+|w_2|(p)\leq C_{11}.$$
\end{lemma}
{\bf Proof.} Note that
$$|z_1|^2=e^{\frac{w_1+\bar w_1}{2}}=e^{x_1},\; |w_2|^2\leq |x_2|^2+16\pi^2.$$
We only need  to prove that  $x_1$ and $|x_2|$ are  uniformly  bounded above.
Set  $$\tilde f=f-\frac{\partial f}{\partial x_2}(q_1)(x_2-x_2(q_1))-f(q_1),\;\;\;\;\hat f=  f- \frac{\partial f}{\partial x_2}(z^\ast)(x_2-x_2(z^\ast))-f(z^\ast).$$
 Then
$$\tilde f (q_1)=\inf \tilde f=0,\;\;\;\;\;\hat f (z^\ast)=\inf \hat f=0.$$
Using  Lemma \ref{lemma_6.2.1}   we have
 \begin{equation}\label{eqn_7.a.3}
|\tilde f( q)-\tilde f(q_1) |\leq N_1,\;\;\;\;
|\hat f( q)-\hat f(z^\ast) |\leq N_1,
\end{equation}
for any $q\in B_8(\fkz^o),$ where $N_1= \exp\{20\sqrt{\mff C_{10}}\}.$
Note that  $$\tilde f=\hat{f} +(\xi_2(p^\ast)-\xi_2(p_1))x_2+d_0,$$
where $d_0= \xi_2(p_1) x_2(q_1)-f(q_1)-\xi_2(p^\ast) x_2(z^\ast)+f(z^\ast).$
 Then by  \eqref{eqn_7.a.3} we have
\begin{eqnarray}\label{eqn_7.a.4}
2N_1&\geq&  |\tilde f( q)-\tilde f(q_1)|+|\tilde f(z^\ast)-\tilde f(q_1)|
 \nonumber\\ &\geq&  |\tilde f( q)-\tilde f(z^\ast)| \nonumber\\
&\geq & |\xi_2(p^*)-\xi_2(p_1)||x_2(z^\ast)-x_2(q)| -|\hat f( q)-\hat f(z^\ast) |\nonumber\\
&\geq&   \frac{1}{10}|x_2(z^\ast)-x_2(q)|-N_1
\end{eqnarray}
where we use \eqref{eqn_7.13} and \eqref{eqn_7.a.3} in the last inequality. It follows from \eqref{eqn_7.a.4} and  Corollary \ref{corollary_7.3.7} that
 \begin{equation}\label{eqn_7.a.7}
 |x_2(q)|\leq |x_2(z^\ast)| +30N_1 \leq C_6 +30N_1=N_2\end{equation} for any $q\in B_8(\fkz^o). $

 Let  $q\in B_8(p^o)$ be a point with $u_1(q)=\max\limits_{B_8(p^\circ)}u_1.$ Then $$x_1(\nabla^u(q))=\max\limits_{p\in B_8( \nabla^u(p^\circ))}x_1(p).$$  Set $\tilde u=u-u_1(q)\xi_1-u_2(q)\xi_2.$
 Then $\tilde u(q)=\inf \tilde u.$  As in the Step 1 of Lemma \ref{lemma_7.3.6}, we have
  $\tilde u(p^\ast)-\tilde u(q) \leq C_4$.
Since $\tilde u(q)=\inf \tilde u,$ we conclude that
\begin{eqnarray}\label{eqn_7.a.5} C_4 &\geq& \tilde u(p^\ast)-\tilde u(p^\circ)\nonumber\\
&=& (u(p^\ast)- u(p^\circ))+u_1(q)\xi_1(p^\circ)-u_2(q)(\xi_2(p^\ast)-\xi_2(p^\circ))\nonumber\\
&\geq & u_1(q)\xi_1(p^\circ)-N_2|\xi_2(p^\ast)-\xi_2(p^\circ)|,
\end{eqnarray}
where  we use \eqref{eqn_7.a.7} and  $u(p^\circ)=\inf u$ in the last inequality.
By  Lemma \ref{lemma_7.3.9}  we have
\begin{equation}\label{eqn_7.a.6}
|\xi_2(p^\ast)|\leq R ,\;\;\;|\xi_2(p^\circ)|\leq R,\;\;\;|\xi_2(p^\ast)-\xi_2(p^\circ)|\leq 2R.
\end{equation}
The claim in the proof of Lemma \ref{lemma_7.3.5} gives us  $\xi_{1}(p^{\circ})\geq C_{2}^{-1}>0.$ Substituting \eqref{eqn_7.a.6} into \eqref{eqn_7.a.5}   we have
$$x_1(\nabla^u( q))\leq  C_{12}$$ for some constant $C_{12}>0$ independent of $k.$
 Hence on $B_8(\fkz^o),$  $|z| $ is uniformly bounded from above.  $\blacksquare$

\subsection{A convergence theorem}\label{sect_7.4}
We use Theorem \ref{theorem_7.3.4} to prove a convergence
theorem.
\begin{theorem}\label{theorem_7.4.1} Let $(u\indexm, p^\circ\indexm,
\check p\indexm)$ be a sequence of
bounded-normalized triples. Suppose that $\lim_{k\to
0} \mc |S(u\indexm)|\to 0$.
 Then
 \begin{enumerate}
 \item[(1).]  $p^\circ\indexm$ converges to a point $p^\circ_\infty$ and $u_k$
$C^{3,\alpha}$-converges to a strictly convex function $u_\infty$
in $D_{s}(p^\circ_{\infty})
\subset\halfplane$, where $s$ is a constant independent of $k.$
\item[(2).] 
 $\fkz^\ast\indexm$ converges to a point $\fkz^\ast_\infty$ and $f \indexm$
   $C^{3,\alpha}$-converges to a
 function $f_\infty $ in $D_{a_1}(\fkz_\infty^\ast)$. 
  Here  $a_1$ is the constant in Theorem \ref{corollary_6.4.4}.
 \end{enumerate}
\end{theorem}
\n{\bf Proof.}
Applying Theorem \ref{theorem_7.3.4} and Lemma \ref{lemma_3.1.4},
we conclude that
$$\Omega_k:=S_{u\indexm}(p^\circ\indexm,\sigma\indexm)$$ are $L$-normalized,
where $L$ is independent of $k$.
 Then by Theorem \ref{theorem_3.2.7} and Remark \ref{remarka_3.2.7} we conclude (1).

 By Lemma \ref{lemma_7.a.9} and $B_7(\fkz^\ast)\subset B_8(\fkz^o)$  we have,   in $B_7(\fkz^\ast)$,
\begin{equation}\label{eqn_7.a.9a}
|z|\leq C_{11},\;\;\;\mc K(z)\leq 4, \;\;\; C\inv\leq
 {W(z)} \leq C
 \end{equation} for some constant $C>0$ independent of k.  Here we use the fact that $W(\fkz^\circ)$ is bounded due to
(1).
By adding a linear function  we can assume that
   $f(\fkz^\ast)=\inf_{\mathbb C\times \mathbb C^*}f=0.$


 Note that the geodesic distance and \eqref{eqn_7.a.9a} are invariant under the transformations of adding  linear functions.
 We know  that $f$ satisfies the assumption in  Theorem \ref{corollary_6.4.4}.
By Theorem \ref{corollary_6.4.4} with $C_1=C_{11}+4+C,a=4$, we conclude that   $f \indexm$
   $C^{3,\alpha}$-converges to a
 function $f_\infty $ in $D_{a_1}(\fkz^\ast_\infty)$ for some constant $a_1>0$ independent of $k$.
 $\blacksquare$

\subsection{Proof of Theorem \ref{theorem_7.0.1}} \label{sect_7.5}

{
Put} $$
\mc W_f=\frac{W}{\max_{B_a(\fkz_o)}W},\;\;\;\mc R_f= \mc K(f) +\|\nabla\log |\mc S(f)|\|^2+\Psi(f).
$$
Suppose that the theorem is not true, then there is a sequence of functions
 $f\indexm$  and a sequence of points  $\fkz_k'\in B_{a/2}(\fkz_o)$  such that
  \begin{equation}\label{eqn_7.18}
   \mc W_k^{\half} \mc
  R_k(\fkz_k')a^2   \to \infty\;\;\; \mbox{ as } k\to \infty,
\end{equation}
where $\mc W_k:=\mc W_{f\indexm}, \mc R_k=\mc R_{f\indexm}$. Note that $\mc W_k\leq 1$  in $B_a(\fkz_o)$.
We apply the argument in Remark \ref{remark_4.1.1} to the function
$$F\indexm(z ):=  \mc W_k^{\half} \mc R_k(z)
 [ d_{f\indexm}(z,\partial
B_a(\fkz_o))]
 ^2$$
defined in $B_a(
 \fkz_o),$ the geodesic ball with respect to the metric $\omega_{f\indexm}$. Suppose that it
attains its maximum at $z^o\indexm$. By \eqref{eqn_7.18} we have
\begin{equation}\label{eqn_7.20}
\lim_{k\to\infty}F\indexm(z^o\indexm)\to +\infty, \;\;\;\;\lim_{k\to\infty}\mc W_k^{\half} \mc R_k(z^o\indexm)=+\infty, \;\;\;\lim_{k\to\infty} \mc R_k(z^o\indexm)=+\infty.
\end{equation}
Put $$d\indexm=\frac{1}{2} d_{f\indexm}(z^o\indexm,\partial B_a
(\fkz_o)).$$ Then 
in $B_{d_k}(z_k^o)$
$$
{W_k}^{\frac{1}{2}}\mc R_k \leq 4 {W_k}^{\frac{1}{2}}\mc R_k(z_k^o).
$$
Let $z_k^\ast\in Z$ be the point such that $ d(z_k^o,
z_k^\ast)=d(z_k^o,Z ).$ Denote by $q_k^o,q_k^\ast, \ldots$ the  images
of $ z_k^o, z_k^\ast, \ldots$ under the moment map $\tau_f$.

  Now we perform the affine
  blowing-up analysis to derive a contradiction.
We take an
affine transformation as in Lemma \ref{lemma_7.1.2}  by setting
   $$\alpha_k= \beta_k^2=\mc R_k(z_k^o),\;\;\;
   \eta_k= \log \alpha_k+\log W(z_k^o) ,$$
   and $b_k=c_k=\gamma_k=0$.
Then we assume that the notations are changed  as in $u_k\to \tilde u_k,
f_k\to\tilde f_k, d_k\to \tilde d_k$ and etc.
\v \def \half{\frac{1}{2}}
 \n{\em Claim 1:}
\begin{enumerate}
\item[(1)] $\alpha_k\to \infty$, $ \tilde d_k\to \infty$ as $k\to \infty$;
\item[(2)]$\lim\limits_{k\to\infty}\max\limits_{B_{\tilde d_k}(\tilde z^o_k)}
|\mc {\tilde S}_k| = 0$,
\item[(3)] $\tilde W(\tilde z^o_k)=1$;
{\item[(4)] ${\tilde W_k}^\half\tilde {\mc R}_k(\tilde z_k^o)=1$ and ${\tilde W_k}^\half\tilde {\mc R_k}\leq 4{\tilde W_k}^\half\tilde{\mc R}_k(z_k^o)=4$ in
$B_{\tilde d_k}(\tilde z^o_k)$;
\item[(5)] $\tilde W_k^{\half}\tilde\Psi_k\to 0$ in
$B_{\tilde d_k/2}(\tilde z^o_k)$.}
 \end{enumerate}
\v\n
{\bf Proof of claim.} (1-4) follows from Lemma \ref{lemma_7.1.2} and \eqref{eqn_7.20}. 
 Now we prove (5). From (4) and $\tilde S\neq 0$ we see that $\tilde f_k$ satisfies the assumption of Lemma \ref{lemma_7.16a}
 with $a=\tilde d_k$ and $o=\tilde z_k^o.$
 Set $$A_k:=\max_{B_{\tilde d_k}(\tilde z_k^o)} (\tilde W_k)^{\half}\max_{B_{\tilde d_k}(\tilde z_k^o)}
  |\tilde S_k|,  \;\;\;\;\;\;B_k=\max_{B_{\tilde d_k}(\tilde z_k^o)} (\tilde W_k)^{\frac{1}{4}}\tilde d_k\inv.$$
 Then by Lemma \ref{lemma_7.16a} in $B_{\tilde d_k/2}(\tilde z_k^o),$ we have
$$
\tilde W_k^{\half}\tilde\Psi_k \leq C\left[A_k+A_k^{\frac{2}{3}}+B_k+B_k^{2}\right].
$$
A direct calculation tells us
$$
A_k=\max_{B_{\tilde d_k}(\tilde z_k^o)}  \tilde W_k ^{\half}\max|\tilde S_k|\leq\left( \mc W_k\right)^{-\half}\mc R(z_k^o)\inv \max|S_k|
\to 0;
$$
$$
B_k^2=\max_{B_{\tilde d_k}(\tilde z_k^o)}  \tilde W_k ^{\half}\tilde d_k^{-2}
\leq  (\mc W_k^{\half}\mc R(\tilde z^o)d_k^2) \inv \to 0,
$$
 and so is $B_k$.  Hence
(5) follows. $\blacksquare$

\v
\begin{lemma}\label{lemma_modify_7.16}
$\tilde z_k^\ast \in B_{\sqrt{\mff C_5}}(\tilde z_k^o).$
\end{lemma}
{\bf Proof.} By Corollary \ref{corollaryc_4.2.2} and Corollary \ref{corollaryc_4.2.3} we have
$\mc R_k(z_k^o) d^2(z_k^o,Z)\leq \mff C_5.$
 Since $\tilde {\mc R}_k(\tilde z_k^o) \tilde {d}^2(\tilde z_k^o,Z)=\mc R_k(z_k^o) d^2(z_k^o,Z)$ and $\tilde {\mc R}_k(\tilde z_k^o)=1,$  we conclude that $\tilde d_{k}(\tilde z_k^o,Z)\leq \sqrt{\mff C_5}$ and $\tilde z_k^\ast \in B_{\sqrt{\mff C_5}}(\tilde z_k^o).$
$\blacksquare$
\v

By (5), we know that for any fixed $R$ and for any small constant $\epsilon>0$
\begin{equation}\label{new}
1-\epsilon \leq \tilde W_k(z)\leq 1+\epsilon,\;\;\;\;\;\tilde \Psi\indexm\leq \epsilon,\;\;\;\;\;  \forall z\in B_R(\tilde z^o_k),
\end{equation}
when $k$ large enough. By Lemma \ref{lemma_modify_7.16}, \eqref{new} also holds in $B_R(\tilde z^\ast_k).$ It follows from (4) and \eqref{new} that
\begin{equation}\label{eqn_tilde_R}
\tilde{\mc R}_{k}\leq 5, \;\;\;\;\;\;\;  \forall z\in B_R(\tilde z^o_k).
\end{equation}
To derive a contradiction we need the  convergence of $f_k.$ We discuss two cases.
\v\n{\bf Case 1.} There is a constant $C'>0$ such that
$\tilde d\indexm(\tilde z^o\indexm,\tilde z^\ast\indexm)\geq C'.$  By Lemma \ref{lemma_modify_7.16} we have
$$C'\leq \tilde d(\tilde z^o\indexm,\tilde z^\ast\indexm)=\tilde d(\tilde z^o\indexm,Z)\leq \sqrt{\mff C_5}.$$
We omit again the index $k$.  Let $\tilde p^\circ=\tau_{f}(\tilde z^\circ).$ To apply Theorem \ref{theorem_7.4.1}, we need the  following preparations:
\begin{itemize}
\item[(i)]  by affine transformations in \eqref{eqnc_7.1}
with $\lambda=\alpha=1,$ we can  minimal-normalize $(\tilde u, \tilde p^\circ, \tilde{\check p})$ (cf.\;Section \ref{sect_7.3}),
 \item[(ii)]  since  $\tilde{\mathcal R},\tilde{\mathcal K},\tilde\Psi$ and  $\tilde W(z)/\tilde W(z')$ for any $z,z'\in B_{N}(\tilde z^o)$  are invariant under these transformations,  the statements of {\em Claim 1} and \eqref{new}, \eqref{eqn_tilde_R} remain true.
\end{itemize}
After this transformation,  we assume that the notations are changed as  in $\tilde u\to \hat
u,\;\tilde d\to \hat
d$ and etc.   In this way we get a sequence $(\hat u\indexm, \hat p^\circ\indexm, \hat{\check p}\indexm)$.
 Applying  Theorem \ref{theorem_7.4.1} we
conclude that $\hat p^\circ\indexm$ converges to a
point $\hat p^\circ_\infty $ and $\hat u_k$
$C^{3,\alpha}$-converges to a  function $\hat u_\infty $ in a neighborhood of $\tau_f(\hat p^\circ_\infty);$ and $\hat z^*\indexm$ converges to a
point $\hat z^*_\infty $ and $\hat f_k$ locally $C^{3,\alpha}$-converges to a  function $\hat f_\infty$  in the ball
$D_{a_1}(\hat z^*_\infty)$. Then $\hat z^\circ\indexm$ converges to a
point $\hat z^\circ_\infty $ and $\tau_{f}(\hat z^\circ_\infty)=\hat p^{\circ}_{\infty}.$

\v\n{\bf Case 2.} $\lim\limits_{k\to\infty}\tilde d\indexm(\tilde z^o\indexm,\tilde z^\ast\indexm)=0.$
  We omit again the index $k$. By Lemma \ref{lemma_7.5.1} below, we can find a point
$\fkz^\circ$ in $B_2(\tilde z^\ast)$ such that $d(\fkz^\circ,Z)=c$.
For simplicity, without loss of generality, we assume that $c=1$ and
$d(\fkz^\circ,\tilde z^\ast)=1$.  Let $\tilde p^\circ=\tau_{f}(\fkz^\circ).$

We take a transformation as (i) and (ii) in {\bf Case 1} such that  $(\hat u\indexm, \hat p^\circ\indexm, \hat{\check p}\indexm)$ satisfies the conditions in Theorem \ref{theorem_7.4.1}.
Applying (2) of Theorem \ref{theorem_7.4.1} we
conclude that $\hat z^\ast\indexm$ converges to a
point $\hat z^\ast_\infty  $    and $\hat f_k$ locally $C^{3,\alpha}$-converges to a  function $\hat f_\infty$  in the ball
$D_{a_1}(\hat z^\ast_\infty)$.   By $\lim\limits_{k\to\infty}\hat d\indexm(\hat z^o\indexm,\hat z^\ast\indexm)=0,$ we have $\hat z^o\indexm$ converges to a
point $\hat z^\circ_\infty $ and $\hat f_k$ locally $C^{3,\alpha}$-converges to a  function $\hat f_\infty$  in the ball
$D_{a_1}(\hat z^\circ_\infty)$.   Then $ \hat z^\circ_\infty =\hat z^*_\infty.$

\v Hence for both cases we have $\hat z^o\indexm$ converges to a
point $\hat z^\circ_\infty $ and $\hat f_k$  $C^{3,\alpha}$-converges to a  function $\hat f_\infty$  in a neighborhood of $D_{b_1}(\hat z^\circ_\infty),$ and
\begin{equation}\label{eqn_7.29}
\hat W\equiv const.,\;\;\;\;\;\;C_1\inv\leq \hat
f_{i\bar j}\leq C_1,\;\;\;\;\;\;\mbox{ in }\;\;\;D_{b_1}(\hat z^\circ_\infty)
\end{equation}  where $C_1>0$ is a constant and $b_1$ is a constant independent of $k$. Here
\eqref{eqn_7.29}
follows from \eqref{new} and the $C^{3,\alpha}$-convergence of $\hat f_k$.
\v\n{\em Claim.} $\lim\limits_{k\to\infty}\max\limits_{D_{b_1}(\hat z^\circ_\infty)}\|\nabla \log |\mc S(\hat f_k)|\|^2_{\hat f_k}=0.$   
\v\n Proof of claim. Suppose the affine transformations from $u$ to $\hat u$ are
 ( we omit the index $k$)
 $$\hat \xi_1 =\hat \alpha \xi_1,\;\; \hat \xi_2=\hat \beta (\xi_2+\gamma),\;\; \hat u(\hat \xi_1, \hat \xi_2)=\hat \alpha u\left((\hat \alpha)\inv\hat \xi_1, (\hat \beta)\inv(\hat \xi_2-\gamma) \right)-l(\hat \xi), $$
where $  \gamma$ is a constant and $l(\hat \xi)$  is a linear function. Obviously,
\begin{equation}\label{eqna_7.30}
\lim\limits_{k\to\infty}\hat \alpha_k=\lim\limits_{k\to\infty}\alpha_k=+\infty.
\end{equation}
Next we prove
\begin{equation}\label{eqna_7.31}\lim_{k\to\infty} \hat\beta_k=+\infty.\end{equation}
To prove this we consider $\hat f\indexm$  in the Eucildean ball $D_{a_1}(\hat z^*_\infty).$ From the $C^{3,\alpha}$-convergence of $\hat f_k,$  we have
\begin{equation}
|\hat z_1|+|\hat w_2|\leq C_1',\;\;\;\; (C_1')\inv\leq \hat{f}_{i\bar j}\leq C_1',\;\;\;\;\;\;\mbox{  in }\;\;\;D_{a_1}(\hat z^*_\infty)
\end{equation} for some constant $C_1'>1.$ Then in $\t\cap D_{a_1}(\hat z^*_{\infty}),$ in terms of coordinates $\hat{x}_1, \hat{x}_2$,
\begin{equation}
C_2\inv \leq \frac{\partial^2 \hat f}{\p \hat x_2^2}\leq C_2,\;\; |\frac{\partial^2 \hat f}{\p \hat x_i\p \hat x_j}|\leq C_2.
\end{equation}
It follows that
\begin{equation}\label{eqnc_7.7}
C_2\inv \leq \frac{\partial^2 \hat u}{\partial {\hat \xi}_2^2} \leq C_2,\;\;\;|\hat u^{ij}|\leq C_2,\end{equation} for some constant $C_2>0 $ independent of $k$, where $(\hat u^{ij})=(\hat u_{ij})^{-1}$.
By   $$  \frac{\partial^2 \hat u}{\partial {\hat \xi}_2^2} =\frac{\hat \alpha}{\hat \beta^2}\frac{\partial^2 u}{\partial {\xi}_2^2},
\mbox{  and  }\frac{\partial^2 h}{\partial \xi_2^2}|_{\ell\cap \tau_f(B_a(\fkz_o))} \geq \mff N_5\inv$$  we conclude that
$$\frac{\hat{\beta}^2}{\hat{\alpha}}\geq \mff N_5\inv C_2\inv.$$ Then \eqref{eqna_7.31} follows from \eqref{eqna_7.30}. Note that \eqref{eqnc_7.7} also holds in $D_{b_1}(z^\circ_{\infty}).$ Since
\begin{equation*}\left\|\nabla \log |\mc S(\hat f)|\right\|^2_{\hat f}=\left\|\nabla \log |\mc S( f)|\right\|_{\hat f}^2=\sum {\hat u}^{ij}\frac{\partial \xi_k}{\partial \hat \xi_i} \frac{\partial \xi_l}{\partial \hat \xi_j} \frac{\partial \log |\mc S(f)|}{\partial \xi_k}    \frac{\partial \log |\mc S(f)|}{\partial \xi_l}, \end{equation*}
and
$$\frac{\partial \xi_1}{\partial \hat \xi_1}=\frac{1}{\hat \alpha},\;\;\frac{\partial \xi_2}{\partial \hat \xi_2}=\frac{1}{\hat\beta},\;\;\frac{\partial \xi_1}{\partial \hat \xi_2}=\frac{\partial \xi_2}{\partial \hat \xi_1}=0,$$
by    \eqref{eqna_7.30}, \eqref{eqna_7.31} and \eqref{eqnc_7.7} the claim is proved.
\v\n

Combining this claim  and
$\mc {\tilde R}_k(\tilde z_k^o)=\mc {\hat R}_k(\hat z_k^o)  $ we get, for $k$ large enough,
\begin{equation}
\label{eqna_7.33}\mc {\hat K}_k(\hat z_k^o)\geq 1-\epsilon .\end{equation}
Note that
$$\frac{\partial^{i+j}  \mc S(\hat u)}{\partial \hat \xi_1^i\partial \hat \xi_2^j}=\left(\hat{\alpha}\right)^{-i-1}\left(\hat{\beta}\right)^{-j}\frac{\partial^{i+j}  \mc S(u  )}{\partial   \xi_1^i\partial \xi_2^j}.$$
Then by \eqref{eqna_7.30}, \eqref{eqna_7.31} and $\|\mc S(u_k)\|_{C^3(\Delta)}\leq \mff N_5$ we have
$$\lim_{k\to\infty}\| \mc S(\hat u_k)\|_{C^{3}}=0,$$
where $\|.\|_{C^3}$ denotes the Euclidean  $C^3$-norm in $(\hat \xi_1,\hat\xi_2).$
 Then,  by \eqref{eqn_7.29} and Theorem \ref{corollary_6.4.4}, $\hat f\indexm$ $C^{6,\alpha}$-converges to $\hat f_\infty$ in a neighborhood of $\hat z^o_\infty$.
 Then by \eqref{eqn_7.29}
we have $ \hat{\mc K}\equiv 0. $   This contradicts
\eqref{eqna_7.33}.
The theorem is proved. $\;\;\;\blacksquare$

Now we prove the following lemma that is needed in the above proof (cf.\;the paragraph right after \eqref{new}).
\begin{lemma}\label{lemma_7.5.1}
Let $z^\ast\in Z$. Suppose that $ \mc K\leq 4$ in $B_2(z^\ast)$. Then there exists a constant $c>0$, independent of $f$,  such that there exists a point $z^o$ in $\partial B_1(z^\ast)$ satisfying
$$d(z^o, B_2(z^\ast)\cap Z)=c.$$
Obviously $c\leq 1$. Hence $d(z^o,Z)=c$.  Here $d(\cdot,\cdot)$ denotes  the Riemann distance.
\end{lemma}
{\bf Proof.} Without loss of generality, we assume that the
$w$-coordinate of $z^\ast$ is
$$
(x_1,y_1)=(-\infty, 0),\;\;\; ( x_2,y_2)=(0,0).
$$
Recall that $\t$ is identified with $\t\times 2\sqrt{-1}\{1\}$ of
$\CHART_{\halfplane}$. Then when we consider $B_2(z^\ast)$, we
restrict ourself to  $B_2(z^\ast)\cap\t$. Similarly, when we consider
$Z$ we treat it as the line $\t_2=\{x_1=-\infty\}$, which is the
dual to $\t^\ast_2$ (when treat $Z$ as a 1-dimensional toric
manifold).

If the claim of the lemma does not hold, then there exists a
sequence of functions $f\indexm$ such that:
 for any point $z\in \partial B_1\indexn(z^\ast)$
$$
d_{f_k}(z, B_2\indexn(z^\ast)\cap Z)\leq \frac{1}{k}.
$$
 Let $\hat z\indexm\in B_2\indexn(z^\ast)\cap Z$ such that
 $$
 d_{f\indexm}(z,\hat z\indexm)=d_{f\indexm}(z,B_2\indexn(z^\ast)\cap
 Z)\leq \frac{1}{k}.
 $$
 Then
 $$
 1-\frac{1}{k}\leq d_{f\indexm}(z^\ast,\hat z\indexm)\leq
 1+\frac{1}{k}.
 $$
Let $A\indexm= B_1\indexn(z^\ast)\cap \t_2$ and
 split it into $A^\pm\indexm$ according to the sign of the $x_2$-coordinate.
Define
$$
\mc A\indexm^\pm=\{z\in \partial B_1(z^\ast)| d_{f_k}(z, A\indexm^\pm)\leq 1/k\}.
$$
Then $\mc A\indexm^\pm$ are non-empty and are closed subsets of
$\partial B_1\indexn(z^\ast)$. Since they cover $\partial B_1\indexn(z^\ast)$,
their intersection is non-empty. Choose a point $z\indexm^\circ\in \mc
A\indexm^+\cap \mc A\indexm^-$. Choose two points $z\indexm^\pm\in A\indexm^\pm$ such that
$d(z\indexm^\circ,z\indexm^\pm)\leq 1/k$. We normalize the function $f\indexm$ by adding a linear function (cf.\;Section \ref{sect_7.1}) such
that it achieves its minimum at $z\indexm^\pm$. We denote the function
by $f\indexm^\pm$ respectively.

Then we have following facts: \v {\em Fact 1:
$|f\indexm^\pm(z^+)-f\indexm^\pm(z^-)|\leq Ck\inv$ for some constant $C$.} In
fact, by \eqref{eqn_6.10} and the assumption, we have
$$|\log (1+f\indexm^\pm)(z\indexm^+)-\log (1+f\indexm^\pm)(z\indexm^-)|\leq \sqrt{\mff C_{10}}d(z\indexm^+,z\indexm^-).$$
Then
$$
|f\indexm^\pm(z\indexm^+)-f\indexm^\pm(z\indexm^-)|\leq C'd(z\indexm^+,z\indexm^-)\leq Ck\inv.
$$
\v {\em Fact 2: let $\bar f\indexm^\pm=f\indexm^\pm|_{Z}$ and $\bar d$ denote the geodesic  distance on $\t_2$ with
respect to $G_{\bar f\indexm^\pm}$. Then $\bar d(z\indexm^+,z\indexm^-)\geq 1$
when $k$ large.} In fact,
$$
\bar d(z\indexm^\pm, z^\ast)\geq d(z\indexm^\pm, z^\ast)\geq 1-1/k.
$$

We omit the index $k$ for simplicity. Now we focus on $\bar f^\pm$. By changing coordinate on $x_2$ we
may assume that
$$
x_2(z^+)-x_2(z^-)=1.
$$
Set $x_2^\pm=x_2(z^\pm)$. Then
$$
1 \leq \bar d^2 ( z^+,z^-) \leq
\left(\int_{x_2^-}^{x_2^+} \sqrt{\bar f^\pm_{22}} dx_2\right)^2
\leq\int_{x_2^-}^{x_2^+} {\bar f^\pm_{22}} dx_2 \leq   |{\bar
f^\pm}_2(z^+)-\bar f^\pm_2(z^-)|.
$$
We summarize that we have two convex functions $\bar f^\pm$  on the {\em unit} interval
$[x^-_2,x^+_2]$ that
are different up to a linear function and they have the following  properties:
\begin{itemize}
\item $ |\nabla\bar f^+(x_2^+)-\nabla \bar f^+(x_2^-)|\geq 1
$, $ |\nabla\bar f^-(x_2^+)-\nabla \bar f^-(x_2^-)|\geq 1;
$ \item $ |\bar f^+(x_2^+)-\bar f^+(x_2^-)|\to 0$,\;\;$ |\bar f^-(x_2^+)-\bar f^-(x_2^-)|\to 0, $ as\;$ k\to \infty$;
\item $\nabla\bar f^+(x_2^+)=0$,$\nabla\bar f^-(x_2^-)=0$.
\end{itemize}
This is impossible when $k$ is large enough. We get a contradiction. $\blacksquare$

\end{document}